\documentclass{amsart}
\usepackage{bbm}
\usepackage{graphicx}
\usepackage{verbatim}
\usepackage{stackrel}
\usepackage{color}
\graphicspath{{Fig/}}
\usepackage{amsmath}
\usepackage{amssymb}
\usepackage{array}
\usepackage{mathtools}
\usepackage{xcolor}
\usepackage{mathdots}
\usepackage{pgf}
\usepackage{bm}
\usepackage{hyperref}

\usepackage{enumitem}

\usepackage{tikz}
\usetikzlibrary{shapes.arrows, fadings}
\usetikzlibrary{arrows}
\usetikzlibrary{matrix}
\usetikzlibrary{shapes,intersections}

\usepackage[latin1]{inputenc}

\hypersetup{
    colorlinks,
    linkcolor={blue!40!black},
    citecolor={pink!40!red},
}

\newtheoremstyle%
 {bluethm}%
 {}{}%
 {\color{blue}\itshape}
 {}%
 {\color{blue}\bfseries}%
 {\color{blue}.}%
 { }{}

\newtheorem{theorem}{Theorem}
\newtheorem{prop}{Proposition}
\newtheorem{lemma}{Lemma}
\newtheorem{corollary}{Corollary}

\theoremstyle{definition}

\newtheorem{defi}{Definition}

\newtheorem{remark}{Remark}

\newtheorem{assumption}{Assumption}[part]

\def\N{{\mathbb N}}
\def\R{{\mathbb R}}

\def\P{{\mathbb P}}
\def\E{{\mathbb E}}

\newcommand{\diff}{\mathop{}\mathopen{}\mathrm{d}}

\newcommand\ind[1]{\mathbbm{1}_{\left\{#1\right\}}}

\newcommand{\bfx}{\mathbf{x}}
\newcommand{\bfy}{\mathbf{y}}
\newcommand{\bfz}{\mathbf{z}}

\newcommand{\tridiff}{\diff s\diff\bfx\diff\bfy}
\def\cal{\mathcal}
\def\eps{\varepsilon}

\makeatletter
\newcommand{\proofstep}[1]{%
  \par
  \addvspace{\medskipamount}
  \textit{#1\@addpunct{.}}\enspace\ignorespaces
}
\makeatother

\title[LDP for the pure jump $k$-nary interacting particle systems]{Pathwise Large deviations for the pure jump $k$-nary interacting particle systems}

\author{Wen Sun}\thanks{The  author's work is supported by the  Dirichlet Postdoc Fellowship funded by the Berlin Mathematical School.}
\email{Wen.Sun@TU-Berlin.de}
\address        {TU Berlin Mathematics Building, Str. des 17. Juni 136, 10587 Berlin}

\date{\today}

\begin{document}

\begin{abstract}
  A pathwise large deviation result is proved for the pure jump models of $k$-nary interacting particle system introduced by Kolokoltsov~\cite{kololln,kolokoltsov2010nonlinear} that generalize classical Boltzmann's collision model, Smoluchovski's coagulation model and many others. The upper bound is obtained by following the standard methods (KOV~\cite{kipnis1989hydrodynamics}) of  using a process ``perturbed'' by a regular function. To show the lower bound, we propose a family of orthogonal martingale measures and prove a coupling for the general perturbations. The rate function is studied based on the idea of L\'eonard~\cite{leonard1995large} with a simplification by considering the conjugation of integral functionals on a subspace of $L^{\infty}$. General ``gelling'' solutions in the domain of the rate function are also discussed.
\end{abstract}

\maketitle

\bigskip

\hrule

\vspace{-3mm}

\tableofcontents

\vspace{-1cm}

\hrule

\bigskip

\section{Introduction}
We study the large deviation problem for the pure jump models of $k$-nary interacting particle system proposed by Kolokoltsov~\cite{kololln,kolokoltsov2010nonlinear}.

\subsection{The pure jump $k$-nary interacting particle system}

We first recall that a configuration of the finite  $k$-nary interacting particle systems is of the set form
$ {\bfx}=\{x_1,x_2,\dots,x_{n}\}\subset X$
where  $n$ is the number of particles, $x_i$ is the state of the particle taking values in the state space $X$. The multinary interaction describe for example  the  chemical reactions  among any  subset of particles, \emph{i.e.}, $\bfx_I=\{x_{i_1},x_{i_2},\dots,x_{i_\ell}\}$  of $\bfx$. The products of this reaction is a set of particles at states $\bfy=\{y_1,y_2,\dots, y_m\}$. In this case, the state of the system jumps from the configuration $\bfx$ to the configuration $(\bfx\setminus \bfx_I)\cup\bfy$.
When the number of interacting particles $\ell$ and the number of product particles $m$ are not larger than  a constant $k$, we call it a $k$-nary interacting particle system. The transition of such a reaction $\bfx_I\mapsto \bfy$ is captured by a time-inhomogeneous transition kernel $P(t,\bfx_I;\diff \bfy)$.

The most classical binary interacting particle system is the \emph{spatially homogeneous Boltzmann collision model} where $X=\R^d$, the  space of velocities. The binary collision occurs among any two particles in a closed system with velocities changing from $(z,z')$ to $(z^*,z'^*)$ that preserve total momentum and energy,
\[
z+z'=z^*+z'^*,\qquad |z|^2+|z'|^2=|z^*|^2+|z'^*|^2.
\]
with a collision kernel $B(z,z',\diff z^*,\diff z'^*)$ that is a symmetric function on the product state space $(\R^d)^2$ and the probability space of the product space $\cal{P}((\R^d)^2)$. In our framework, the kernel $P$ then can be represented as
\[
P(\bfx;\diff\bfy)=B(z,z',\diff z^*,\diff z'^*)\qquad \textrm{if }\bfx=\{z,z'\},\bfy=\{z^*,z'^*\}.
\]

A typical example with time changing number of particles is the mass-preserving \emph{Smoluchovski coagulation model} where $X=\R^+$, the space of  mass. Each two particles with mass $\{x_{i_1},x_{i_2}\}$ can coagulate into a larger particle with mass $y=x_{i_1}+x_{i_2}$ (binary coagulation) at rate $K(x_{i_1},x_{i_2})$. In this case, the kernel $P$ writes
\[P(\bfx;\diff \bfy)=K(x_{i_1},x_{i_2})\delta_{\{x_{i_1}+x_{i_2}\}}(\diff \bfy)\qquad\textrm{if }\bfx=\{x_{i_1},x_{i_2}\}.\]

The  pure jump $k$-nary interacting particle system also generalizes the reversible and mass-preserving \emph{Becker-D\"oring coagulation and fragmentation model} where $k=2$, $X=\N^*$ and the $x_i$ represent the mass of the particles. Each particle with mass $x$ may interact with a monomer (particles with mass $1$) merging into a particle with mass $x+1$ at rate $a_x$; each particles with a mass larger than $2$ may break into a monomer and a particle with mass $x-1$ at rate $b_x$. Then the kernel $P$ writes 
\begin{align*}
&\textrm{(coagulation)} \hfill &P(\bfx;\diff \bfy)=a_x\delta_{\{x+1\}}(\diff \bfy),& \qquad \textrm{if }\bfx=\{x,1\};&\\
&\textrm{(fragmentation)}  \hfill &P(\bfx;\diff \bfy)=b_x \delta_{\{1,x-1\}}(\diff \bfy),& \qquad \textrm{if }\bfx=\{x\}.&
\end{align*}

\subsection{The hydrodynamic equation}

Let $\bfx(t)=\{x_1(t),\dots,x_{n(t)}(t)\}$ be the jump process described above.
It is proved by Kolokoltsov~\cite{kololln}, under certain conditions and a proper scaling $h$, for any finite time interval, the empirical process
$$t\in[0,T]\mapsto h\delta^+_{\bfx(t)},$$
where $\delta^+$ maps any finite subset $\bfz=\{z_1,z_2\dots,z_j\}\subset X$ to $\N$ by
\[
\delta^+_\bfz=\sum_{i=1}^j\delta_{z_i},
\]
satisfies a law of large numbers as $h\to 0$ and that the limiting deterministic dynamics is characterized as a weak solution of the nonlinear equation
  \begin{equation}\label{llneq}
\diff \langle g,\sigma_t\rangle=\int_{\cal{J}}\langle g,\delta^+_{\bfy}-\delta^+_\bfx\rangle \cal{K}^{ P}[\sigma_t](\diff t\diff \bfx\diff\bfy),\qquad g\in\cal{C}_b(X)
  \end{equation}
  where \begin{itemize}
  \item $\cal{C}_b(X)$ denotes the set of bounded continuous functions on $X$;
  \item $\cal{J}$ denotes the state space of all jump pairs $(\bfx,\bfy)$ (specified later in Definition~\ref{jumpspace});
  \item $\cal{K}^P[\cdot]$ denotes a  mapping,  depending on transition kernel $P$, from  a finite measure on $X$ to a finite measure on the  trajectory space of jumps $[0,T]{\times}\cal{J}$ (specified later in Definition~\ref{jumptraj}).
  \end{itemize}

    For instance,  the Smoluchovski coagulation equation writes
    \[
\diff \langle g,\sigma_t\rangle=\frac{1}{2}\int_{\R^+\times \R^+} \langle g,\delta_{x_1+x_2}-\delta_{x_1}-\delta_{x_2}\rangle K(x_1,x_2)\sigma_t(\diff x_1)\sigma_t(\diff x_2)\diff t
\]
where $\frac{1}{2}K(x_1,x_2)\sigma_t(\diff x_1)\sigma_t(\diff x_2)$ can be seen as a measure on the (pair) jump space
    \[\cal{J}=\left\{\big(\{x_1,x_2\},\{x_1+x_2\}\big)\bigg|x_1,x_2\in\R^+\right\}.\]

    \subsection{Main Contributions}
    We obtain a large deviation estimates for the path of  empirical measures $h\delta^+_{\bfx(t)}$ passing to the the hydrodynamic equation~\eqref{llneq} as $h\to 0$. By using  the standard ``good'' perturbation -- the  particle system with $f$-perturbed kernel
    \[
e^{\langle f_t,\delta^+_\bfy-\delta^+_{\bfx}\rangle}P(t,\bfx;\diff\bfy),
\]
for some class of  regular functions $f$, we obtain a rate function for the upper bound in a variational form
\begin{multline*}
  \cal{R}^{P}_{\rm upper}(\pi)=\sup_{f} \bigg\{\langle \pi_T,f_T \rangle-\langle \pi_0,f_0\rangle-\int_0^T\langle \pi_t,\partial_t f_t\rangle \diff t\\
-\int_0^T\int_{\cal{J}}\left(e^{\langle f_t,\delta^+_\bfy-\delta^+_{\bfx}\rangle}-1\right)\cal{K}^P[\pi](\diff t\diff\bfx\diff\bfy)\bigg\}.
\end{multline*}
An alternative non-variational formula for this rate function,
\[
\cal{R}^{P}_{\rm upper}(\pi)=  \inf_{\eta \ge 0} \int_0^T\int_{\cal{J}}\left(\eta\log(\eta)-\eta+1\right)(t,\bfx,\bfy)\cal{K}^P[\pi](\diff t\diff\bfx\diff\bfy),
\]
where $\eta $ running over all non-negative measurable function on $[0,T]{\times}\cal{J}$ such that $\pi$ satisfies the  hydrodynamic equation~\eqref{llneq} with perturbed kernel
$$\eta (t,\bfx,\bfy)P(t,\bfx;\diff\bfy),$$
is given in the papers  L\'eonard~\cite{leonard1995large,leonard2001convex}  by considering the duality on two Orlicz spaces. We obtain the same formula with a simpler proof by considering the conjugation of integral functionals on a subspace of $\cal{L}^\infty$. Our proof is based on the paper Rockafellar~\cite{rock1}. We claim that in general, the perturbed kernel is not unique determined by the path $\pi$ due to the lack of uniqueness in the general Hahn-Banach Theorem.

Hence unfortunately, we may not be able to expect an approximation from the $f$-perturbed kernel to this general perturbation that is valid in many $1$-nary interacting models, for example, in Kipnis~et~al.~\cite{kipnis1989hydrodynamics} and Jona-Lasinio et al.~\cite{jona1993}. In order to find a perturbed particle system  passing to the path in the domain of the rate function $\cal{R}^{P}_{\rm upper}(\pi)$, we use the idea of martingale measures in the sense of Walsh~\cite{walsh} to keep track of all jumps of the measure valued process $h\delta^+_{\bfx(\cdot)}$. In this frame work, inspired by Shiga and Tanaka~\cite{Shiga1985CentralLT}, we are able to give a coupling between the original particle system with kernel $P$ and the perturbed particle system with kernel $\eta P$ when the Radon-Nikodym derivative $\eta$ is bounded. Based on this coupling, we obtain the rate function of lower bound for the large deviation on a class of  ``good'' paths  given by
\[
\cal{R}^{P}_{\rm lower}(\pi)=  \inf_{\eta \ge 0} \int_0^T\int_{\cal{J}}\left(\eta\log(\eta)-\eta+1\right)(t,\bfx,\bfy)\cal{K}^P[\pi](\diff t\diff\bfx\diff\bfy)
\]
where $\eta $ running over all non-negative \emph{bounded} measurable function on $[0,T]{\times}\cal{J}$ such that $\pi$ satisfies the  hydrodynamic equation~\eqref{llneq} with perturbed kernel $\eta P.$ Naturally, we have the inequality $$\cal{R}^{P}_{\rm lower}(\pi)\ge \cal{R}^{P}_{\rm upper}(\pi).$$ If the Radon-Nikodym derivative $\eta$ is bounded and uniquely determined by the path $\pi$, then $$\cal{R}^{P}_{\rm lower}(\pi)=\cal{R}^{P}_{\rm upper}(\pi).$$ We shall see this is usually not true in the example of Becker-D\"oring model (Section~\ref{exbd}). The uniqueness of $\eta$ highly depends on the structure of the jump space $\cal{J}$, we offer a sufficient condition in this paper.

In general, the well-posedness of equation~\eqref{llneq} with kernel $\eta P$ in a shape of the path in the domain of the rate function $\cal{R}^{P}_{\rm upper}(\pi)$ is not guaranteed. For example, if the coagulation kernel $K(x,y)$ is too large, a phenomenon called ``gelation'' may occur in the Smoluchovski coagulation and cause problems related to the existence and uniqueness properties. The stochastic Smoluchovski coagulation processes, also called \emph{Marcus-Lushnikov processes}, may lead to a limit that has a correction term in addition to the  Smoluchovski equation. Nevertheless, some ``gelling'' solution lies in the domain of our rate function and can be approximated by the solution with truncated Radon-Nikodym derivative. We provide a existence results on these type of ``gelling'' solution by extending the results of the Smoluchovski equation by Norris~\cite{norris1999} to our $k$-nary interacting equation. For more about the stochastic approximation to Smoluchovski equation, we refer to the papers Jeon~\cite{jeon1998existence}, Norris~\cite{norris1999,norris2000}, Fournier et al.~\cite{FOURNIER2009}, Kolokoltsov~\cite{kolokoltsov2010central} and also the survey Aldous~\cite{Aldous1999}.

The main theorems are stated in Section~\ref{sec:main}. We also apply our results to
\begin{itemize}
\item the Becker-D\"oring  coagulation and fragmentation process with reaction rates that are  strictly sublinear with respect to the mass;
\item the Boltzmann collision model with a collision kernel that is asymptotically upper controlled by the energy function,
\end{itemize}
in the end of Section~\ref{sec:main}. For the large deviation of Boltzmann collision models, we refer to L\'eonard~\cite{leonard1995large} and two very recent works by Basile et al.~\cite{basile2021large} and Heydecker~\cite{heydecker2021large}. In a forthcoming work~\cite{ldpgel}, we investigate the pathwise large deviation of  Smoluchovski's coagulation model with gelling kernels in a similary framework.

\subsection*{Outline of the paper}
Section~\ref{sec:process} introduces the basic notations and the measure-valued Markov process describing the evolution  of the $k$-nary interacting particle system.
Section~\ref{sec:main} presents our main results on the pathwise large deviation, rate functions and the corresponding  ``good'' paths. The ``good'' gelling solution  is also discussed.
Section~\ref{sec:coupling} gives an reformation of the measure-valued Markov process by using a family of martingale measures on the space of jumps. In this framework, a general coupling between two particle systems with different kernels is established that is essential in the proof of the lower bound for the large deviation.
Section~\ref{sec:lln} provides estimations of the measures on the path space of jumps and on the moments of the Markov process. The proofs of the law of large numbers (Theorem~\ref{lln}), the limit of the martingale measures (Corollary~\ref{mb}) and the limit of the relative entropy (Corollary~\ref{appsub}) are also given.
Section~\ref{sec:ldp} presents the proof of large deviation (Theorem~\ref{main:ldp}) and the study of rate functions (Theorem~\ref{ratethm}).
Section~\ref{sec:cm} gives the proof of Theorem~\ref{gel}, the existence of the path before ``gelation''.
Appendix~\ref{sec:app} gives a quick overview and an extension of the conjugations of integral functionals on $\cal{L}^\infty$.

\section{Measure-valued Markov process}\label{sec:process}
In this section we first introduce the notations and assumptions used throughout this paper. Then we define a measure-valued Markov process that describing the evolution of the $k$-nary interacting particle system.

\subsection{Basic Notations}

The state space of one particle is a Polish space $X$ endowed with its Borel sigma algebra  $\cal{B}(X)$. For each $j\in\N$, let  $X^j=X{\times}X{\times}\dots{\times}X$ and  $\cal{X}:=\cup_{j=0}^\infty X^j$ equipped with their product topologies.
The state space of the system is thus given by $\cal{X}$. If $\vec{\bfx}=(x_1,\dots,x_j)\in \cal{X}$, then the system has $j$ particles at states $x_i$, $1\le i\le j$.

Since the interaction is assumed to be exchangeable, for the later convenience, we introduce the space $S\cal{X}=\cal{X}\setminus \sim$ (\emph{resp.} $S{X}^j={X}^j\setminus \sim$) be the quotient space with a equivalence relation $\sim$  defined by all permutations. For example $(x_1,x_2,\dots,x_n)\sim (x_{p(1)},x_{p(2)},\dots,x_{p(n)})$ for any $p$ who is a permutation on $\{1,2,\dots,n\}$.
The elements in $\cal{X}$ will be denoted by vectors $\vec{\bfx},\vec{\bfy}$ \emph{etc.} The corresponding elements in $S\cal{X}$ will be denoted by $\bfx,\bfy$ \emph{etc.} Let $|\bfx| (\emph{resp.}  |\vec{\bfx}|$) be the number of elements (\emph{resp.} length) of $\bfx (\emph{resp.}  \vec{\bfx}$).

Clearly, the empirical measure of the state is given by a mapping  $\delta^+_\cdot:S\cal{X}\mapsto \cal{M}^+_\delta(X)$ 
\[
\mathbf{x}\sim(x_1,x_2,\dots,x_n)\mapsto \delta^+_{\mathbf{x}}=\delta_{x_1}+\delta_{x_2}+\dots+\delta_{x_n}
\]
where $\delta_x$ for all $x\in X$ is the Dirac's measure and
\[
\cal{M}_\delta^+(X)=
\left\{\delta^+_{\mathbf{x}}\big|\mathbf{x}\in S \cal{X}\right\}.
\]
is the space of finite sum of $\delta$-measures on $X$. Clearly, $\delta^+$ defines a bijection between $S\cal{X}$ and $\cal{M}^+_\delta(X)$.

Throughout this paper, we fixed the time $T>0$. Let $\cal{M}^+(X)$ be the set of non-negative finite measures on $X$ endowed with weak topology. 
Our particle system can be described by the c\`adl\`ag stochastic process
\[
\delta^+_\bfx:[0,T] \mapsto \delta^+_{\bfx(t)}\in D_T(\cal{M}^+(X)),
\]
where $D_T(\cal{M}^+(X))$ is the space of all c\`adl\`ag paths from $[0,T]$ to $\cal{M}^+_\delta(X)$ endowed with the Skorokhod $(J_1)$ topology.

In order to study the measure-valued path, we also need the following notations.
Let $\cal{C}_b(X)$ (\emph{resp.} $\cal{C}_b(S\cal{X})$) denote the space of all bounded continuous functions on $X$  (\emph{resp.} $S\cal{X}$) .
Let $\cal{C}_{\rm sym}(\cal{X})$ (\emph{resp.} $\cal{C}_{\rm sym}(X^j)$) be the Banach space of symmetric bounded continuous functions on $\cal{X}$ (\emph{resp.} on $X^j$). Clearly $\cal{C}_{\rm sym}(\cal{X})=\cal{C}(S\cal{X})$ and $\cal{C}_{\rm sym}({X}^j)=\cal{C}(S{X}^j)$.  We denote by $\cal{C}^{1,0}_b([0,T]{\times} X)$  the bounded continuous functions on $[0,T]{\times}X$ which have first continuous derivatives with respect to the time variable. Let $\cal{M}_{\rm sym}^+(\cal{X})$ (\emph{resp.} $\cal{M}_{\rm sym}^+(X^j)$) be the space of symmetric non-negative finite measures on $\cal{X}$ (\emph{resp.} on $X^j$). With some abuse of notation, for space $V=X,\cal{X},S\cal{X}, X^j,SX^j$, we use the pair $\langle f,\mu \rangle$ represents the integral $\int f\diff \mu$ for a measurable function $f$ on $V$ and a measure $\mu$ on $V$.

We can define a metric on $\cal{M}^+(X)$ by
\[
d(\mu,\nu):=\sum_{m=0}^\infty\frac{1}{4^m}\left(1\land \left|\langle f_m,\mu\rangle-\langle f_m,\nu\rangle\right|\right)
\]
where $\{f_m\}$ is  a sequence of non-negative bounded functions on $X$ that are convergence determining and
$\langle f,\mu \rangle$ stands for the integral $\int_{X}f(x)\mu(\diff x)$. Then $d(\mu,\nu)$ metrizes the weak topology. See Section~3.2 in Dawson~\cite{dawson}.

In order to study the measures on $S\cal{X}$, for any $\ell\in\N$, for all $\pi\in\cal{M}^+(X),$ we introduce a measure $\pi^{\tilde{\otimes} \ell}$ on $SX^\ell$ by, for all measurable set $U\in \cal{B}({SX^\ell})$,
\[
\pi^{\tilde{\otimes} \ell}(U):=\frac{1}{\ell!}\int_{\overline{U}}\prod_{j=1}^\ell \pi(\diff y_j)
\]
where
\[
\overline{U}:=\left\{\vec{\bfy}\in X^\ell\big|{\bfy}\in U \right\}.
\]
For all Borel measurable function $f$ on $X$, we define a function $f^{\tilde{\otimes} \ell}$ on $SX^\ell$: for any $\bfx\in SX^\ell$
\[
f^{\tilde{\otimes} \ell}(\bfx):=\prod_{j=1}^\ell f(x_j)
\]
by picking any  representation $\vec{\bfx}=(x_1,\dots,x_\ell)\sim \bfx$.

To control the  growth behavior of the particles at very large states, we introduce a function $E$, a non-negative continuous function from $X$ to $\R^+$ with compact level sets, \emph{i.e.}, for all $a>0$, $\{x|E(x)\le a\}$ are compacts in $X$. Let $\cal{M}_{1+E}$ be the set of all finite non-negative measures on $X$ that integrate $1+E$:
  \[
  \cal{M}_{1+E}(X):=
  \left\{\pi\in\cal{M}^+(X)\big|\langle 1+E,\pi\rangle<\infty\right\}.
  \]
Let $\cal{L}^\infty_{1+E}(X)$ denotes the set of all  real-valued measurable functions $f$ on $X$ with the property that
  \[
\limsup_{E(x)\to \infty}\frac{|f(x)|}{1+E(x)}<\infty.
\]
Equip this space $\cal{L}^\infty_{1+E}(X)$ with the norm
\[
\|f\|_{(1+E),\infty}:=\sup_{x\in X}\frac{f(x)}{1+E(x)}.
\]
In particular, for all $f\in \cal{L}^\infty_{1+E}(X)$
\[
|\langle \mu,f \rangle|\le \|f\|_{(1+E),\infty}\langle \mu,1+E\rangle<\infty.
\]
For $S=\cal{M}^+(X),\cal{M}_{1+E}(X)$, let $C_T(S)$ denote the space of $S$-valued continuous functions on $[0,T]$; $D_T(S)$ denotes the space of $S$-valued c\`adl\`ag functions on $[0,T]$  endowed with the Skorokhod $(J_1)$ topology.

\subsection{The dynamic of the $k$-nary interacting particles}

To describe the dynamic of the interaction among a set of particles $\bfx$ with a set of products $\bfy$, we introduce the following time in-homogeneous transition kernel.

\begin{assumption}[Kernel: basic assumptions]\label{kernel}
  We denote by $P(t, \bfx,\diff \bfy)$ the $k$-nary interaction kernel which is a measurable function
  \[      
      [0,T]\times S{\cal{X}}~ \mapsto \cal{M}^+(S \cal{X}),
  \]
  where $ \cal{M}^+(S \cal{X})$ equipped with its weak topology. Recall that $E$ is a non-negative continuous function from $X$ to $\R^+$ with compact level sets.
  
  Throughout out the paper, we assume 
  \begin{enumerate}
  \item The number of particles that participating one interaction is not larger than $k$.
  \item The increment of number of particles (\emph{fragmentation}) only happens in the $1$-nary interaction and this increment is not larger than $k-1$, \emph{i.e.} for all $1\le \ell \le k$, let
    \[
d_\ell:= \left\{ \begin{array}{ll}
k,& \ell=1;\\
\ell, & 2\le\ell\le k,
\end{array} \right.
    \]
and for   any $t\in[0,T]$, any $\bfx\in SX^\ell$,
    \[
{\rm supp}P(t,\bfx,\cdot)\subset \cup_{m=0}^{d_\ell} SX^m;
\]

\item For all $t\in[0,T]$, $1\le \ell\le k$, $0\le m\le d_\ell $,  let
  \[
P_m^\ell(t,\mathbf{x},\diff \mathbf{y}):=P(t,\bfx,\diff \bfy)\ind{\bfx\in SX^\ell,\bfy\in SX^m}.
\]
They are \emph{$E$-non-increasing}, that is,  for all a.e.-$t\in[0,T]$ and all  $\bfx\in SX^\ell$,
\begin{equation*}
\langle E,\delta^+_\bfy-\delta^+_\bfx\rangle \le 0\qquad P^{\ell}_m(t,\bfx,\diff\bfy)-a.s.
\end{equation*}
\item We assume $P$ satisfies a ``No dust condition'':  the kernel $P$ is either $1$-non-increasing (no fragmentation in the $1$-nary interaction), or there exists a constant $\eps_0>0$, such that $E(x)\ge \eps_0$ for all $x\in X$.

\item   For all $1\le \ell\le k$, for all $\mathbf{x}\in SX^\ell$,  let
  \[
  P(t,\mathbf{x}):=
  \sum_{m=0}^{d_\ell}\int_{SX^m}P(t,\mathbf{x},\diff \mathbf{y}).
  \]
  We define a norm for the kernels by
   \[
\|P\|_{(1+E)^{\tilde{\otimes}},\infty}:=\sup_{0\le t\le T,\bfx\in SX^\ell,1\le \ell\le k}\frac{P(t,\bfx)}{(1+E)^{\tilde{\otimes}\ell}(\bfx)}.
\]
We assume that the kernel $P$ is \emph{$(1+E)^{\tilde{\otimes}}$-bounded}, \emph{i.e.} $\|P\|_{(1+E)^{\tilde{\otimes}},\infty}<\infty$.
  \end{enumerate}

\end{assumption}

\begin{remark}
  To avoid technical difficulties in the proof for large deviation, we  keep the ``no dust'' condition holds.
  Nevertheless, with higher  moment control of the initial condition, we can extend our results to the more general case. See Kolokoltsov~\cite{kololln} for instance.
\end{remark}

\begin{remark}
With an abuse of notation, kernel $P_m^\ell(t,\mathbf{x},\diff \mathbf{y})$ is sometimes be written as $P_m^\ell(t,x_1,\dots,x_\ell,\diff y_1,\dots,\diff y_m)$, that is a symmetric function on $X^\ell$ and a symmetric measure on $X^m$.
\end{remark}

For any $j\in\N$, let $\bfx\sim\vec{\bfx}=(x_1,\dots,x_j)\in\cal{X}$ be any configuration, then the generator of the markov process $\bfx(t)$ at state $\bfx$ is given by, for any  continuous function $f$ on $S\cal{X}$
 \begin{equation}\label{gen1}
  G^P_tf(\bfx):=\sum_{\ell=1}^{k\land j}\sum_{I\subset \{1,\dots,j\},|I|=\ell}\sum_{m=0}^{d_\ell}\int_{SX^m}[f({\bfx}_{\overline{I}}\cup{\bfy})-f({\bfx})]P_m^\ell(t,{\bfx}_I,\diff {\bfy})
  \end{equation}
  where $\overline{I}=\{1,2,\dots,j\}\setminus I$, $\bfx_I=\{x_j,j\in I\}$, $\bfx_{\overline{I}}=\{x_j,j\in \overline{I}\}$ and $\bfx_{\overline{I}}\cup{\bfy}$ is set union. It is clear that the definition does not depend on the choice of $\vec{\bfx}$.

\subsection{The Infinitesimal Generator of the Measure-valued Process}

In~\cite{kololln,kolokoltsov2010nonlinear},   Kolokoltsov has shown that with a proper scaling parameter $h>0$, the process $h\delta^+_{\bfx(t)}$ tends to the hydrodynamic equation~\eqref{llneq} as $h\to 0$. For the purpose of large deviation, we need a precise description of the jump space and the measures on the jumps. We form the measure-valued process
\[\mu^h_t:=h\delta^+_{\bfx(t)}\]
as follows.

 Let
\[
\cal{M}_{h\delta}^+(X)=
\left\{h\delta^+_{\mathbf{x}}\big|\mathbf{x}\in S \cal{X}\right\},
\]
be a subset of $\cal{M}^+(X)$
  For any Borel measurable function $F$ on the set $\cal{M}^+_{h\delta}(X)$, any $j\in\N$, for any ${\bfx}\in SX^j$, let
  \begin{multline}\label{gen}
    \Lambda^{h,P}_tF(h\delta^+_{{\bfx}})
    :=\\
    \sum_{\ell=1}^{k\land j}h^{\ell-1}
  \sum_{I\subset \{1,\dots,j\},|I|=\ell}\sum_{m=0}^{d_\ell}\int_{SX^m}
  [F(h\delta^+_{{\bfx}}-h\delta^+_{{\bfx}_I}+h\delta^+_{\bfy})-F(h\delta^+_{{\bfx}})]P(t;{\bfx}_I,\diff \bfy),
  \end{multline}
 where $\bfx_I$, $\bfx_{\overline{I}}$ are determined by fix any $\vec{\bfx}\sim\bfx$. By using~\eqref{gen1}, it is not hard to see that~\eqref{gen} is the generator of the process $\mu^h_t$.

\begin{defi}[Measure on the Quotient Space]
  For all $1\le \ell\le k$, we define a mapping
  $$
  \kappa^\ell[\delta^+_\bfx]:S\cal{X}\mapsto \cal{M}^+(SX^\ell)
  $$
  by, for any Borel measurable function $f$ on $SX^\ell$ and all $\bfx\in S\cal{X}$, one pick any $\vec{\bfx}\in\cal{X}$ such that $\vec{\bfx}\sim \bfx$, let
  \[
\int_{SX^\ell}f(\bfz)\kappa^\ell[\delta^+_\bfx](\diff \bfz):=\sum_{I\subset\{1,\dots,|\bfx|\},|I|=\ell}f(\vec{\bfx}_I),
\]
where $\vec{\bfx}_I=(x_{i_1},\dots,x_{i_j})$  if $\bfx=(x_1,\dots,x_n)$ and $I=\{i_1,\dots,i_j\}\subset \{1,\dots,n\}$.
By convention, for all $|\bfx|<\ell$, let $\sum_{I\subset\{1,\dots,|\bfx|\},|I|=\ell}=0$.

\end{defi}

\begin{defi}[Jump space]\label{jumpspace}
  For all $1\le \ell\le k$, $0\le m\le \ell$, we introduce a space for all the interaction among $\ell$ particles that resulting $m$ particles,
  \[
\cal{J}_{m}^{\ell}:=S{X}^\ell\times S {X}^m;
\]
a  space for all possible interactions among $\ell$ particles,
  \[
\cal{J}^\ell:=\cup_{m=0}^{d_\ell} \cal{J}_{m}^{\ell};
\]
and the space of all possible jumps
\[
\cal{J}:=\cup_{\ell=1}^k\cal{J}^\ell.
\]
Let $\cal{B}(\cal{J}_{m}^{\ell})$, $\cal{B}(\cal{J}^{\ell})$ and $\cal{B}(\cal{J})$ be the corresponding Borel sets.
\end{defi}

By using the above notation, the generator~\eqref{gen} has an expression
  \begin{equation}\label{gen2}
\Lambda^{h,P}_tF(h\delta^+_{\bfx})= h^{\ell-1}\sum_{\ell=1}^{k\land |\bfx|}
\sum_{m=0}^{d_\ell}\int_{\cal{J}^\ell_m}
  [F(h\delta^+_{\bfx}-h\delta^+_{\bfz}+h\delta^+_{\bfy})-F(h\delta^+_{\bfx})]P(t,\bfz,\diff \bfy) \kappa^\ell[\delta^+_{\bfx}](\diff \bfz)
  \end{equation}
  which says that at state $(t,\bfx)$, the jumps from $\bfx$ to $(\bfx\setminus\bfz)\cup\bfy$ happens at rate $$h^{\ell-1}P(t,\bfz,\diff\bfy)\kappa^{\ell}[{\delta^+_{\bfx}}](\diff \bfz)$$ where $\ell=|\bfz|$.

  Moreover, for any linear functional $F_g[\mu]=\langle g,\mu \rangle $,
    \begin{equation}
\Lambda^{h,P}_t \langle g,h\delta^+_{\bfx}\rangle= \sum_{\ell=1}^{k}h^\ell
\sum_{m=0}^{d_\ell}\int_{\cal{J}^\ell_m}
  \langle g,\delta^+_{\bfy}-\delta^+_{\bfz}\rangle P(t,\bfz,\diff \bfy) \kappa^\ell[\delta^+_{\bfx}](\diff \bfz).
  \end{equation}

\begin{prop}[Well-posedness of the finite particle system]\label{wellp}
  For any $\mu_0^h\in\cal{M}_{h\delta}^+(X)$, any kernel $P$ satisfies Assumption~\ref{kernel}, there exists a probability measure $\P^{h,P}_{\mu_0^h}$ on the  space
  $$\left(D_T(\cal{M}^+(X)),(\cal{F}_{s,t})_{0\le s\le t\le T}\right),$$
  where $\cal{F}_{s,t}:=\cap_{\eps>0}\sigma\{(\mu_u):s\le u\le t+\eps\}$ is the nature filtration generated by the coordinate  process $(\mu_t)$
  , such that 
  for all $F\in\cal{C}_b^{1,0}([0,T]\times \cal{M}^+(X))$,
  \[
F(t,\mu_t)-F(s,\mu_s)-\int_s^t\partial F_r(r,\mu_r)\diff r-\int_s^t\Lambda^{h,P}_rF(r,\mu_r)\diff r
\]
is a $((\cal{F}_{s,t})_{0\le s\le t\le T},\P^{h,P}_{\mu_0^h})$-martingale.
\end{prop}

\begin{proof}
We refer to Section~5.1 in Dawson~\cite{dawson} and Section~4.7.A in Ethier and Kurtz~\cite{ethier} for the definitions and proofs of  time in-homogeneous martingale problems. 
\end{proof}

\section{Main Results}\label{sec:main}

In this section, we state the main results (Theorem~\ref{main:ldp}) on the large deviation for the measure valued process $(\mu_t^h)$, whose evolution is governed by the generator $\Lambda_t^{h,P}$ defined in~\eqref{gen2}. The rate function of the upper bound is obtained in a variational form. The rate function of the lower bound is given in a non-variational form.  In Theorem~\ref{ratethm}, we give the relation between these two representations of rate functions.

Our proof for the lower bound relies on the existence of the ``good'' path in a ``Cameron-Martin'' space, who is the domain of the variational rate function (upper bound), in the sense of Fleischmann, G\"artner and Kaj~\cite{fgk}.  Nevertheless, the existence of the paths in such a space is not guaranteed.  Most existence results are proved with kernels that are not growing too fast in the infinity. See DiPerna and Lions~\cite{diperna1989cauchy} for the Boltzmann equation; Jeon~\cite{jeon1998existence}, Norris~\cite{norris1999,norris2000} for the Smoluchovski coagulation equations and Kolokoltsov~\cite{kolokoltsov2010nonlinear} for the general $k$-nary kinetic equations. In Theorem~\ref{gel}, we prove a existence result before gelation for the general $k$-nary kinetic equation~\eqref{llneq} with a $(1+E)^{\tilde{\otimes}}$-bounded kernel. See Norris~\cite{norris1999} and  Fournier et al.~\cite{FOURNIER2009} for the ``gelling'' solution of the Smoluchovski coagulation equation. Our large deviation lower bound thus holds before the gelling time.

We will consider the convergence under two types of initial conditions and two types of assumptions on the asymptotic behaviors of the kernels.

\begin{assumption}\label{condition:kernel}[Kernels: mild growth condition]
  We shall consider the following two  conditions  on the asymptotic behaviors of the kernels in addition to Assumption~\ref{kernel}.
   \begin{itemize}
   \item We say the kernel $P$ satisfying the $(\otimes)$ condition if it is \emph{strongly $(1+E)^{\tilde{\otimes}}$-bounded}:
     \[
\lim_{\langle E,\delta^+_\bfx \rangle\to \infty}\frac{ P(t,\mathbf{x})}{ (1+E)^{\tilde{\otimes} \ell}(\bfx)}=0.
     \]
   \item We say the kernel $P$ satisfies $(\oplus)$ condition if it is \emph{$(1+E)^{+}$-bounded}, \emph{i.e.} $\|P\|_{(1+E)^{+},\infty}<\infty$ where
  \[
\|P\|_{(1+E)^{+},\infty}:=\sup_{0\le t\le T,\bfx\in SX^\ell,1\le \ell\le k}\frac{P(t,\bfx)}{\langle 1+E,\delta^+_\bfx\rangle}.
\]
For the proof of large deviation, we assume additionally in this $(\oplus)$ case that there exists a constant $\lambda_0>0$, such that for all $1\le \ell\le k$, all $0\le m\le d_\ell$
\begin{equation}\label{jumpcontrol}
  \int_{\cal{J}^{\ell}_m}  \left(e^{  |\langle(1+E^\beta),\delta^+_\bfy-\delta^+_\bfx\rangle|}-1\right)P(s,\bfx,\diff\bfy)\le \lambda_0\langle 1+E^\beta,\delta^+_\bfx\rangle,
\end{equation}
where  $\beta>1$ is a constant.

  In the case that the Chaotic Initial Condition holds, we assume $P$ is $1$-non increasing.
  \end{itemize}
\end{assumption}

\begin{assumption}[Initial Conditions]\label{condition:initial}
  We define a constant  $\beta=1$ if the kernel $P$ is \emph{strongly $(1+E)^{\tilde{\otimes}}$-bounded} and $\beta>1$ if the kernel $P$ is \emph{ $(1+E)^{+}$-bounded}.
    We shall consider the following two initial conditions.
    \begin{itemize}
    \item[$\bullet$](Deterministic Initial Condition)
\item[]     Suppose there are a sequence of measures $(\mu_0^h)\in\cal{M}^+_{h\delta}(X)$ and a measure $\nu$ in $\cal{M}^+(X)$, which satisfy
      \[
\sup_{h>0}\langle 1+E^\beta,\mu_0^h\rangle<\infty,\qquad\&\qquad \langle 1+E^\beta,\nu\rangle<\infty,
\]
and $\mu_0^h\to \nu$ in  $\cal{M}^+(X)$ in the weak topology. In this case we define a rate function on $\cal{M}^+(X)$ by
\begin{displaymath}
I_0(\pi) = \left\{ \begin{array}{ll}
0 &\textrm{if }\pi= \nu\\
+\infty &\textrm{if }\pi\neq \nu.
\end{array} \right.
\end{displaymath}

\item[$\bullet$](Chaotic Initial Condition)
\item[] We suppose that $1/h$ are positive integers and there are a sequence of \emph{i.i.d.} $X$-valued random variables $(x_i)_{i\in\N}$ following a probability law $\nu\in \cal{M}^+(X)$  with
\[
\langle e^{\alpha(1+E)^{k}},\nu\rangle<\infty\qquad\textrm{and}\qquad \langle e^{\alpha(1+E^\beta)},\nu\rangle<\infty,\qquad\textrm{for some }\alpha>0.
\]
Let
\[
\mu_0^h=h\sum_{i=1}^{1/h}\delta_{x_i}.
\]
 In this case we define a rate function on $\cal{M}^+(X)$ by the Kullback information
\begin{displaymath}
I_0(\pi) = H(\pi|\nu)=\left\{ \begin{array}{ll}
\langle \log(\frac{\diff \pi}{\diff \nu}),\nu\rangle &\textrm{if }\pi\ll \nu\\
+\infty &\textrm{else. }
\end{array} \right.
\end{displaymath}
    \end{itemize}
  \end{assumption}

\subsection{Law of Large Numbers}

Now we give a proper definition of the limit equation~\eqref{llneq} and state the law of large numbers result. We remark that the proofs in Kolokoltsov~\cite{kololln} for the continuous time-independent kernels applies for our time in-homogeneous measurable kernel case with no essential change provided these kernels satisfies Assumptions~\ref{kernel} and~\ref{condition:kernel}. Nevertheless  we give a short proof of the law of large numbers Theorem~\ref{lln} as a byproduct of the large deviation at the end of Section~\ref{sec:lln}.

\begin{defi}[Limit Measure on the path of jumps]\label{jumptraj}
 For  any measure valued path $\pi\in D_T(\cal{M}^+(X))$, we define a measure $\cal{K}^P[\pi]$ on $[0,T]{\times}\cal{J}$ by
  \[
\cal{K}^{P}[\pi]((0,t]\times U):= \sum_{\ell=1}^k\sum_{m=0}^{d_\ell} \int_0^t\int_{U\cap \cal{J}^\ell_m}  P(s,\bfz,\diff\bfy) \pi^{\tilde{\otimes} \ell}_s(\diff \bfz)\diff s,~~\forall~t\in[0,T],U\in\cal{B}({\cal{J}}).
  \]
  Note that when kernel $P$ is $(1+E)^{\tilde{\otimes}}$-bounded and if
  \[
\sup_{t\in[0,T]} \langle 1+E,\pi_t\rangle<\infty,
\]
then $\cal{K}^{P}[\pi]$ is a  finite measure. For any measurable function $f$ on $[0,T]{\times}\cal{J}$, with an abuse of notation, we denote
\[
\langle f, \cal{K}^{P}[\pi]\rangle:=\int_{0}^T\int_{\cal{J}}f(s,\bfx,\bfy)\cal{K}^{P}[\pi](\tridiff).
\]
\end{defi}

 \begin{theorem}[Existence of the limit equation (Kolokoltsov~\cite{kololln})]\label{exi}
   Suppose that  $P$ is either strongly $(1+E)^{\tilde{\otimes}}$-bounded or  $(1+E)^{+}$-bounded with a initial measure $\nu$ satisfies
  \[
\langle 1+E^\beta,\nu\rangle<\infty
\]
 for some $\beta>1$.
   Then for all $\nu\in\cal{M}_{1+E}(X)$, there exists a continuous curve $t\in[0,T]\mapsto \sigma_t\in\cal{M}_{1+E}(X)$ where $\cal{M}_{1+E}(X)$ endowed with Banach topology, such that the equation~\eqref{llneq},
\[
         \langle g,\sigma_t\rangle=\langle g,\nu\rangle+\int_0^t\int_{\cal{J}}\langle g,\delta^+_{\bfy}-\delta^+_\bfx\rangle \cal{K}^{P}[\sigma](\diff s\diff \bfx,\diff\bfy)
         \]
       holds for all bounded Borel measurable function $g$ on $X$.
 \end{theorem}

 \begin{remark}
   The limit equation~\eqref{llneq} also writes in the following form
      \begin{multline*}
      \langle g,\mu_t\rangle=\langle g,\nu\rangle
      +\sum_{\ell=1}^k\sum_{m=0}^{d_\ell}\frac{1}{\ell!}\int_0^t\diff s\\
  \int_{X^\ell}\prod_{i=1}^\ell\mu_s(\diff x_i)\int_{X^m}    \left(\sum_{j=1}^mg(y_j)-\sum_{i=1}^\ell g(x_i)\right)
      P(s,x_1,\dots,x_\ell,\diff y_1,\dots,\diff y_m).
  \end{multline*}
\end{remark}

\begin{theorem}[Law of Large numbers (Kolokoltsov~\cite{kololln})]\label{lln}
  Suppose that one of the initial Assumption~\ref{condition:initial} holds. If the kernel $P$ satisfies either $(\otimes)$ or $(\oplus)$ condition in Assumption~\ref{condition:kernel}, then the sequence of laws of $\mu^{h}$ is tight on $D_T(\cal{M}^+(X))$, endowed with Skorokhod topology associated with the weak topology on $\cal{M}^+(X)$. Moreover,  almost surely, the limit  $\sigma$ of any converging sub-sequence of $(\mu^h)$ solves~\eqref{llneq} with initial value $\nu$.
\end{theorem}

\subsection{``Cameron-Martin'' Spaces}
It is well-known that in the theory of path-wise large deviations, the domain of the variational rate function is in a ``Cameron-Martin'' space in the sense of Fleischmann, G\"artner and Kaj~\cite{fgk}. See also Dawson and G\"artner~\cite{dawsont1987large}, Kipnis~et~al.~\cite{kipnis1989hydrodynamics} and L\'eonard~\cite{leonard1995large}. We begin to introduce the path in this ``Cameron-Martin'' Space.

\begin{defi}[Absolutely Continuous Kernels (I)]\label{kernel1}
We say that a kernel $P'$ is \emph{absolutely continuous} with respect to the kernel $P$ (denote as $P'\ll P$), if for any measurable set $(\bfx,V)\subset\cal{J}$ and  for almost all $t\in[0,T]$
    \[
P(t,\bfx,V)=0\qquad \Rightarrow \qquad P'(t,\bfx,V)=0.
\]
The \emph{Radon-Nikodym derivative} of kernel $P'$ with respect to the kernel $P$ is given by
    \[
(t,\bfx)\mapsto \eta^{(P'|P)}(t,\bfx,\bfy):=\frac{P'(t,\bfx,\diff \bfy)}{P(t,\bfx,\diff \bfy)}\ge 0 \qquad a.e.P(t,\bfx,\diff \bfy),
    \]
    which is a non-negative measurable function on $[0,T]\times \cal{J}$.

\end{defi}

\begin{defi}[Functional spaces on the path of jumps]
 We shall use the following functional spaces on the measurable space $([0,T]{\times}\cal{J};\cal{K}^P[\pi])$,
  \[
\cal{L}^1(([0,T]{\times}\cal{J});\cal{K}^{P}[\pi]):=\left\{f:([0,T]{\times}\cal{J})\to \R,\textrm{measurable}\bigg|\langle f, \cal{K}^{P}[\pi]\rangle<\infty\right\},
\]
  \[
\cal{L}^\infty(([0,T]{\times}\cal{J});\cal{K}^{P}[\pi]):=\left\{g:([0,T]{\times}\cal{J})\to \R,\textrm{measurable}\bigg| |g |<\infty, \cal{K}^{P}[\pi]\textrm{-a.s.}\right\}.
  \]
 For any $f\in\cal{L}^1(([0,T]{\times}\cal{J};\cal{K}^{P}[\pi])$ and any $g\in \cal{L}^\infty(([0,T]{\times}\cal{J});\cal{K}^{P}[\pi])$, define the duality
 \[
\langle f,g\rangle_{\cal{K}^{P}[\pi]}:=\int_{[0,T]{\times}\cal{J}}f(t,\bfx,\bfy)g(t,\bfx,\bfy) \cal{K}^{P}[\pi](\diff t,\diff \bfx,\diff \bfy).
 \]
\end{defi}

\begin{defi}[Convex Integral Functional]\label{deftau}
  Let $\tau$, $\tau^*$ be the pair of Young functions on $\R$,
  \[
\tau(u)=e^u-u-1,
\]
and 
   \[
  \tau^{*}(u)= \left\{ \begin{array}{ll}
  (u+1)\log(u+1)-u & \textrm{if $u\ge -1$,}\\
+\infty & \textrm{else.}
  \end{array} \right.
  \]
 For $g=\tau,\tau^*$, define a convex integral functional $I_g^{P}(\pi,\cdot)$ by
 \begin{equation}\label{ratei}
I_g^{P}(\pi,f):=\langle g\circ f,  \cal{K}^{P}[\pi]\rangle.
\end{equation}

\end{defi}

\begin{defi}[Cameron-Martin Space]\label{cmspace}
  Fix $\nu\in\cal{M}_{1+E}(X)$. Let $H^P[\nu]$ denote  the set of all paths $\sigma\in C_T(\cal{M}^+(X))$ with
  \begin{enumerate}
  \item $\sigma_0=\nu$;
  \item the map $t\mapsto \sigma_t$ is absolutely continuous;
  \item 
    there exists a kernel $P'$ which is absolutely continuous with respect to $P$ such that $\sigma$ is a solution of equation~\eqref{llneq} with kernel $P'$ and initial value $\nu$, \emph{i.e.}, for almost every $t\in[0,T]$,
    \[
\langle g,\sigma_t\rangle=\langle g,\nu\rangle+\int_0^t\int_{\cal{J}}\langle g,\delta^+_{\bfy}-\delta^+_\bfx\rangle \cal{K}^{P'}[\sigma](\diff s\diff \bfx,\diff\bfy),\qquad \forall g\in \cal{C}_b(X).
    \]
  \item the Radon-Nikodym derivative
   $\eta^{(P'|P)}$
belongs to the  set
\[
\left\{\eta\in\cal{L}^{1}([0,T]\times \cal{J};\cal{K}^P[\sigma])\bigg| \left\langle \tau^*\left(\eta-1\right),\cal{K}^P[\sigma]\right\rangle<\infty\right\}.
\]
  \end{enumerate}
\end{defi}

As we already mentioned, in our $k$-nary interacting case, the existence of the path in $H^P[\nu]$ is not guaranteed. For example, the Smoluchowski equations with  instantaneous gelling kernels does not have a mass conserving solution. See Carr et al.~\cite{carr1992instantaneous} and a Aldous's survey~\cite{Aldous1999}. 
 In general, the existence results are proved by using truncated kernels. Hence, we introduce some sub-spaces of $H^P[\nu]$ such that their paths could be approximated by the paths with bounded Radon-Nikodym derivatives under certain hypotheses.

\begin{defi}[Some Sub-Cameron-Martin Spaces]
  Fix $\nu\in\cal{M}_{1+E}(X)$. We shall consider the following sub-spaces of $H^P[\nu]$ in this paper.
  \begin{enumerate}[label=(\Alph*)]
  \item  Denote by $H_0^P[\nu]$ be  the subset of $H^P[\nu]$, such that the associated Radon-Nikodym derivative $\eta^{(P'|P)}$ is a bounded non-negative Borel measurable function on $[0,T]{\times}\cal{J}$.
 
\item  Denote by $H_{(1+E)^{\tilde{\otimes}}}^{P}[\nu]$  the subset of $H^P[\nu]$ with $(1+E)^{\tilde{\otimes}}$-bounded absolutely continuous kernels, \emph{i.e.}
\begin{equation*}
\|P'\|_{(1+E)^{\tilde{\otimes}},\infty}<\infty.
\end{equation*}
For any $w>0$, let $H_{(1+E)^{\tilde{\otimes}}}^{P,w}[\nu]$ be a subset of $H_{(1+E)^{\tilde{\otimes}}}^{P}[\nu]$ whose absolutely continuous kernel $P'$ satisfies $\|P'\|_{(1+E)^{\tilde{\otimes}},\infty}\le w$. Moreover, we have
\[
    H_{(1+E)^{\tilde{\otimes}}}^{P}[\nu]={\bigcup_{w\in\N}    H_{(1+E)^{\tilde{\otimes}}}^{P,w}[\nu]}.
\]
\item  Denote by $H_{(1+E)_0^{\tilde{\otimes}}}^{P}[\nu]$  the subset of $H^P[\nu]$ with strongly $(1+E)^{\tilde{\otimes}}$-bounded absolutely continuous kernels, \emph{i.e.}
\[
   \|P'\|_{(1+E)^{\tilde{\otimes}},\infty}<\infty,~\&~
    \lim_{\langle E,\delta^+_\bfx \rangle\to \infty,\bfx\in SX^\ell}\frac{ P'(t,\mathbf{x})}{ (1+E)^{\tilde{\otimes} \ell}(\bfx)}=0,\forall 1\le \ell\le k.
\]

\item Denote by $H_{(1+E)^{+}}^{P}[\nu]$  the subset of $H^P[\nu]$ with $(1+E)^{+}$-bounded kernels, \emph{i.e.}
\[
\|P'\|_{(1+E)^{+},\infty}<\infty.
\]
\end{enumerate}
Clearly,
\[
H_0^P[\nu] \subset H_{(1+E)_0^{\tilde{\otimes}}}^{P}[\nu]\subset H_{(1+E)^{+}}^{P}[\nu]\subset H_{(1+E)^{\tilde{\otimes}}}^{P}[\nu]\subset  H^P[\nu].
 \]
\end{defi}

We remark that, according to Theorem~\ref{exi} and its proof, the existence of the path in $H_{(1+E)_0^{\tilde{\otimes}}}^{P}[\nu]$  and  $H_{(1+E)^{+}}^{P}[\nu]$ is proved by using approximation of paths with  truncated kernels. These truncated kernels has finite  Radon-Nikodym derivatives. In other words, the path in $ H_{(1+E)_0^{\tilde{\otimes}}}^{P}[\nu]$  and  $H_{(1+E)^{+}}^{P}[\nu]$ can be approximated by the paths lies in $H_0^P[\nu]$. 

\begin{corollary}[Sub-Cameron-Martin Spaces]
Under the hypothesis of Theorem~\ref{exi}, the paths in $ H_{(1+E)_0^{\tilde{\otimes}}}^{P}[\nu]$  and  $H_{(1+E)^{+}}^{P}[\nu]$  exist. Moreover, we have
\[
H_{(1+E)_0^{\tilde{\otimes}}}^{P}[\nu] \subset \overline{H_0^P[\nu]},
\]
and 
  \[   H_{(1+E)^{+}}^{P}[\nu]\subset \overline{H_0^P[\nu]},
  \]
  where $\overline{H_0^P[\nu]}$ is the completion of $H_0^P[\nu]$ in $C_T(\cal{M}^+(X))$.

 \end{corollary}

For the paths in $ H_{(1+E)^{\tilde{\otimes}}}^{P}[\nu]$, we prove the existence before a finite ``gelling'' moment. The approach is close to  Norris~\cite{norris1999} for the Smoluchowski's equation with product kernel.

 \begin{theorem}[Existence before gelation]\label{gel}
 Let $w>0$ fixed.  Suppose that kernel $P$ is  $(1+E)^{\tilde{\otimes}}$-bounded with $\|P\|_{(1+E)^{\tilde{\otimes}},\infty}\le w$ and the initial $\nu$ satisfies
   \[
\langle 1+E^2,\nu\rangle<\infty.
\]
Then there exists a constant $T_*(w,\nu)\in(0,\infty]$ and a continuous curve $t\in[0,T_*)\mapsto \sigma_t\in\cal{M}_{1+E}(X)$ such that~\eqref{llneq} holds  for all bounded Borel measurable function $g$ on $X$.
Moreover, if $T<T_*(w,\nu)$, then
\[
   H_{(1+E)^{\tilde{\otimes}}}^{P,w}[\nu]\subset \overline{H_0^P[\nu]}.
\]

 \end{theorem}

\begin{remark}
  In general, for any  $(1+E)^{\tilde{\otimes}}$-bounded kernel, the law of large numbers limit of the process may lie outside the space $H^P[\nu]$. For instance, in the limit of Marcus-Lushnikov process, gelation may happen and an correction term may appear in additional to the Smoluchowski's equation. See Norris~\cite{norris1999,norris2000}, Fournier et al.~\cite{Fournier2004} and Fournier et al.~\cite{FOURNIER2009} for examples. Nevertheless, with a modification, our proof on large deviation can be applied to the Marcus-Lushnikov process with certain gelling kernels. The details will be developed in the forthcoming paper~\cite{ldpgel}.
\end{remark}

By using Theorem~\ref{lln}, we can easily obtain the law of large numbers of the process $(\mu^{h,\eta P})$ with kernel $\eta P$ where the  Radon-Nikodym derivative $\eta$ is finite. We will show later in Section~\ref{sec:coupling}, there exists a coupling between the process $\mu^{h,P}$ and  $\mu^{h,\eta P}$ when $\eta$ is finite. We remark that this coupling is essential for the proof of the lower bound and it is more general than the usual coupling between a $f$-perturbed kernel $e^{L[f]}P$ and $P$. We will define this $f$-perturbed kernel later in Section~\ref{sec:coupling}. See also  Kipnis~et~al.~\cite{kipnis1989hydrodynamics}, Jona-Lasinio et al.~\cite{jona1993} and L\'eonard~\cite{leonard1995large}.

\begin{corollary}[Particle approximation in the (sub)-Cameron-Martin Space]\label{appsub}
Under the hypothesis of Theorem~\ref{lln}, for any non-negative $\eta\in\cal{L}^\infty([0,T]{\times}\cal{J})$, let $(\mu^{h,\eta})$ be the  coordinate processes under $\P^{h,\eta P}$, then $\mu^{h,\eta}\Rightarrow\sigma^\eta$ in the same sense as in Theorem~\ref{lln} where $\sigma^\eta\in H_0^P[\nu]$ with Radon-Nikodym derivative $\eta$. Moreover,
    \begin{equation*}
    \lim_{h\to 0}h\cdot H(\P^{h,\eta P}|\P^{h,P})
      = \int_0^T\int_{\cal{J}} \tau^{*}\left(\eta(s,\bfx,\bfy)-1\right)\cal{K}^{P}[\sigma^\eta](\diff s\diff\bfx\diff\bfy),
    \end{equation*}
    where $H(\P^{h,\eta P}|\P^{h,P})$ is the relative entropy between $\P^{h,\eta P}$ and $\P^{h,P}$.
\end{corollary}

\subsection{Main Results on large deviation}
 
We first define the rate functions.

\begin{defi}
 Define a functional $\gamma^{P}: [0,T]{\times}{D_T(\cal{M}_{1+E}(X))}{\times}{ \cal{C}^{1,0}_b([0,T]{\times} X)}{\to}{\R}$ by
  \begin{multline}\label{rategamma}
    \gamma_t^{P}(\pi,g):=\langle \pi_t,g_t \rangle-\langle \pi_0,g_0\rangle-\int_0^t\langle \pi_s,\partial_s g_s\rangle \diff s\\
-\int_0^t\int_{\cal{J}}L[g_s](s,\bfx,\bfy)\cal{K}^P[\pi](\diff s\diff\bfx\diff\bfy),
  \end{multline}
where for any function $f$ on $[0,T]{\times} X$,  $L[f]$ is a function on $[0,T]{\times} \cal{J}$ defined by
  \[
L[f](t,\bfx,\bfy):=\langle f_t,\delta^+_\bfy-\delta^+_\bfx\rangle.
  \]

\end{defi}

\begin{defi}
For any $\pi\in D_T(\cal{M}^+(X))$, let
  \begin{equation}\label{defrate}
  \cal{R}^{P}_{\rm upper}(\pi):=\sup_{g\in \cal{C}_b^{1,0}([0,T];X)} \left\{\gamma_{T}^{P}(\pi,g)
  -I_\tau^{P}(\pi,L[g])\right\}.
  \end{equation}

  For any $\nu\in\cal{M}_{1+E}(X)$ any $\pi\in H_0^P[\nu]$, let
    \begin{equation}\label{defrate2}
  \cal{R}^{P}_{\rm lower}(\pi):=\inf_\eta I_{\tau^*}^P(\pi,\eta)
    \end{equation}
    where $\eta$ running over all possible bounded non-negative measurable functions on $[0,T]{\times}\cal{J}$ such that $\pi$ solves~\eqref{llneq} with kernel $\eta P$ starting from $\nu$. In other words,
\[
\gamma^{\eta P}_t(\pi,g)=0\qquad \forall a.e.-t\in[0,T],\forall g\in \cal{C}_b^{1,0}([0,T]{\times} X).
\]
  \end{defi}

  \begin{theorem}[Large Deviation]\label{main:ldp}
    Suppose that one of the initial Assumption~\ref{condition:initial} holds. If the kernel $P$ satisfies either $(\otimes)$ or $(\oplus)$ condition in Assumption~\ref{condition:kernel}, then for any closed set $\mathcal{C}$ and open set $\mathcal{O}$ of $D_T(\cal{M}^+(X))$, we have
  \[
  \limsup_{h\to 0} h\log \P^{h,P}(\mathcal{C})\le -\inf_{\pi\in\mathcal{C}}
  \left\{I_\nu(\pi_0)+ \cal{R}^{P}_{\rm upper}(\pi)\right\},
\]
\[
\liminf_{h\to 0} h\log \P^{h,P}(\mathcal{O})\ge-\inf_{\pi\in\mathcal{O}\cap H^P_0}
 \left\{I_\nu(\pi_0)+ \cal{R}^{P}_{\rm lower}(\pi)\right\},
\]
where
\[
H^P_0={{\bigcup_{ \substack{\nu'\in\cal{M}^+(X)\\\langle 1+E^\beta,\nu'\rangle<\infty}}H_0^P[\nu'] }}.
\]
Recall that $\beta=1$ if $P$ satisfies $(\otimes)$, $\beta>1$ if $P$ satisfies $(\oplus)$ and $I_\nu$ is defined in Assumption~\ref{condition:initial}.
  \end{theorem}

  In the $1$-nary interacting systems, it is possible to prove the two rate functions are equivalent. Unfortunately, in our general $k$-nary interacting system we could not prove that $ \cal{R}^{P}_{\rm upper}= \cal{R}^{P}_{\rm lower}$ due to the lack of uniqueness in the general Hahn-Banach extension theorem. It highly depends on the structure of the jumps. Inspired by L\'eonard~\cite{leonard1995large}, we can obtain a non-variational representation of the rate function $ \cal{R}^{P}_{\rm upper}$ by considering the conjugation of integral functionals $I_\tau^P(\pi,f)$ on a subspace of $\cal{L}^\infty(([0,T]{\times}\cal{J});\cal{K}^{P}[\pi])$. This non-variational representation has a similar form with the rate function $\cal{R}^{P}_{\rm lower}$. We can give a condition on the jumps $L[g]$ for all $g\in \cal{C}^{1,0}_b([0,T]{\times} X)$, such that
  \[
 \cal{R}^{P}_{\rm upper}= \cal{R}^{P}_{\rm lower} \qquad\textrm{holds on }H_0^P[\nu].
  \]

  \begin{theorem}[Representations of the Rate Functions]\label{ratethm}
    For any $(1+E)^{\tilde{\otimes}}$-bounded kernel $P$, any path $\pi\in D_T(\cal{M}^+(X))$, such that
    \[
\sup_{0\le t\le T}\langle 1+E,\mu_t\rangle<\infty,
    \]
we have
    \begin{equation}\label{rate1}
      \cal{R}^{P}_{\rm upper}(\pi)=\inf_{z\in\cal{O}^P[\pi]}I_{\tau^*}^{P}(\pi,z),
    \end{equation}
    where
\begin{multline*}
  \cal{O}^P[\pi]\\
  =\left\{z\in \cal{L}^1([0,T]{\times} \cal{J};\cal{K}^{P}[\pi])\bigg|
\gamma^{P}_T(\pi,g)=\langle z,L[g]\rangle_{\cal{K}^{P}[\pi]},~\forall g\in \cal{C}^{1,0}_b([0,T]{\times} X)
\right\},
\end{multline*}
and convention $\inf_{\emptyset}=+\infty$. 
Moreover, $ \cal{R}^{P}_{\rm upper}(\pi)<\infty$ if and only if $\pi\in H^P[\pi_0]$.
If there is no nonzero linear functional on $\cal{L}^\infty(([0,T]{\times}\cal{J});\cal{K}^{P}[\pi])$ of the form
\begin{equation}\label{vanish}
  \langle \cdot, w \rangle_{\cal{K}^P[\pi]},\qquad w\in \cal{L}^1([0,T]{\times} \cal{J};\cal{K}^{P}[\pi])
\end{equation}
vanishes throughout $L[g]$ for all $g\in \cal{C}^{1,0}_b([0,T]{\times} X)$, then
    \begin{equation}\label{rate1}
      \cal{R}^{P}_{\rm upper}(\pi)=I_{\tau^*}^{P}(\pi,z),
    \end{equation}
where $\cal{O}^P[\pi]=\{z\}$.

Specially,
\begin{enumerate}[label=$(\alph*)$]
\item if $\pi\in H_0^P[\pi_0]$, then $\cal{R}^{P}_{\rm upper}(\pi)\le \cal{R}^{P}_{\rm lower}(\pi)$;
\item  if $\pi\in H_0^P[\pi_0]$ and condition~\eqref{vanish} holds, then
  \[\cal{R}^{P}_{\rm upper}(\pi)= \cal{R}^{P}_{\rm lower}(\pi)=I_{\tau^*}^{P}(\pi,z);\]
\item if there exists a function $f\in \cal{C}_b^{1,0}([0,T];X)$, such that
\[
\gamma_t^{e^{L[f]}P}(\pi,g)=0,\qquad \forall a.e.~t\in[0,T],~\forall g\in \cal{C}_b^{1,0}([0,T];X),
\]
then
\begin{equation}\label{rate2}
  \cal{R}^{P}_{\rm upper}(\pi)=  \cal{R}^{P}_{\rm lower}(\pi)=I_{\tau^*}^{P}(\pi,e^{L[f]}-1).
\end{equation}
\end{enumerate}
  \end{theorem}

\subsection{Examples}
We apply our results to two examples of binary interacting particle systems: the Becker-D\"oring coagulation-fragmentation model and the Boltzmann collision model. We shall study the similar large deviation problem in Smoluchovski's coagulation in the forthcoming paper~\cite{ldpgel}. Our results provide a lower bound of the Boltzmann model studied in L\'eonard~\cite{leonard1995large}.

\subsubsection{Becker-D\"oring model}\label{exbd}
We take the Becker-D\"oring model as an illustration to the case with deterministic initial condition and  $(\otimes)$-type kernels.

Let the state space $X=\N^*$, the positive integers and $E(i)=i$, the mass of the particle. Given a sequence of coagulation rates $(a_i,i\in\N^*)$ and a sequence of fragmentation rates $(b_i,i\in\N^*)$ with convention $b_1=0$ and assumption
\begin{equation}\label{bdrates}
\inf_{i\ge 1} a_i>0,\qquad \limsup_{i\to\infty}\frac{a_i}{i}<\infty,\qquad
\inf_{i\ge 2} b_i>0,\qquad \limsup_{i\to\infty}\frac{b_i}{i}<\infty,
\end{equation}
the weak form of the  Becker-D\"oring ODE $(c_i(t),i\in\N^*)$ is given by
  \begin{multline}\label{bdode}
   \diff \langle g,c(t)\rangle =\frac{1}{2}\sum_{i=1}^\infty\left(g(i+1)-g(i)-g(1)\right)a_i(1+\ind{i=1}) c_i(t)c_1(t)\diff t\\
  +
 \sum_{i=2}^\infty\left(g(i-1)+g(1)-g(i)\right)b_i c_i(t)\diff t,
  \end{multline}
  for some regular test function $g$ on $\N^*$. We refer to Ball~et al~\cite{Ball} for the well-posedness of this equation.
 Let
  \[ \cal{J}_{\rm Coag}=
  \left\{
          \left(   \{x,1\}   ,     \{x+1\}    \right)\big|x\in\N^*
          \right\},\qquad
          \cal{J}_{\rm Frag}=
  \left\{
          \left(   \{x+1\}   ,     \{x,1\}    \right)\big|x\in\N^*
  \right\},
  \]
  then the jump space of  Becker-D\"oring model can be seen as $\cal{J}=\cal{J}_{\rm Coag}\cup \cal{J}_{\rm Frag}$ with  discrete measure $\diff \bfx\diff \bfy$. For any finite measure $c(t)=(c_i(t),i\in\N^*)$ on $\N^*$,
 the measure $\cal{K}^P[c]$ on $[0,T]{\times}\cal{J}$ is then given by,
  \begin{displaymath}
  \cal{K}^P[c](\diff t\diff \bfx\diff \bfy) = \left\{ \begin{array}{l}
    a_x c_x(t)c_1(t)\delta_{\{x+1\}}(\diff \bfy)\diff \bfx\diff t, \\
  \qquad\qquad\qquad\qquad\qquad  \textrm{if } (\bfx,\bfy)=  \left(   \{x,1\}   ,     \{x+1\}    \right)\in \cal{J}_{\rm Coag};\\
  b_x\delta_{\{1,x-1\}}(\diff \bfy)c_x(t)\diff \bfx\diff t \\
\qquad\qquad\qquad\qquad\qquad \textrm{if } (\bfx,\bfy)=  \left(   \{x-1,1\}   ,     \{x\}    \right)\in \cal{J}_{\rm Frag}.
\end{array} \right.
\end{displaymath}
With these notations,  the hydrodynamic equation~\eqref{llneq} reads,
  \begin{multline*}
  \diff \langle g,c(t)\rangle=\int_{\cal{J}_{\rm Coag}}\langle g,\delta^+_{\{x+1\}}-\delta^+_{\{x,1\}}\rangle \cal{K}^{ P}[c(t)](\diff t\diff \bfx\diff\bfy)\\
  \hfill +\int_{\cal{J}_{\rm Frag}}\langle g,\delta^+_{\{x,1\}}-\delta^+_{\{x+1\}}\rangle \cal{K}^{ P}[c(t)](\diff t\diff \bfx\diff\bfy)
  \end{multline*}
  that coincides with the Becker-D\"oring ODE.

  Now we can give the definition of the Cameron-Martin Space $H^P[c(0)]$. Suppose that $\sum_{i=1}^\infty (1+i)c_i(0)<\infty$, let $H^P[c(0)]$ denote the set of all path $\tilde{c}\in C_T(\cal{M}^+(\N^*))$ such that $\tilde{c}$ is the solution of the Becker-D\"oring equation with initial value $c(0)$, coagulation rates $(\tilde{a}_i(t),i\ge 1)$ and fragmentation rates $(\tilde{b}_i(t),i\ge 2)$ where for almost everywhere $t\in[0,T]$,
  \[
  \tilde{a}_i(t)\ge 0,\forall i\ge 1;\qquad \tilde{b}_i(t)\ge 0,\forall i\ge 2;
  \]
  \[
\frac{1}{2}\int_0^T\sum_{i=1}^\infty\tilde{a}_i(t)\left(1+\ind{i=1}\right)\tilde{c}_i(t)\tilde{c}_1(t)\diff t
+\int_0^T\sum_{i=2}^\infty\tilde{b}_i(t)\tilde{c}_i(t)\diff t<\infty;
  \]
  and
\begin{multline*}
\frac{1}{2}\int_0^T\sum_{i=1}^\infty\tau^*\left(\frac{\tilde{a}_i(t)}{a_i}-1\right)a_i\left(1+\ind{i=1}\right)\tilde{c}_i(t)\tilde{c}_1(t)\diff t\\
+\int_0^T\sum_{i=2}^\infty\tau^*\left(\frac{\tilde{b}_i(t)}{b_i}-1\right)b_i\tilde{c}_i(t)\diff t
<\infty.
\end{multline*}

If in additional that
\[
\sup_{i\ge 1,t\in[0,T]}\frac{\tilde{a}_i(t)}{1+i}<\infty,\qquad \&\qquad \sup_{i\ge 2,t\in[0,T]}\frac{\tilde{b}_i(t)}{1+i}<\infty,
\]
then the path lies in  $H_{(1+E)^{\tilde{\otimes}}}^{P}[c(0)]$.
Clearly, in this model, the  $(1+E)^{\tilde{\otimes}}$-bounded property and the $(1+E)^{+}$-bounded property are equivalent, \emph{i.e.},  $H_{(1+E)^{\tilde{\otimes}}}^{P}[c(0)]$ coincides with $H_{(1+E)^{+}}^{P}[c(0)]$. Hence by using the result and proofs in Ball~et al~\cite{Ball}, we have  the well-posedness of the paths in $H_{(1+E)^{\tilde{\otimes}}}^{P}[c(0)]$ and also
  \[
H_{(1+E)^{\tilde{\otimes}}}^{P}[c(0)]\subset \overline{H_0^P[c(0)]}.
  \]
We remark that for any $f\in\cal{C}^{1,0}_b([0,T]{\times}\N^*)$, the usual $f$-perturbed  kernel (we define rigorously later in Section~\ref{sec:coupling}) is given by
  \[
a^f_i(t):=e^{f(t,i+1)-f(t,i)-f(t,1)}a_i,
  \]
  \[
b^{f}_i(t):=e^{f(t,i-1)+f(t,1)-f(t,i)}b_i,
\]
that satisfies an identity
  \[
a^f_i(t)b_{i+1}^f(t)\equiv a_i b_{i+1}.
\]
It indicates that the hydrodynamic limit of the $f$-perturbed systems are far smaller than the  Cameron-Martin Space $H^P[c(0)]$.

Now we consider a stochastic Becker-D\"oring system with  coagulation rates $(a_i)$ and fragmentation rates $(b_i)$ satisfies the  condition~\eqref{bdrates} and are also sub-linear in the infinity,
\[
\lim_{i\to \infty}\frac{a_i+b_i}{1+i}=0,
\]
starting from $N$ monomers. The Markov process is given by $X^N(t)=(X^N_i(t),k\in\N^*)$, where $X^N_i(t)$ represents the number of particles of mass $i$ at time $t$, with initial condition $X^N(0)=(N,0,\dots)$ and
transition  matrix $Q^N{=}(q^N(\cdot,\cdot))$ 
\begin{displaymath}
\begin{cases}
 q^N(\vec{x},\vec{x}{-}\vec{e}_1{-}\vec{e}_k{+}\vec{e}_{k{+}1}){=}a_k x_1(x_k{-}\ind{k{=}1})/N, \\
 q^N(\vec{x},\vec{x}{+}\vec{e}_{1}{+}\vec{e}_k{-}\vec{e}_{k{+}1}){=}b_{k{+}1}x_{k{+}1}, 
\end{cases}
\end{displaymath}
where  $\vec{x}{\in}\N^{\N^*}$ and $(\vec{e}_k,k{\in}\N^+)$ is the standard orthonormal basis of $\N^{\N^*}$.

Clearly, the kernel satifies the $(\otimes)$ condition and as $N\to\infty$, the initial value
\[
\frac{X^N(0)}{N}\to \delta_1
\]
that  satisfies the Initial Assumption~\ref{condition:initial}. As an application of the results by Jeon~\cite{jeon1998existence}, the hydrodynamic limit of the scaled process $(X^N(t)/N)$ is the Becker-D\"oring ODE~\eqref{bdode}. For the a central limit theorem with bounded kernels, we refer to Sun~\cite{WenBD}. Now we can state the large deviation results.
\begin{corollary}[LDP in Becker-D\"oring Model]
 For any closed set $\mathcal{C}$ and open set $\mathcal{O}$ of $D_T(\cal{M}^+(\N^*))$, we have
  \[
  \limsup_{N\to\infty} \frac{1}{N}\log \P\left(\left(\frac{X^N(t)}{N},0\le t\le T\right)\in\mathcal{C}\right)\le -\inf_{\pi\in\mathcal{C}}
  \left\{I_{\delta_1}(\pi(0))+ \cal{R}^{a,b}_{\rm upper}(\pi)\right\},
\]
\[
\liminf_{N\to\infty}  \frac{1}{N}\log \P\left(\left(\frac{X^N(t)}{N},0\le t\le T\right)\in\mathcal{O}\right)\ge-\inf_{\pi\in\mathcal{O}\cap H^P_0[\delta_1]}
 \left\{I_{\delta_1}(\pi(0))+ \cal{R}^{a,b}_{\rm lower}(\pi)\right\},
 \]
 where   $I_{\delta_1}(\pi(0))=0$ if $\pi(0)=(1,0,\dots)$ and $\infty$ if else; for any $\pi(t)=(\pi_i(t),i\ge 1)\in D_T(\cal{M}^+(\N^*))$,
  \begin{multline*}
  \cal{R}^{a,b}_{\rm upper}(\pi):=\sup_{g\in \cal{C}_b^{1,0}([0,T];\N^{*})} \bigg\{\langle \pi(t),g_t \rangle-\langle \pi(0),g_0\rangle-\int_0^T\langle \pi(t),\partial_t g_t\rangle \diff t\\
  -\frac{1}{2}\int_0^T\sum_{i=1}^\infty\left(\exp\left\{g_t(i+1)-g_t(i)-g_t(1)\right\}-1\right)a_i\left(1+\ind{i=1}\right)\pi_i(t)\pi_1(t)\diff t\\
-\int_0^T\sum_{i=2}^\infty\left(\exp\left\{g_t(i-1)+g_t(1)-g_t(i)\right\}-1\right)b_i\pi_i(t)\diff t\bigg\};
  \end{multline*}
for any  $\pi\in H_0^P[\delta_1]$, 
\begin{multline*}
  \cal{R}^{a,b}_{\rm lower}(\pi):=
  \inf_{\tilde{a}(t),\tilde{b}(t)} \bigg(
\frac{1}{2}\int_0^T\sum_{i=1}^\infty\tau^*\left(\frac{\tilde{a}_i(t)}{a_i}-1\right)a_i\left(1+\ind{i=1}\right)\pi_i(t)\pi_1(t)\diff t\\
+\int_0^T\sum_{i=2}^\infty\tau^*\left(\frac{\tilde{b}_i(t)}{b_i}-1\right)b_i\pi_i(t)\diff t,
  \bigg)
\end{multline*}
where $\tilde{a}(t)=(\tilde{a}_i(t),i\ge 1)$ and $\tilde{b}(t)=(\tilde{b}_i(t),i\ge 2)$ running over all two pair of non-negative measurable functions such that
\[
\sup_{t\in[0,T],i\ge 1}\frac{\tilde{a}_i(t)}{a_i}+\sup_{t\in[0,T],i\ge 2}\frac{\tilde{b}_i(t)}{b_i}<\infty,
\]
and $\pi$ is a solution of the Becker-D\"oring equation with (time in-homogeneous)  parameters $\tilde{a}(t)$ and $\tilde{b}(t)$ starting from $(1,0,\dots)$.
\end{corollary}

We remark that for any $\pi\in H_0^P[\delta_1]$ fixed, the parameters $\tilde{a}(t)$ and $\tilde{b}(t)$ are not unique determined due to the reversity of the jumps, \emph{i.e.} the condition~\eqref{vanish} does not hold. For example,  for any $k\ge 2$, we define parameters $\tilde{a}'(t)$ and $\tilde{b}'(t)$ by
\[
\tilde{a}'_i(t)=\tilde{a}_i(t),~\tilde{b}'_{i+1}(t)=\tilde{b}_{i+1}(t),\qquad \textrm{if }i\neq k,
\]
and 
\[
\tilde{a}'_k(t)=\tilde{a}_k(t)+2\pi_{k+1}(t),~\tilde{b}'_{k+1}(t)=\tilde{b}_{k+1}(t)+\pi_k(t)\pi_1(t),
\]
then $\pi$ is also a solution of the Becker-D\"oring equation with (time in-homogeneous)  parameters $\tilde{a}'(t)$ and $\tilde{b}'(t)$ starting from $(1,0,\dots)$.

\subsubsection{Spatially homogeneous Boltzmann collisions}

Let the state space of velocity be $X=\R^d$. Let $E(z)=|z|^2$ be the kinetic energy. Assume that  collision is described by a $E$-conserved kernel $B(t,z,z',\diff z^*,\diff z'^{*})$ satisfies Assumption~\ref{kernel} and the $(\oplus)$  condition in Assumption~\ref{condition:kernel}. Clearly, it covers Arkeryd's condition~\cite{arkeryd1972boltzmann}
\[
B(z,z'):=\int B(t,z,z',\diff z^*,\diff z'^{*})\le const.(1+|z|^{\lambda}+|z'|^{\lambda})\qquad {\rm with }~\lambda\in(0,2].
  \]
  We remark that our condition~\eqref{jumpcontrol} on the jumps is an analog to the condition~($A_1$) in L\'eonard~\cite{leonard1995large}.

  For all $N\ge 1$, the $N$-particle collision system is described by the Markov process
  $(X_i^N(t),1\le i\le N)$ evolving in $(\R^d)^N$. The law of large numbers and central limit theorems of the associated empirical processes
  \[
\frac{1}{N}\sum_{i=1}^N\delta_{X_i^N(t)}\in\cal{P}(\R^{d})
\]
where $\cal{P}(\R^{d})$ is the probability measure space on $\R^d$ and the hydrodynamic equation of weak type
 \begin{multline}\label{boeq}
\langle \pi_t,g_t \rangle=\langle \nu,g_0\rangle+\int_0^t\langle \pi_s,\partial_s g_s\rangle \diff s\\
+\frac{1}{2}\int_0^t\int_{(\R^d)^2}L[g_s](z,z^*,z',z'^{*})B(s,z,z^*,\diff z',\diff z'^{*}) \pi_s(\diff z)\pi_s(\diff z^*)\diff s,
  \end{multline}
  where
  \[
L[g_s]=g_s(z')+g_s(z'^{*})-g_s(z)-g_s(z'),
\]
are widely studied. See for example Sznitman~\cite{sznitman1984equations} and M\'el\'eard~\cite{meleard}. For the related large deviations, we refer to L\'eonard~\cite{leonard1995large} and two very recent works by Basile et al.~\cite{basile2021large} and Heydecker~\cite{heydecker2021large}.
 
 Now we state our large deviation results.
\begin{corollary}[LDP in spatially homogeneous Boltzmann collisions] Suppose the initial condition is either deterministic or chaotic (Assumption~\ref{condition:initial}), then  for any closed set $\mathcal{C}$ and open set $\mathcal{O}$ of $D_T(\cal{P}(\R^d))$, we have
  \begin{multline*}
  \limsup_{N\to\infty} \frac{1}{N}\log \P\left(\left(\frac{1}{N}\sum_{i=1}^N\delta_{X_i^N(t)},0\le t\le T\right)\in\mathcal{C}\right)\\\le -\inf_{\pi\in\mathcal{C}}
  \left\{I_{\nu}(\pi_0)+ \cal{R}^{B}_{\rm upper}(\pi)\right\},
\end{multline*}
\begin{multline*}
\liminf_{N\to\infty}  \frac{1}{N}\log \P\left(\left(\frac{1}{N}\sum_{i=1}^N\delta_{X_i^N(t)},0\le t\le T\right)\in\mathcal{O}\right)\\\ge-\inf_{\pi\in\mathcal{O}\cap H^B_0[\nu']}
 \left\{I_{\nu}(\pi_0)+ \cal{R}^{B}_{\rm lower}(\pi)\right\},
 \end{multline*}
 where $H^B_0={\bigcup_{ \substack{\nu'\in\cal{P}(X)\\\langle E^\beta,\nu'\rangle<\infty}}H_0^{B}[\nu'] }$ and the rate functions are given by
  \begin{multline*}
  \cal{R}^{B}_{\rm upper}(\pi):=\sup_{g\in \cal{C}_b^{1,0}([0,T];\R^{d})} \bigg\{\langle \pi_t,g_t \rangle-\langle \pi_0,g_0\rangle-\int_0^T\langle \pi_t,\partial_t g_t\rangle \diff t\\
  -\frac{1}{2}\int_0^T\int_{(\R^d)^2}\left(\exp\left\{L[g_t](z,z^*,z',z'^{*})\right\}-1\right)B(t,z,z^*,\diff z',\diff z'^{*}) \pi_t(\diff z)\pi_t(\diff z^{*})\diff t\bigg\};
  \end{multline*}
and for any  $\pi\in H_0^B[\nu']$, 
\begin{multline*}
  \cal{R}^{B}_{\rm lower}(\pi):=\\
  \inf_{\eta} \bigg(
\frac{1}{2}\int_0^T\int_{(\R^d)^2}\tau^*\left(\eta(t,z,z^*,z',z'^{*})-1\right)B(t,z,z^*,\diff z',\diff z'^{*}) \pi_t(\diff z)\pi_t(\diff z^{*})\diff t
  \bigg),
\end{multline*}
  where $\eta$ running over all possible bounded non-negative measurable functions on $[0,T]{\times}(\R^d)^4$ such that $\pi$ solves~\eqref{boeq}
with kernel $\eta B$ starting from $\nu'$.
\end{corollary}

\section{Reformulation: Martingale Measures and Coupling}\label{sec:coupling}
In this section, we introduce a family of orthogonal martingale measures on jump space. As a consequence of this framework, the structure of the Dynkin's martingales related to the martingale problem (See Proposition~\ref{wellp}) is depicted. Then we can obtain the exponential martingales associated with the change of measure between two solutions of martingale problems with different kernels. Due to the structure of the multinary interactions, the classic way for the change of measure by using the test functions (See  Kipnis, Olla and Varadhan~\cite{kipnis1989hydrodynamics} or L\'eonard~\cite{leonard1995large}) is not enough for the proof of the lower bound of a large deviation. Our approach is close to Shiga and Tanaka~~\cite{Shiga1985CentralLT} and Djehiche and Kaj~\cite{djehiche1995}. For more references in the topic of martingale measures, we refer to Chapter 2 of Walsh~\cite{walsh}.
\subsection{A Martingale Measure on Jump Space}

\begin{defi}[Measures on the Jump space]\label{mjs}
  For any $\xi\in D([0,T];\cal{M}^+(X))$, let $\Delta{\xi_t}:={\xi_t}-{\xi_{t-}}$ be the jump function associated with $\xi$. One can define a counting measure  $\cal{N}[\xi](\diff s\diff \bfz \diff \bfy)$ on $[0,T]{\times}\cal{J}$ by, for all $t\in(0,T]$ and any Borel measurable set $U$ in $\cal{J}$,
    \begin{equation}
\cal{N}[\xi]((0,t]{\times}U):=\sum_{0<s\le t}\ind{\Delta {\xi_s}\in \delta_U}
    \end{equation}
    where
    \[
\delta_U:=\left\{\delta^+_{\bfy}-\delta^+_{\bfz}\bigg|(\bfz,\bfy)\in U\right\}.
\]

For any $h>0$, any kernel $P$,  we define a measure $\kappa^{h,P}[\xi]$ on $({[0,T]\times {\cal{J}}})$ by, $\forall t\in(0,T],~\forall U\in\cal{B}(\cal{J}),$
\begin{multline}
\kappa^{h,P}[\xi]((0,t]\times U):=\\
 \left\{ \begin{array}{ll}
  \sum_{\ell=1}^k\sum_{m=0}^{d_\ell} \int_0^t\int_{U\cap \cal{J}^\ell_m} P(s,\bfz,\diff\bfy)h^\ell \kappa^\ell[h^{-1}\xi](\diff \bfz)\diff s & \textrm{if }\xi\in D([0,T];\cal{M}^+_{h\delta}(X))\\
0 & \textrm{else.}
\end{array} \right.
\end{multline}
For all $U\in\cal{B}(\cal{J})$, any $\xi\in D([0,T];\cal{M}^+(X)) $, let $  \tilde{\cal{N}}^{h,P}_0[\xi](U)\equiv 0$ and let
\begin{multline}\label{defitildeN}
  \tilde{\cal{N}}^{h,P}_t[\xi](U)
  :=\\
   \left\{ \begin{array}{ll}
  \cal{N}[h^{-1}\xi]((0,t]{\times}U)-h^{-1}\kappa^{h,P}[\xi]((0,t]\times U) & \textrm{if }\xi\in D([0,T];\cal{M}^+_{h\delta}(X))\\
0 & \textrm{else.}
\end{array} \right.
\end{multline}
\end{defi}


In this context, the measure-valued Markov process $(\mu^h_t)$ obeys the following representation,
  \begin{multline}\label{dy1}
    \mu^h_t=  \mu^h_0+h\int_0^t\int_{\cal{J}} (\delta^+_\bfy-\delta^+_\bfx)\cal{N}[h^{-1}\mu_s^h](\diff s\diff\bfx\diff\bfy)\\
    = \mu^h_0+h\int_0^t\int_{\cal{J}} (\delta^+_\bfy-\delta^+_\bfx)\tilde{\cal{N}}^{h,P}[\mu_s^h](\diff s\diff\bfx\diff\bfy)+\int_0^t\int_{\cal{J}} (\delta^+_\bfy-\delta^+_\bfx)\kappa^{h,P}[ \mu^h_s](\diff s\diff\bfx\diff\bfy).
  \end{multline}
Also, for any $F\in\cal{C}_b( \cal{M}^+(X))$
  \begin{multline}\label{fmg}
    F[\mu^h_t]-F[\mu^h_0]-\int_0^t\Lambda^{h,P}F[\mu^h_s]\\
    =\int_0^t\int_\cal{J}\left(F[\mu^h_s+h\delta^+_\bfy-h\delta^+_\bfx]-F[\mu^h_s]\right)\tilde{\cal{N}}^{h,P}[\mu_s^h](\diff s\diff\bfx\diff\bfy),
  \end{multline}
and for any $g\in\cal{C}_b( X)$,
  \begin{equation}\label{linearmg}
\langle g,\mu^h_t\rangle-\langle g,\mu^h_0\rangle-\int_0^t\Lambda^{h,P}_s\langle g,\mu^h_s\rangle\diff s=h\int_0^t\int_\cal{J}\left(\langle g,\delta^+_\bfy-\delta^+_\bfx\rangle\right)\tilde{\cal{N}}^{h,P}[\mu_s^h](\diff s\diff\bfx\diff\bfy).
    \end{equation}
Clearly, the rightsides of~\eqref{fmg} and~\eqref{linearmg} are Dynkin's martingales (See Proposition~\ref{wellp}). In order to prove a coupling between two processes with different kernels, we need the following stronger results.

\begin{prop}[Cox process, Martingale Measure on the Jump Space]\label{cox}
 Fix $\mu_0^h\in\cal{M}^+_{h\delta}(X)$, then under $\P^{h,P}_{\mu_0^h}$, $\cal{N}[h^{-1}\mu]$ defines a random measure on $(0,T]\times \cal{J}$ whose compensator is given by $h^{-1}\kappa^{h,P}[\mu]$.
    Moreover, the process
  \[
  \left\{\tilde{\cal{N}}^{h,P}[\mu]((0,t]{\times}U),\cal{F}_t,U\in\cal{B}({\cal{J}})\right\}
    \]
    is an orthogonal martingale measure on $\cal{J}$ with the predictable process of bounded variation (Covariance Measure) given by, $\forall U_1,U_2\in\cal{B}({\cal{J}})$,
    \[
    \left\langle
    \tilde{\cal{N}}^{h,P}[\mu]((0,t]{\times}U_1),
      \tilde{\cal{N}}^{h,P}[\mu]((0,t]{\times}U_2) \right \rangle=h^{-1}\kappa^{h,P}[\mu]((0,t]\times U_1\cap U_2)
    \]
\end{prop}

As a byproduct, we can obtain the scaling limit of this orthogonal martingale measure in the sense of Walsh~\cite{walsh} along with the law of large number limit. We briefly talk about the central limit theorem here. The rigorous proof of the convergence in the martingale measure approach is beyond the scope of this paper. For more about the central limit results of $k$-nary interating particle system in the semigroup approach, we refer to Chapter~10 of the book Kolokoltsov~\cite{kolokoltsov2010nonlinear} and the paper~\cite{kolokoltsov2010central} for the Smoluchowski's coaguation model with non-gelling kernels.

\begin{defi}[White noise]
For  any $\pi\in C([0,T];\cal{M}_{1+E}(X))$, let $\cal{W}^P[\pi]$ be the time and spatially inhomogeneous white noise on $[0,T]{\times}\cal{J}$ based on the $\sigma$-finite measure $\cal{K}^P[\pi]$.
\end{defi}

\begin{lemma}[Prop~2.10 of Walsh~\cite{walsh}]\label{lemmawalsh}
Let $M$ be an orthogonal martingale measure, and suppose that for each $U\in \cal{B}(X)$, $t\mapsto M_t(U)$ is continuous. Then $M$ is a white noise if and only if its covariance measure is deterministic.
\end{lemma}

\begin{corollary} \label{mb}

 Under the assmuptios in Theorem~\ref{lln}, for all bounded Borel measurable function  $f$ on $[0,T]\times\cal{J}$, we have
  \[
\sqrt{h}\langle f,\tilde{\cal{N}}^{h,P}[\mu^h]\rangle \Rightarrow \langle f,\cal{W}^P[\sigma]\rangle.
\]
where $\sigma$ is the limit of $\mu^h$, for the convergence in probability under $\E^{h,P}$.
Moreover, for any $\delta>0$, we have
  \[
  \lim_{h\to 0} \P^{h,P}\left(\left|\langle f,\tilde{\cal{N}}^{h, P}[\mu^h]\rangle \right|\ge h^{-1}\delta\right)=0.
  \]
\end{corollary}

We will prove this corollary at the end of Section~\ref{sec:lln}.

\begin{remark}[Central limit theorem]
  Let
  \[
  \theta^h_t:=h^{-1/2}(\mu^h_t-\sigma_t)
  \]
  be the fluctuation process. Then  for all bounded measurable function $g$, as $h\to 0$, we can expect the limit of the real valued random process
  \[
  \left(\langle g,  \theta^h_t\rangle,0\le t\le T\right)\Rightarrow
  \left(\langle g,  \xi_t\rangle,0\le t\le T\right)
  \]
  where
    \begin{multline*}
      \langle g,\xi_t\rangle=\langle g,\xi_0\rangle
      +\int_0^t\sum_{\ell=1}^k\frac{1}{(\ell-1)!}\sum_{m=0}^{d_\ell}\int_{X^\ell}\int_{X^m}\\
      \left(\sum_{j=1}^mg(y_j)-\sum_{i=1}^\ell g(x_i)\right)
      P(s,x_1,\dots,x_\ell,\diff y_1,\dots,\diff y_m)\prod_{i=2}^\ell\sigma_s(\diff x_i)\xi_s(\diff x_1)\diff s\\
      +\int_0^t\int_{\cal{J}}\langle g,\delta^+_\bfy-\delta^+_\bfx\rangle\cal{W}^P[\sigma](\diff s\diff \bfx\diff\bfy).
    \end{multline*}
    Note that the last term is a Gaussian process with covariance.
    \begin{multline*}
      \int_0^t\int_{\cal{J}}\left(\langle g,\delta^+_\bfy-\delta^+_\bfx\rangle\right)^2\cal{K}^P[\sigma](\diff s\diff \bfx\diff\bfy)
      =\sum_{\ell=1}^k\sum_{m=0}^{d_\ell}\int_0^t\diff s\\
   \int_{X^\ell}\frac{1}{\ell!}\prod_{i=1}^\ell\sigma_s(\diff x_i)  \int_{X^m} \left(\sum_{j=1}^mg(y_j)-\sum_{i=1}^\ell g(x_i)\right)^2
      P(s,x_1,\dots,x_\ell,\diff y_1,\dots,\diff y_m).
    \end{multline*}
 
\end{remark}

\subsection{Coupling}

 In addition to the absolutely continuous kernel in Definition~\ref{kernel1}, we introduce the following  kernels that will be used in this paper.
\begin{defi}[Absolutely Continuous Kernels (II)]\label{kernel2}
Suppose the kernel $P$ satisfies the assumptions in Definition~\ref{kernel}.
  \begin{enumerate}
   
\item
    For any $f\in\cal{C}_b^{1,0}([0,T]\times X)$, we define a $f$-perturbed kernel $P^f$ by
  \[
  P^f(t,\bfx,\diff \bfy):= e^{\langle f_t,\delta^+_\bfy-\delta^+_\bfx\rangle} P(t,\bfx,\diff \bfy).
  \]
  Then the corresponding  Radon-Nikodym derivative is given by
  \[
\eta^{(P^f|P)}(t,\bfx,\bfy)= e^{\langle f_t,\delta^+_\bfy-\delta^+_\bfx\rangle}.
\]
\item For any $n>0$,  any   kernel $\eta P$, we define a bounded cut-off kernel  $\eta_nP$ by
\[\eta_n:=\left(\eta\land n\right)\ind{\langle E,\delta^+_\bfx\rangle\le n,\langle E,\delta^+_\bfy\rangle\le n}.\]
  Specially, if $\eta=1$, then the truncated kernel $P_n$ is given by
  \[
P_n(t,\bfx,\diff\bfy)=P(t,\bfx,\diff\bfy)\ind{\langle E,\delta^+_\bfx\rangle\le n,\langle E,\delta^+_\bfy\rangle\le n}.
  \]

  \end{enumerate}
\end{defi}

For any two kernels $P,P'$, any $\mu^h_0\in \cal{M}^+_{h\delta}(X)$, let $\P_{\mu^h_0}^{P}, \P_{\mu^h_0}^{P'}$ be the solution of the martingale problem defined in Proposition~\ref{wellp}.
Inspired by Shiga and Tanaka~\cite{Shiga1985CentralLT}, we introduce the following exponential processes for the purpose of change of measure from $\P_{\mu^h_0}^{P}$ to $\P_{\mu^h_0}^{P'}$. We will show that these exponential processes are real martingales in Proposition~\ref{coupling}.
\begin{defi}[Exponential martingales]\label{defmg}
  For any  kernel $P'\ll P$, for any $\mu\in D([0,T];\cal{M}^+(X))$, define
  \begin{multline}\label{expmg1}
    \cal{M}_t^{h,(P'|P)}[\mu]:=\exp\bigg\{
\int_0^t\int_{\cal{J}}\log \left(\eta^{(P'|P)}(s,\bfx,\bfy)\right)\cal{N}[h^{-1}\mu](s\diff s\diff\bfx\diff\bfy )
    \\-h^{-1}\int_0^t\int_{\cal{J}}\left(\kappa^{h,P'}[\mu](\diff s\diff\bfx\diff\bfy)-\kappa^{h,P}[\mu](\diff s\diff\bfx\diff\bfy)\right)\bigg\}.
  \end{multline}
Specially, for the $f$-pertubated kernel $P^f$, we write,
  \begin{multline}\label{expmg2}
    \cal{M}_t^{h,(P^f|P)}[\mu]=\exp\bigg\{
\int_0^t\int_{\cal{J}}\left(\langle f_s,\delta^+_\bfy-\delta^+_\bfx\rangle\right)\cal{N}[h^{-1}\mu](\diff s\diff\bfx\diff\bfy )
    \\-h^{-1}\int_0^t\int_{\cal{J}}\left(e^{\langle f_t,\delta^+_\bfy-\delta^+_\bfx\rangle}-1\right)\kappa^{h,P}[\mu](\diff s\diff\bfx\diff\bfy)\bigg\}.
  \end{multline}
For any $s,t\in[0,T]$, let
\[
\cal{M}_{s,t}^{h,(P'|P)}[\mu]=\frac{\cal{M}_t^{h,(P'|P)}}{\cal{M}_s^{h,(P'|P)}}[\mu].
\]
\end{defi}

\begin{lemma}[$f$-perturbated particle system]\label{fcoup}
 For any $f\in\cal{C}_b^{1,0}([0,T]\times X)$, the probability $\P^{h,P^f}_{\mu_0^h}$ is absolutely continuous with respect to $\P^{h,P}_{\mu_0^h}$  and the Radon-Nikodym derivative is given by
   \begin{multline}\label{rndf}
     \frac{\diff \P^{h,P^f}_{\mu_0^h}}{\diff\P^{h,P}_{\mu_0^h}}[\mu]_{\big|{\cal{F}_t}}=\cal{M}_t^{h,(P^f|P)}[\mu]\\
     =\exp\bigg\{
    h^{-1}\left(\langle f_t,\mu_t\rangle-\langle f_0,\mu_0\rangle-
    \int_0^t\langle \partial_s f_s,\mu_s\rangle\diff s-\int_0^t\Lambda^{h,P}_s\langle f_s,\mu_s\rangle\diff s\right)
    \\-h^{-1}\int_0^t\int_{\cal{J}}\tau\left(\langle f_s,\delta^+_\bfy-\delta^+_\bfx\rangle\right)\kappa^{h,P}[\mu](\diff s\diff\bfx\diff\bfy)\bigg\}.
   \end{multline}
   Moreover, the relative entropy $H(\P^{h,P^f}_{\mu_0^h}|\P^{h,P}_{\mu_0^h})$ is given by
   \begin{equation*}
h\cdot     H(\P^{h,P^f}_{\mu_0^h}|\P^{h,P}_{\mu_0^h})=\E^{h,P^f}_{\mu_0^h}\left(\int_0^T\int_{\cal{J}}\tau^*\left(e^{\langle f_s,\delta^+_\bfy-\delta^+_\bfx\rangle}-1\right)\kappa^{h,P}[\mu](\diff s\diff\bfx\diff\bfy)\right).
   \end{equation*}
\end{lemma}
\begin{proof}
  This alternative  representation of $\cal{M}_t^{h,(P^f|P)}[\mu]$ is obtained by using equation~\eqref{linearmg}.
The rest of proof is standard. See Kipnis, Olla and Varadhan~\cite{kipnis1989hydrodynamics} for example.
\end{proof}

We now state our main coupling result in this section.
\begin{prop}[Change of measure]\label{coupling}
For any two  \emph{$(1+E)^{\tilde{\otimes}}$-bounded} kernels $P$ and $P'$, any $\mu_0^h\in\cal{M}^+_{h\delta}(X)$ such that $\langle 1+E,\mu_0^h\rangle<\infty$, we have
  \begin{equation}\label{fm1}
\E^{h,P}_{\mu_0^h}\left(\kappa^{h,P'}[\mu^h]((0,T]{\times}\cal{J})\right)<\infty.
  \end{equation}
If  $P'\ll P$ with a  Radon-Nikodym derivative $\eta^{(P'|P)}$ satisfying
  \begin{equation}\label{condition:mg}
   \E^{h,P}_{\mu_0^h} \int_0^T\int_{\cal{J}}\eta^{(P'|P)}(s,\bfx,\bfy)\kappa^{h,P'}[\mu^h](\tridiff)<\infty,
  \end{equation}
then  the process
\begin{equation}\label{mg1}
  \left(\left(\cal{M}_t^{h,(P'|P)}[\mu^h]\right)_{t\ge 0},D\left([0,T];\cal{M}^+(X)\right),\left(\cal{F}_t\right)_{t\ge 0},\P^{h,P}_{\mu_0^h}\right)
  \end{equation}
  is a martingale.
  Moreover, under the probability 
   $\cal{M}_t^{h,(P'|P)}[\mu]\cdot\P^{h,P}_{\mu_0^h}(\diff \mu)$,
   the process
  \begin{equation}\label{mg2}
  \left\{\tilde{\cal{N}}^{h,P'}[\mu^h]((0,t]{\times}U),\cal{F}_t,U\in\cal{B}({\cal{J}})\right\}
    \end{equation}
    is an orthogonal martingale measure on $\cal{J}$ with the predictable process of bounded variation (Covariance Measure) given by, $\forall U_1,U_2\in\cal{B}({\cal{J}})$,
    \[
    \left\langle
    \tilde{\cal{N}}^{h,P'}[\mu^h]((0,t]{\times}U_1),
      \tilde{\cal{N}}^{h,P'}[\mu^h]((0,t]{\times}U_2) \right \rangle=h^{-1}\kappa^{h,P'}[\mu^h]((0,t]\times U_1\cap U_2).
    \]
 That is to say, the probability $\P^{h,P'}_{\mu_0^h}$ is locally absolutely continuous with respect to the probability $\P^{h,P}_{\mu_0^h}$ with density process $(\cal{M}_t^{h,(P'|P)}[\mu^h])_{0\le t\le T}$ and  $\P^{h,P'}_{\mu_0^h}$
is the unique solution of the martingale problem associated with the generator $(\Lambda^{h,{P'}}_t)$ on the function space $\cal{C}_b(\cal{M}^+(X))$.

Moreover, the relative entropy between  $\P^{h,P'}_{\mu_0^h}$ and $\P^{h,P}_{\mu_0^h}$ is given by
\begin{equation}\label{re1}
h\cdot H(\P^{h,P'}_{\mu_0^h}|\P^{h,P}_{\mu_0^h})
      =\E^{h,P'}_{\mu_0^h}  \int_0^T\int_{\cal{J}} \tau^{*}\left(\eta^{(P'|P)}(s,\bfx,\bfy)-1\right)\kappa^{h,P}[\mu^h](\diff s\diff\bfx\diff\bfy).
\end{equation}
\end{prop}

\begin{remark}
  Clearly, if the Radon-Nikodym derivative  $\eta^{(P'|P)}$ is bounded, then the condition~\eqref{condition:mg} is natural. Here we give another two examples.
  \begin{itemize}
  \item When $P$ satsifies condition $(\oplus)$ and the $\langle 1+E^\beta, \mu_0^h\rangle<\infty$ holds for some $\beta>1$, later in Section~\ref{sec:lln} we can prove that for any $1\le \ell\le k$
    \[
\sup_{0\le t\le T}\E^{h,P}_{\mu_0^h} (\langle 1+E^\beta,\mu^h_t\rangle)^\ell<\infty.
    \]
Hence for any nonnegative measurable $\eta$ such that $\eta^2P$ is $(1+E^\beta)^{\tilde{\otimes}}$-bounded, the condition~\eqref{condition:mg} holds. For example in Smoluchowski's model if the coaguation kernel $K(x,y)=1$ and the intial condition has finite second moments, then it can be coupled with the finite particle system with kernel $K'(x,y)=xy$ with $\eta(x,y)=xy$ and $\beta=2$.
    \item
Suppose that $\eta^{(P'|P)}$ does not depend on $\bfy$ and
      \begin{equation}\label{epscond}
\eps:=\min_{1\le \ell\le k,\le m\le d_\ell}\inf_{\substack{t\in[0,T],\bfx\in SX^\ell,\\\langle E,\delta^+_{\bfx}\rangle\le \langle E,h^{-1}\mu_0^h \rangle }}\int_{SX^m}P(t,\bfx,\diff\bfy)>0,
\end{equation}
then
\begin{multline*}
  \E^{h,P}_{\mu_0^h} \int_0^t\int_{\cal{J}}\eta^{(P'|P)}(s,\bfx)\kappa^{h,P'}[\mu^h](\tridiff)\\
  \le 
  \frac{\|P'\|_{(1+E)^{\tilde{\otimes}},\infty}^2}{\eps}  \E^{h,P}_{\mu_0^h} \int_0^t\sum_{\ell=1}^k\frac{1}{\ell!}\left(\langle \left(1+E\right)^2,\mu_s^h\rangle\right)^\ell\diff s\\
  \le 
  \frac{\|P'\|_{(1+E)^{\tilde{\otimes}},\infty}^2T}{\eps} \sup_{0\le s\le T}\E^{h,P}_{\mu_0^h}e^{\langle \left(1+E\right)^2,\mu_s^h\rangle}.
\end{multline*}
The last term is finite if the process with kernel $P$ satisfies certain assumptions. We shall see this propogation of moments results in Section~\ref{sec:lln}. We do not need additional assumptions on $P'$ here. The condition~\eqref{epscond} holds in Smoluchowski's model if the coaguation kernel $\inf_{x,y\in\R^+}K(x,y)>0$; it holds in Becker-D\"oring model if the coagulation rates  $\inf_{i\ge 1}a_i>0$ and fragmentation rates  $\inf_{i\ge 2}b_i>0$; it holds in Boltzmann model if the collision kernel $\inf_{\xi,\xi_*\in\R^d}\int_{\R^{2d}} B(\xi,\xi_*;\diff \xi',\diff \xi_*')>0$.
\end{itemize}

    \end{remark}

\subsection{Proofs}
\begin{proof}[Proof of Proposition~\ref{cox}]
The proof follows Theorem~6.1.3 of Dawson~\cite{dawson} or El~Karoui et~al.~\cite{karoui}. See Lemma 2.4 in Shiga et al.~\cite{Shiga1985CentralLT} for another proof.

Let $\cal{M}_F(X)$ denote the set of all finite measures on $X$.  
By definitions of $\delta^+$ and $\cal{J}$, it is clear that for all $(\bfx,\bfy)\in\cal{J}$, $\delta^+_\bfy-\delta^+_\bfx\in \cal{M}_F(X)$ where $\delta^+_\bfy$ (\emph{resp.}  $\delta^+_\bfx$) is  the positive (\emph{resp.} negative) part of this measure. Therefore, under $\P^{h,P}_{\mu_0^h}$,
\[
\int_{\cal{J}}\delta_{\delta^+_\bfy-\delta^+_\bfx}(\diff \zeta)\cal{N}[h^{-1}\mu](\diff s\diff\bfx\diff\bfy)
\]
defines an adapted random point process on $[0,T]\times \cal{M}_F(X)$ where $\delta_\cdot$ is the delta measure in Definition~\ref{mjs}.
Then the equation~\eqref{dy1} can by written as
  \[
  \mu_t=  \mu_0^h+h\int_0^t\int_{\cal{M}_F(X)}\int_{\cal{J}} \delta_{\delta^+_\bfy-\delta^+_\bfx}(\diff \zeta)\cal{N}[h^{-1}\mu](\diff s\diff\bfx\diff\bfy),
\]
and for all $g\in\cal{C}_b(X)$, the equation~\eqref{linearmg} can by written as
  \begin{equation*}
\langle g,\mu_t\rangle-\langle g,\mu_0^h\rangle=h\int_0^t\int_{\cal{M}_F(X)}\int_\cal{J}\langle g,\zeta\rangle \delta_{\delta^+_\bfy-\delta^+_\bfx}(\diff \zeta){\cal{N}}[h^{-1}\mu](\diff s\diff\bfx\diff\bfy).
  \end{equation*}
  Let $\langle \cal{N}\rangle^{h,P}[h^{-1}\mu]$ be the compensator of $\cal{N}[h^{-1}\mu]$ under $\P^{h,P}_{\mu_0^h}$. Then by It\^o's lemma, the Laplacian gives
  \begin{align*}
  &  \diff e^{-h^{-1}\langle g,\mu_t\rangle+h^{-1}\langle g,\mu_0^h\rangle}\\
    =&e^{-\langle g,\mu_{t-}\rangle+h^{-1}\langle g,\mu_0^h\rangle}\int_{\cal{J}} \left(e^{-\langle g,\delta^+_\bfy-\delta^+_\bfx \rangle}-1\right){\cal{N}}[h^{-1}\mu](\diff t\diff\bfx\diff\bfy)\\
    =&e^{-\langle g,\mu_{t-}\rangle+h^{-1}\langle g,\mu_0^h\rangle}\int_{\cal{J}} \left(e^{-\langle g,\delta^+_\bfy-\delta^+_\bfx \rangle}-1\right)\left({\cal{N}}[h^{-1}\mu]-\langle \cal{N}\rangle^{h,P}[h^{-1}\mu]\right)(\diff t\diff\bfx\diff\bfy)\\
    &+e^{-\langle g,\mu_{t-}\rangle+h^{-1}\langle g,\mu_0^h\rangle}\int_{\cal{M}_F(X)}\int_{\cal{J}} \left(e^{-\langle g,\zeta\rangle}-1\right)\delta_{\delta^+_\bfy-\delta^+_\bfx }(\diff \zeta)\langle \cal{N}\rangle^{h,P}[h^{-1}\mu](\diff t\diff\bfx\diff\bfy).
  \end{align*}

  On the other hand, by using the notation~\eqref{expmg2} for $-g$
  \begin{align*}
&\diff e^{-h^{-1}\langle g,\mu_t\rangle+h^{-1}\langle g,\mu_0^h\rangle}\\
=&\diff \left(\cal{M}_t^{h,(P^{-g}|P)}[\mu]\cdot \exp\left\{h^{-1}\int_0^t\int_{\cal{J}}\left(e^{-\langle g,\delta^+_{\bfy}-\delta^+_\bfx\rangle}-1\right)\kappa^{h,P}[\mu](\diff t\diff \bfx\diff\bfy)\right\}\right)\\
=& \exp\left\{h^{-1}\int_0^t\int_{\cal{J}}\left(e^{-\langle g,\delta^+_{\bfy}-\delta^+_\bfx\rangle}-1\right)\kappa^{h,P}[\mu](\diff t\diff \bfx\diff\bfy)\right\}\diff  \cal{M}_t^{h,(P^{-g}|P)}[\mu_t]\\
&+e^{-h^{-1}\langle g,\mu_{t-}\rangle+h^{-1}\langle g,\mu_0^h\rangle}h^{-1}\int_{\cal{M}_F(X)}\int_{\cal{J}}\left(e^{-\langle g,\zeta\rangle}-1\right)\delta_{\delta^+_\bfy-\delta^+_\bfx }(\diff \zeta)\kappa^{h,P}[\mu](\diff t\diff \bfx\diff\bfy),
  \end{align*}
where $(\cal{M}_t^{h,(P^{-g}|P)}[\mu])$ is a martingale by It\^o's lemma.

By using the uniqueness of the representation of  a locally square integrable martingale plus a process of locally integrable variation for the semimartingale, we have $ \P^{h,P}_{\mu_0^h}$ almost surely,
  \[
 \int_{\cal{J}}\delta_{\delta^+_\bfy-\delta^+_\bfx}(\diff \zeta)\langle \cal{N}\rangle^{h,P}[h^{-1}\mu](\diff t\diff\bfx\diff\bfy)= h^{-1}\int_{\cal{J}}\delta_{\delta^+_\bfy-\delta^+_\bfx}(\diff \zeta)\kappa^{h,P}[\mu](\diff t\diff \bfx\diff\bfy).
  \]
  Since for all set $U\in\cal{B}({\cal{J}})$, $\delta_{U}$ is a  measurable set in $\cal{M}_F(X)$. In conclusion, we have
  \begin{multline*}
    \langle \cal{N}\rangle^{h,P}[h^{-1}\mu](\diff t\times U)=\int_{\delta_U}\int_{\cal{J}}\delta_{\delta^+_\bfy-\delta^+_\bfx}(\diff \zeta)\langle \cal{N}\rangle^{h,P}[h^{-1}\mu](\diff t\diff\bfx\diff\bfy)\\
    =h^{-1}\int_{\delta_U}\int_{\cal{J}}\delta_{\delta^+_\bfy-\delta^+_\bfx}(\diff \zeta)\kappa^{h,P}[\mu](\diff t\diff \bfx\diff\bfy)=h^{-1}\kappa^{h,P}[\mu](\diff t\times U)\\
  \P^{h,P}_{\mu_0^h}~a.s.
  \end{multline*}

\end{proof}

\begin{lemma}\label{cdtc0}
  Suppose kernel $P$ satisfies Assumption~\ref{kernel}, then there exists a constant $c_0$ such that for any $h>0$, any $\mu_0^h\in\cal{M}^+_{h\delta}(X)$ fixed, one has
    \begin{equation}\label{1+emoment}
\sup_{0\le t\le T}{\langle 1+E,\mu_t^h\rangle}<c_0{\langle 1+E,\mu_0^h\rangle}, \qquad \P^{h,P}_{\mu_0^h}-a.s.
\end{equation}
where $(\mu_t^h)$ is the Markov process starting from $\mu_0^h$.
\end{lemma}
\begin{proof}
If $P$ is $1$-non-increasing and $E$-non-increasing , then $c_0=1$. If there exists a constant $\eps_0>0$ such that $E(x)\ge \eps_0$ for all $x\in X$, then  for $\P^{h,P}_{\mu_0^h}-a.s.$ and $t-a.e.$ on $[0,T]$, we have
    \[
\langle 1+E,\mu_t^h\rangle\le \langle \frac{E}{\eps_0}+E,\mu_t^h\rangle\le \left(\frac{1}{\eps_0}+1\right)\langle 1+E,\mu^h_0\rangle.
\]

\end{proof}

\begin{proof}[Proof of Proposition~\ref{coupling}]
Since both kernels $P$ and $P'$ are \emph{$(1+E)^{\tilde{\otimes}}$-bounded}, by using relation~\eqref{1+emoment}, we prove relation~\eqref{fm1} by
  \begin{multline*}
  \E^{h,P}_{\mu_0^h}  \kappa^{h,P'}[\mu^h]((0,T]{\times}\cal{J})=
    \E^{h,P}_{\mu_0^h}  \sum_{\ell=1}^k\int_0^T\int_{SX^\ell}P'(s,\bfx)h^{\ell}\kappa^\ell[h^{-1}\mu^h_s]\diff s\\
      \le
      \|P'\|_{(1+E)^{\tilde{\otimes}},\infty} \E^{h,P}_{\mu_0^h} \int_0^T\sum_{\ell=1}^k\frac{1}{\ell!}\left(\langle(1+E),\mu^h_s\rangle\right)^{\ell}\diff s\\
      \le
            \|P'\|_{(1+E)^{\tilde{\otimes}},\infty} T\sum_{\ell=1}^k\frac{1}{\ell!}\left(c_0\langle(1+E),\mu^h_0\rangle\right)^{\ell}\diff s<\infty.
    \end{multline*}
To prove that the process~\eqref{mg1} is a martingale, we first note that, by using It\^o's formula under $\P^{h,P}_{\mu_0^h}$, 
  \[
  \cal{M}_t^{h,(P'|P)}[\mu]-1=\int_0^t    \cal{M}_{s-}^{h,(P'|P)}[\mu]\int_{\cal{J}}\left(\eta^{(P'|P)}(s,\bfx,\bfy)-1\right)\tilde{\cal{N}}^{h,P}[\mu](\tridiff).
  \]
  According to relation~\eqref{fm1}, for all $t\in[0,T]$, we have
  \begin{multline*}
    \E^{h,P}_{\mu_0^h}\int_0^t\int_{\cal{J}}\left|\eta^{(P'|P)}(s,\bfx,\bfy)-1\right| \kappa^{h,P}[\mu]\\
    \le \E^{h,P}_{\mu_0^h}\int_0^T\int_{\cal{J}}\kappa^{h,P'}[\mu](\tridiff)+ \E^{h,P}_{\mu_0^h}\int_0^T\int_{\cal{J}}\kappa^{h,P}[\mu](\tridiff)<\infty,
  \end{multline*}
  which implies that
  \[
t\mapsto \int_0^t\int_{\cal{J}}\left(\eta^{(P'|P)}(s,\bfx,\bfy)-1\right)\tilde{\cal{N}}^{h,P}[\mu](\tridiff)
\]
is a $\cal{F}_t$-martingale under  $\P^{h,P}_{\mu_0^h}$. See Chapter II Section 3 Page 62 in Ikeda and Watanabe~\cite{ikeda} for example. Moreover, under the condition~\eqref{condition:mg}, for all $t\in[0,T]$,
\[
    \E^{h,P}_{\mu_0^h} \int_0^t\int_{\cal{J}}\eta^{(P'|P)}(s,\bfx,\bfy)\kappa^{h,P'}[\mu](\tridiff)<\infty,
    \]
    then the continuous exponential super martingale $ \cal{M}_t^{h,(P'|P)}[\mu]$ is a real martingale. See Page 144 in Ikeda and Watanabe~\cite{ikeda} for example.

To show that
   the process~\eqref{mg2} is an orthogonal martingale measure on $\cal{J}$  under the probability $\cal{M}_t^{h,(P'|P)}[\mu]\cdot\P^{h,P}_{\mu_0^h}(\diff \mu)$, we first note that the process
  \[
  \tilde{\cal{N}}^{h,P}[\mu  ]((0,t]\times U)-\int_0^t\frac{1}{\cal{M}_{s-}^{h,(P'|P)}[\mu]}\diff \langle \tilde{\cal{N}}^{h,P}((0,s]\times U), \cal{M}_{s}^{h,(P'|P)}[\mu] \rangle
\]
is a  martingale under $ \cal{M}_t^{h,(P'|P)}[\mu]\cdot\P^{h,P}_{\mu_0^h}(\diff \mu)$ by  applying the Girsanov's Theorem (See Theorem III.3.11 in~\cite{jacod}).
To calculate the minus term, we use again It\^o's lemma,  under $\P^{h,P}_{\mu_0^h}(\diff \mu)$,
\begin{multline*}
\left\langle \tilde{\cal{N}}^{h,P}[\mu]((0,t]\times U), \cal{M}_{t}^{h,(P'|P)}[\mu] \right\rangle\\
  =h^{-1}\int_0^t\int_{U}\cal{M}_{s-}^{h,(P'|P)}[\mu]
  \left(\eta^{P'|P}(s,\bfx,\bfy)-1\right)\kappa^{h,P}[\mu](\diff s\diff \bfx\diff\bfy)\\
  =h^{-1}\int_0^t\cal{M}_{s-}^{h,(P'|P)}[\mu]\left(\kappa^{h,P'}[\mu](\diff s\times U)-\kappa^{h,P}[\mu](\diff s\times U)\right).
\end{multline*}
Hence, the process
  \begin{multline*}
    \tilde{\cal{N}}^{h,P'}[\mu]((0,t]\times U)\\
      =
    \tilde{\cal{N}}^{h,P}[\mu  ]((0,t]\times U)-\left(h^{-1}\kappa^{h,P'}[\mu]((0,t]\times U)-h^{-1}\kappa^{h,P}[\mu]((0,t]\times U)\right)
      \end{multline*}
  is a martingale under $\cal{M}_t^{h,(P'|P)}[\mu]\cdot\P^{h,P}_{\mu_0^h}(\diff \mu)$.

  To calculate the covariance measure, for any two measurable sets $U_1,U_2$ in $\cal{J}$,
      similarly, by It\^o's lemma, we have
      \begin{multline*}
        \left\langle \tilde{\cal{N}}^{h,P}_t[\mu](U_1)\tilde{\cal{N}}^{h,P}_t[\mu]( U_2), \cal{M}_{t}^{h,(P'|P)}[\mu] \right\rangle\\
                   =h^{-1}\int_0^t\tilde{\cal{N}}^{h,P}_{s-}[\mu](U_2)\cal{M}_{s-}^{h,(P'|P)}[\mu]
                \left(\kappa^{h,P'}[\mu](\diff s \times U_1)-\kappa^{h,P}[\mu](\diff s \times U_1)\right)\\+
                h^{-1}\int_0^t\tilde{\cal{N}}^{h,P}_{s-}\mu](U_1)\cal{M}_{s-}^{h,(P'|P)}[\mu]
                  \left(\kappa^{h,P'}[\mu](\diff s\times U_2)-\kappa^{h,P}[\mu](\diff s\times U_2)\right)\\
                  +h^{-1}\int_0^t\cal{M}_{s-}^{h,(P'|P)}[\mu]
  \left(\kappa^{h,P'}[\mu](\diff s\times U_1\cap U_2)-\kappa^{h,P}[\mu](\diff s\times U_1\cap U_2)\right).
      \end{multline*}
      By intergration by parts,
         \begin{multline*}
        h\int_{0}^t\frac{1}{\cal{M}_{s-}^{h,(P'|P)}[\mu]}\diff  \left\langle \tilde{\cal{N}}^{h,P}_s[\mu](U_1)\tilde{\cal{N}}^{h,P}_t[\mu]( U_2), \cal{M}_{s}^{h,(P'|P)}[\mu] \right\rangle\\
                  =\tilde{\cal{N}}^{h,P}_{t}[\mu](U_2)
                  \left(\kappa^{h,P'}_t[\mu]-\kappa^{h,P}_t[\mu]\right)( U_1)
                  + \tilde{\cal{N}}^{h,P}_{t}[\mu](U_1)
                    \left(\kappa^{h,P'}_t[\mu]-\kappa^{h,P}_t[\mu]\right)( U_2)\\
                     +
                    \left(\kappa^{h,P'}_t[\mu]-\kappa^{h,P}_t[\mu]\right)( U_1\cap U_2)
                    -
                \int_0^t\tilde{\cal{N}}^{h,P'}[\mu](\diff s \times U_1)
                  \left(\kappa^{h,P'}_s[\mu]-\kappa^{h,P}_s[\mu]\right)( U_2)\\
                   -\int_0^t\tilde{\cal{N}}^{h,P'}[\mu](\diff s \times U_2)
                   \left(\kappa^{h,P'}_s[\mu]-\kappa^{h,P}_s[\mu]\right)( U_1)\\
                   -h^{-1} \left(\kappa^{h,P'}_t[\mu]-\kappa^{h,P}_t[\mu]\right)( U_2)\left(\kappa^{h,P'}_t[ \mu]-\kappa^{h,P}_t[ \mu]\right)( U_1).
         \end{multline*}

On the other hand, by noting that,
\begin{align*}
&h\left(\tilde{\cal{N}}^{h,P'}_t[\mu]( U_1)\tilde{\cal{N}}^{h,P'}_t[\mu]( U_2)-
        \tilde{\cal{N}}^{h,P}_t[\mu]( U_1)\tilde{\cal{N}}^{h,P}_t[\mu](U_2)\right)\\
        =&\left(\kappa^{h,P}_t[ \mu]-\kappa^{h,P'}_t[ \mu]\right)(U_1)\tilde{\cal{N}}^{h,P}_t[\mu](U_2)
        +\left(\kappa^{h,P}_t[ \mu]-\kappa^{h,P'}_t[ \mu]\right)(U_2)\tilde{\cal{N}}^{h,P}_t[\mu](U_1)\\
    &    +h^{-1}\left(\kappa^{h,P}_t[ \mu](U_1)-\kappa^{h,P'}_t[ \mu](U_1)\right)\left(\kappa^{h,P}_t[ \mu](U_2)-\kappa^{h,P'}_t[ \mu](U_2)\right)
\end{align*}
we have
\begin{align*}
 & \tilde{\cal{N}}^{h,P'}_t[\mu]( U_1)\tilde{\cal{N}}^{h,P'}_t[\mu]( U_2)-h^{-1}\kappa^{h,P'}_t[ \mu](U_1\cap U_2)\\
  = & \tilde{\cal{N}}^{h,P}_t[\mu]( U_1)\tilde{\cal{N}}^{h,P}_t[\mu]( U_2)-h^{-1}\kappa^{h,P}_t[ \mu](U_1\cap U_2)\\
 & + \int_{0}^t\frac{1}{\cal{M}_{s-}^{h,(P'|P)}[\mu]}\diff  \left\langle \tilde{\cal{N}}^{h,P}_s[\mu](U_1)\tilde{\cal{N}}^{h,P}_t[\mu]( U_2), \cal{M}_{s}^{h,(P'|P)}[\mu] \right\rangle\\
   & +
                h^{-1}\int_0^t\tilde{\cal{N}}^{h,P'}[\mu](\diff s \times U_1)
                  \left(\kappa^{h,P'}_s[ \mu]( U_2)-\kappa^{h,P}_s[ \mu]( U_2)\right)\\
             &      +h^{-1}\int_0^t\tilde{\cal{N}}^{h,P'}[\mu](\diff s \times U_2)
                   \left(\kappa^{h,P'}_s[ \mu]( U_1)-\kappa^{h,P}_s[ \mu]( U_1)\right).
  \end{align*}
In conclusion, it is a martingale under $\cal{M}_t^{h,(P'|P)}[\mu]\cdot\P^{h,P}_{\mu_0^h}(\diff \mu)$.

   Finally, by using the notation $  \tilde{\cal{N}}^{h,P'}[\mu]$,
the martingale~\eqref{expmg1} can be written as,
  \begin{multline*}
 \cal{M}_t^{h,(P'|P)}[\mu] =\exp\bigg\{\int_{0}^t\int_{\cal{J}}\log\left(\eta^{(P'|P)}(s,\bfx,\bfy)\right)\tilde{\cal{N}}^{h,P'}[\mu](\diff s\diff \bfx\diff \bfy)\\
  +h^{-1}\int_0^t\int_{\cal{J}} \tau^{*}\left(\eta^{(P'|P)}(s,\bfx,\bfy)-1\right)\kappa^{h,P}[ \mu](\diff s\diff\bfx\diff\bfy)
    \bigg\}.
  \end{multline*}
  Thus, the relative entropy between the two finite particle system is given by
  \begin{multline*}
h\cdot H(\P^{h,P'}_{\mu_0^h}|\P^{h,P}_{\mu_0^h})=h\E^{h,P'}_{\mu_0^h}  \log\left( \cal{M}_T^{h,(P'|P)}[\mu] \right)\\
      =\E^{h,P'}_{\mu_0^h}  \int_0^T\int_{\cal{J}} \tau^{*}\left(\eta^{(P'|P)}(s,\bfx,\bfy)-1\right)\kappa^{h,P}[ \mu](\diff s\diff\bfx\diff\bfy).
    \end{multline*}
\end{proof}

\section{A Prior Bounds and Passing to the Limit}\label{sec:lln}

In this section, we begin with a study of measures  $\kappa^{h,P}[\mu]$ and $\cal{K}^P[\mu]$ on the path space of jumps $[0,T]{\times}\cal{J}$ for any fixed path $\mu\in D_T(\cal{M}^+(X))$. We show that the difference between these two measures can be controlled by the moments of $\mu$. The continuity of the functional $\langle f,\cal{K}^P[\cdot]\rangle$, for any bounded measurable function $f\in\cal{L}^\infty([0,T]{\times}\cal{J})$, in the subspaces of $D_T(\cal{M}^+(X))$ is also investigated.

We next provide the estimates on the $\cal{M}^+(X)$-valued Markov process $\mu^h$, such as propagation of moments, exponentially boundedness and convergences of path and fluctuations. For any two kernels $P,P'$, the convergence of $\cal{M}^+(\cal{J})$-valued path $\cal{K}^{P'}[\mu^h]$ where $\mu^h$ is the Markov process with kernel $P$ is also discussed.

\subsection{Some Estimates on the measures $\kappa^{h,P}[\mu]$ and $\cal{K}^P[\mu]$}

The following result gives an explicit expression and bounds for  $\kappa^\ell[\delta^+_\bfx]$. It implies that the measure
\[
\cal{K}^P[\mu]-\kappa^{h,P}[\mu]
\]
defines a non-negative measure on $[0,T]{\times}\cal{J}$.

\begin{prop}[Proposition 2.1~of Kolokoltsov~\cite{kololln}]\label{youngsch}
  For all $\ell\in \N^*$ fixed, there exists constants $\alpha_\Gamma$ (only depends on $\ell$), parameterized by all Young schemes, \emph{i.e.} by all elements in the set
  \[
\Upsilon^\ell:=\left\{\Gamma=(\gamma_1,\gamma_2,\dots,\gamma_r)\bigg| 1\le r\le \ell-1,\gamma_j\in\N^*, \sum_{j=1}^r\gamma_j=\ell,1\le \gamma_1 \le \dots\le \gamma_r\right\},
  \]
such that for all $\bfx\in S\cal{X}$
  \[
\int_{SX^\ell}f(\bfz)\kappa^\ell[\delta^+_\bfx](\diff \bfz)=\int_{SX^\ell}f(\bfz)(\delta^+_{\bfx})^{\tilde{\otimes} \ell}(\diff \bfz)+\sum_{\Gamma\in \Upsilon^\ell}\alpha_\Gamma \int_{X^r}f_\Gamma(y_1,\dots,y_r)\prod_{j=1}^r\delta^+_{\bfx}(\diff y_j)
\]
where $r=|\Gamma|$ and $f_{\Gamma}$ is a Borel measurable function on $X^r$ defined by,
\[
f_{\Gamma}(y_1,y_2\dots,y_r):=f(\underbrace{y_1\dots,y_1}_{\gamma_1},\underbrace{y_2\dots,y_2}_{\gamma_2},\dots,\underbrace{y_r\dots,y_r}_{\gamma_r}).
\]

  For all $f\in\cal{L}^\infty(SX^\ell)$, any $0\le h\le 1$
  \[
\left|h^{\ell}\int_{SX^\ell}f(\bfz)\kappa^\ell[\delta^+_\bfx](\diff \bfz)-\int_{SX^\ell}f(\bfz)(h\delta^+_{\bfx})^{\tilde{\otimes} \ell}(\diff \bfz)\right|\le h\cdot\|f\|_\infty\sum_{\Gamma\in\Upsilon^\ell}|\alpha_\Gamma| (h\ell)^{|\Gamma|}
\]

If additionally $f\ge 0$, then
  \[
h^{\ell}\int_{SX^\ell}f(\bfz)\kappa^\ell[\delta^+_\bfx](\diff \bfz)\le \int_{SX^\ell}f(\bfz)(h\delta^+_{\bfx})^{\tilde{\otimes} \ell}(\diff \bfz)
\]
in other words, for all $\bfx\in S\cal{X}$ fixed, 
$((h\delta^+_{\bfx})^{\tilde{\otimes} \ell}-h^{\ell}\kappa^\ell[\delta^+_\bfx])$ is a  non-negative measure on $SX^\ell$.

\end{prop}

The next result address a rather detailed estimation of the difference between measure $\kappa^{h,P}[\cdot]$ and measure $\cal{K}^P[\cdot]$.
  \begin{lemma}\label{diffk}
    For all bounded measurable function $g:\cal{J}\mapsto \R$ and any measure-valued path $\mu\in D_T(\cal{M}^+(X))$ fixed.
    If $P$ is $(1+E)^{+}$-bounded, 
    then there exists a constant $C_{\oplus}>0$ only depends on $g,T,c_0,P$ and the Young schemes $(\alpha_\Gamma)$ in Proposition~\ref{youngsch} such that  for any $b>0$ and all $0<h\le 1$,  one has 
 \begin{multline*}
     \frac{1}{C_{\oplus}}  \left|  \left\langle g,\left(\kappa^{h,P}[ \mu]-\cal{K}^P[\mu]\right)\right\rangle\right|\\
    \le\frac{1}{b^{(\beta-1)}}
\int_0^T\langle 1+E^{\beta}, \mu_s\rangle\left(1+\left( \langle 1, \mu_s\rangle\right)^{k-1}\right)\diff s
  +  \left(\int_0^T\left(\left(\langle 1,\mu_s\rangle \right)^k+1\right)\diff s\right)(1+b)^kh.
 \end{multline*}

 If kernel $P$ is strongly $(1+E)^{\tilde{\otimes}}$-bounded, 
 then there exists a constant $C_{\otimes}$,  only depends on $g,T,c_0,P$ and the Young schemes $(\alpha_\Gamma)$, such that for any $\eps>0$,  there exists $b(\eps)$ and for all $0<h\le 1$, one has
   \begin{multline*}
     \frac{1}{C_{\otimes}}  \left|  \left\langle g,\left(\kappa^{h,P}[ \mu]-\cal{K}^P[\mu]\right)\right\rangle\right|\\
    \le \left(\int_0^T\left( \left(\langle 1+E,\mu_s\rangle\right)^k+1\right)\diff s\right)\eps
    +  \left(\int_0^T\left(\left(\langle 1,\mu_s\rangle \right)^k+1\right)\diff s\right)(1+b(\eps))^kh.
    \end{multline*}

  \end{lemma}

\begin{proof}[Proof of Lemma~\ref{diffk}]

When the  kernel $P$ is strongly $(1+E)^{\tilde{\otimes}}$-bounded, for any $\eps>0$, there exists some constant $b(\eps)>0$, such that the relation
\[
\sup_{(t\in[0,T],\bfx\in SX^\ell}\frac{P(t,\bfx)}{(1+E)^{\tilde{\otimes}\ell}(\bfx)}\ind{\langle E,\delta^+_\bfx\rangle\ge b}\le \eps.
\]
holds for all $1\le \ell\le k$.
Then on the domain $\{\langle E,\delta^+_\bfx\rangle>b(\eps)\}$, we can control the difference between the two integrals by
  \begin{multline}\label{beps}
 \left|\int_0^T\int_{\cal{J}}g(s,\bfx,\bfy)\ind{\langle E,\delta^+_\bfx\rangle>b}\left(\kappa^{h,P}[ \mu](\tridiff)-\cal{K}^{P}[\mu](\tridiff)\right)\right|
      \\\le\sum_{\ell=1}^k\sum_{m=0}^{d_\ell}\int_0^T \int_{\cal{J}^\ell_m}\left|g\right|(s,\bfx, \bfy)\ind{\langle E,\delta^+_\bfx\rangle>b}P(s,\bfx,\diff \bfy)
   \left(   h^\ell \kappa^{\ell}[h^{-1}\mu_s](\diff \bfx)+\mu_s^{\tilde{\otimes}\ell}(\diff \bfx)\right)\diff s\\
   \le 2\eps \|g\|_\infty \sum_{\ell=1}^k\frac{1}{\ell!}\int_0^T \left(\langle (1+E),\mu_s\rangle\right)^\ell\diff s\\
    \le 2^k\eps \|g\|_\infty k\left(T+\int_0^T \left(\langle (1+E),\mu_s\rangle\right)^k\diff s\right).
  \end{multline}
According to Proposition~\ref{youngsch}, the measure $\cal{K}^{P}[\mu]-\kappa^{h, P}[ \mu]$ is non-negative.
Hence the difference 
\begin{equation*}
  \left|\int_{\cal{J}}g(s,\bfx,\bfy)\ind{\langle E,\delta^+_{\bfx}\rangle\le b}\left(\kappa^{h, P}[ \mu](\diff s\diff\bfx\diff\bfy)-\cal{K}^{P}[\mu](\tridiff)\right)\right|
\end{equation*}
does not exceed
\[
 2\|g\|_\infty\sum_{\ell=1}^k\int_{SX^\ell}P(s,\bfx)\ind{\langle E,\delta^+_{\bfx}\rangle\le b}\left((\mu_s)^{\tilde{\otimes}\ell}(\diff\bfx)-\kappa^{\ell}[h^{-1}\mu_s](\diff \bfx)\right).
\]
Again by using Proposition~\ref{youngsch}, there exists a constant $\alpha$, such that for all $2\le \ell\le k$,
\begin{multline*}
\frac{1}{\|P\|_{(1+E)^{\tilde{\otimes}},\infty}}\int_{SX^\ell}P(s,\bfx)\ind{\langle E,\delta^+_{\bfx}\rangle\le b}\left((\mu_s)^{\tilde{\otimes}\ell}(\diff\bfx)-\kappa^{\ell}[h^{-1}\mu_s](\diff \bfx)\right)
\\
\le 
\sum_{r=1}^{\ell-1}h^{\ell-r}\sum_{1\le \gamma_1\le \dots\le \gamma_r,\sum_{j=1}^r\gamma_j=\ell}
\int_{X^r}\prod_{j=1}^r(1+E(x_j))^{\gamma_j}\ind{\sum_{j=1}^r\gamma_jE(x_j)\le b}\mu_s(\diff x_j)\\
\le 
(1+b)^\ell\alpha\sum_{r=1}^{\ell-1}h^{\ell-r}
\left(\langle 1,\mu_s\rangle \right)^r
\le 
(1+b)^k\alpha k
2^{k-1}\left(\left(\langle 1,\mu_s\rangle \right)^k+1\right)h
\end{multline*}
holds when $h\le 1$.

In conclusion, for any $\eps>0$, there exists $b(\eps)>0$ such that
for all $h\le 1$,
 \begin{multline*}
    \left|  \left\langle g,\left(\kappa^{h,P}[ \mu]-\cal{K}^P[\mu]\right)\right\rangle\right|
    \le  2^k \|g\|_\infty k\left(\int_0^T\left( \left(\langle 1+E,\mu_s\rangle\right)^k+1\right)\diff s\right)\eps\\
    + \|g\|_\infty k^2 \alpha\|P\|_{(1+E)^{\tilde{\otimes}},\infty}
2^{k} \left(\int_0^T\left(\left(\langle 1,\mu_s\rangle \right)^k+1\right)\diff s\right)(1+b(\eps))^kh.
 \end{multline*}

 Now suppose the kernel $P$ satisfies condition~$(\oplus)$.
For any $\bfx\in SX^\ell$, and $(x_1,\dots,x_\ell)\sim\bfx$, let
\[
\max(E(\bfx)):=\max\{E(x_1),E(x_2),\dots,E(x_\ell)\}.
\]

 We can first estimate the difference between the two integrals on the set $\{\max\{E(\bfx)\}> b\}$ for any $b>0$ by
\begin{multline}\label{bbeta}
\frac{1}{2\|g\|_{\infty} \|P\|_{(1+E)^{+},\infty}}  \left|\int_0^T\int_{\cal{J}} g(s,\bfx,\bfy)\ind{\max\{E(\bfx)\}> b}\left(\kappa^{h,P}[ \mu]-\cal{K}^P[\mu]\right)(\tridiff)\right|\\
\le \sum_{\ell=1}^k\frac{1}{(\ell-1)!}\int_0^T\left(\langle (1+E)\ind{E> b}, \mu_s\rangle\left( \langle 1, \mu_s\rangle\right)^{\ell-1} \right)\diff s\\
\le k\left(\frac{2}{1+b}\right)^{\beta-1}
\int_0^T\langle 1+E^{\beta}, \mu_s\rangle\left(1+\left( \langle 1, \mu_s\rangle\right)^{k-1}\right)\diff s.
\end{multline}
For the difference on the set $\{\max\{E(\bfx)\}\le b\}$, we can repeat the estimation as in the  $(\otimes)$ case.
In conclusion, for all $b>0$, all $h\in(0,1]$, we have
 \begin{multline*}
    \left|  \left\langle g,\left(\kappa^{h,P}[ \mu]-\cal{K}^P[\mu]\right)\right\rangle\right|\\
    \le 2k\left(\frac{2}{1+b}\right)^{\beta-1}\|g\|_{\infty} \|P\|_{(1+E)^{+},\infty}
\int_0^T\langle 1+E^{\beta}, \mu_s\rangle\left(1+\left( \langle 1, \mu_s\rangle\right)^{k-1}\right)\diff s\\
    \|g\|_\infty k^2 \alpha\|P\|_{(1+E)^{\tilde{\otimes}},\infty}
2^{k} \left(\int_0^T\left(\left(\langle 1,\mu_s\rangle \right)^k+1\right)\diff s\right)(1+b)^kh.
 \end{multline*}

  \end{proof}

Next, we investigate the continuity of functional $\langle f, \cal{K}^P[\cdot]\rangle$ on the subspaces $\cal{D}_T(\cal{M}^+(X))$.

  \begin{prop}[Continuity]\label{ctnfk}
  For all $n\in\N$,   we define a subset of $\cal{D}_T(\cal{M}^+(X))$ by
    \begin{equation}\label{dn}
      \cal{D}_T^n:=\left\{\mu\in\cal{D}_T(\cal{M}^+(X))\bigg| \sup_{0\le t\le T}\langle 1+E^\beta,\mu_t \rangle\le n\right\},
    \end{equation}
    and a set
    \begin{equation}\label{dn2}
      \tilde{\cal{D}}_T^n:=\left\{\mu\in\cal{D}_T(\cal{M}^+(X))\bigg| \int_0^T\langle 1+E^\beta,\mu_t \rangle\diff t\le n,\sup_{0\le t\le T}\langle 1,\mu_t\rangle\le C_\#\right\},
    \end{equation}
    where $C_\#$ is a  constant.
    
      For any $f\in\cal{L}^\infty([0,T]{\times}\cal{J})$ fixed, the functional
      $\langle f,\cal{K}^P[\cdot]\rangle$ is continuous on $\cal{D}_T^n$, in the sense that,
for any sequence
 of $(\mu^{h})$  and $\sigma$ in $ \cal{D}_T^n$,
 if  for all bounded Borel measurable function $g$ on $X$,
  \[
 \lim_{h\to 0}\sup_{t\in[0,T]}\left|\langle g,\mu^{h}_t \rangle-\langle g,\sigma_t \rangle\right|=0,
 \]
then  one has
 \[
\lim_{h\to 0}\langle  f,\cal{K}^P[\mu^h] \rangle = \langle  f,\cal{K}^P[\sigma] \rangle.
\]
Specially, if kernel  $P$ satisfies condition~$(\oplus)$, then the functional
      $\langle f,\cal{K}^P[\cdot]\rangle$ is also continuous on $\tilde{\cal{D}}_T^n$ in the same sense.
  \end{prop}

  \begin{proof}[Proof of Proposition~\ref{ctnfk}]
Suppose that the kernel $P$ is strongly $(1+E)^{\tilde{\otimes}}$-bounded. Recall that in this case $\beta=1$. It is easy to  see that in this case for all $\bfx\in SX^\ell$, all $0\le m\le d_\ell$ all $b>0$, the function
\begin{multline}\label{bdtest}
\left|\int_{SX^m} f(s,\bfx,\bfy)P(s,\bfx,\diff\bfy)\ind{\langle E,\delta^+_{\bfx}\rangle\le b}\right|\\
\le \|f\|_{\infty}\|P\|_{(1+E)^{\tilde{\otimes}},\infty}(1+E)^{\tilde{\otimes}\ell}(\bfx)\ind{\langle E,\delta^+_{\bfx}\rangle\le b}\le \|f\|_{\infty}\|P\|_{(1+E)^{\tilde{\otimes}},\infty}(1+b)^\ell,
\end{multline}
that is a bounded function on $[0,T]\times\cal{J}$.
Thus, as $h\to 0$, we have
\begin{multline*}
\int_0^T\int_{\cal{J}}f(s,\bfx,\bfy)\ind{\langle E,\delta^+_\bfx\rangle\le b}\cal{K}^{P}[\mu^{h}](\tridiff)\\
\to
\int_0^T\int_{\cal{J}}f(s,\bfx,\bfy)\ind{\langle E,\delta^+_\bfx \rangle\le b}\cal{K}^{P}[\sigma](\tridiff).
\end{multline*}
By a small modification of the inequality~\eqref{beps} in the proof of Lemma~\ref{diffk}, we can have that for any $\eps>0$, there exists $b(\eps)$, such that
\begin{multline}\label{tail1}
\left|\int_0^T\int_{\cal{J}}f(s,\bfx,\bfy)\ind{\langle E,\delta^+_\bfx\rangle> b(\eps)}\cal{K}^{P}[\mu^{h}](\tridiff)\right|\\
\le 2^{k-1}\eps\|f\|_\infty kT\left(1+\left(\sup_{0\le t\le T}\langle 1+E,\mu_t\rangle\right)^k\right),
\end{multline}
that does not exceed
\[
2^{k-1}\eps\|f\|_\infty kT\left(1+n^k\right),
\]
on $\cal{D}_T^n$.
The proof completes by first taking $\eps\to 0$ and then taking $h\to 0$.

 When the kernel $P$ satisfies condition~$(\oplus)$, similarly, the relation~\eqref{bdtest} holds on the set $\ind{\max\{E(\bfx)\le b}$.
Thus,
\begin{multline*}
\int_0^T\int_{\cal{J}}f(s,\bfx,\bfy)\ind{\max\{E(\bfx)\}\le b}\cal{K}^{P}[\mu^{h}](\tridiff)\\
\to
\int_0^T\int_{\cal{J}}f(s,\bfx,\bfy)\ind{\max\{E(\bfx)\}\le b}\cal{K}^{P}[\sigma](\tridiff).
\end{multline*}
By a small modification of the inequality~\eqref{bbeta}, we obtain that for any $b>0$, such that
\begin{multline}\label{tail2}
\left|\int_0^T\int_{\cal{J}}f(s,\bfx,\bfy)\ind{\max\{E(\bfx)\}> b}\cal{K}^{P}[\mu^{h}](\tridiff)\right|\\
\le \frac{C}{(1+b)^{(\beta-1)}} \int_0^T\langle 1+E^{\beta}, \mu_s\rangle\left(1+\left( \langle 1, \mu_s\rangle\right)^{k-1}\right)\diff s.
\end{multline}
that does not exceed
\[
\frac{C}{(1+b)^{(\beta-1)}} Tn\left(1+n^{k-1}\right)
\]
on $\cal{D}_T^n$ and does not exceed
\[
\frac{C(1+C_\#^{k-1})}{(1+b)^{(\beta-1)}} n
\]
on $\tilde{\cal{D}}_T^n$.
The proof completes by first taking $b\to \infty$ and then taking $h\to 0$.

\end{proof}

  \subsection{Some Estimates of the Process $(\mu^h_t)$}
  Recall Lemma~\ref{cdtc0}, there exists a constant $c_0$ such that for any $h>0$, any $\mu_0^h\in\cal{M}^+_{h\delta}(X)$ fixed, one has for any constant $a>0$, any $\ell$,
 \[
\exp\{a\sup_{0\le t\le T}(\langle 1+E,\mu_t^h\rangle)^\ell\}<\exp\{a(c_0{\langle 1+E,\mu_0^h\rangle})^\ell\}, \qquad \P^{h,P}_{\mu_0^h}-a.s.
\]
where $(\mu_t^h)$ is the Markov process starting from $\mu_0^h$.

  \begin{lemma}[Propagation of Moments]\label{boundness}
If $P$ satisfies $(\oplus)$ condition,
  we define a constant $C_\#$ by
\begin{displaymath}
C_\# = \left\{ \begin{array}{ll}
1 & \textrm{if  the chaotic initial condition holds;}\\
c_0\sup_{h}\langle 1+E,\mu_0^h\rangle & \textrm{if if the deterministic initial condition holds}
\end{array} \right.
\end{displaymath}
then for all $h>0$, $\P^{h,P}$ almost surely, one has
\begin{equation}\label{1bound}
\sup_{0\le t\le T}\langle 1,\mu_t^h\rangle\le C_{\#}.
\end{equation}
  There exists  constants $c_\beta$ and $b_\beta$, such that for any $h>0$, any $\mu_0^h\in\cal{M}^+_{h\delta}(X)$ fixed, we have
  \begin{equation}\label{betamoment}
\E^{h,P}_{\mu_0^h}\sup_{t\in[0,T]}\langle E^\beta,\mu^h_t\rangle\le
c_\beta\langle E^\beta,\mu^h_0\rangle+b_\beta,
  \end{equation}
 and  for any $\gamma$ small enough, there exists a constant $c'$ such that
  \begin{equation}\label{expbetamoment}
\sup_{0\le t\le T}\E^{h,P}_{\mu_0^h}\exp\left\{\gamma h^{-1}\langle 1+E^\beta,\mu^h_t \rangle\right\}
\le \exp\left\{c' \gamma h^{-1}\langle 1+E^\beta,\mu^h_0 \rangle\right\}.
\end{equation}
\end{lemma}

  \begin{proof}[Proof of Lemma~\ref{boundness}]
 
To prove~\eqref{betamoment}, we  apply the Dynkin's formula and use the fact that  the kernel $P$ is $E$-non-increasing and then we obtain
\begin{multline*}
  \langle E^\beta,\mu_t^h\rangle=\langle E^\beta,\mu_0^h\rangle+
  h\int_0^t\int_{\cal{J}}\langle E^\beta,\delta^+_{\bfy}-\delta^+_{\bfx}\rangle\cal{N}[h^{-1}\mu_s^h](\tridiff)\\
  \le \langle E^\beta,\mu_0^h\rangle+
  h\int_0^t\int_{\cal{J}}\left(\langle E,\delta^+_{\bfy}\rangle^\beta-\langle E^\beta,\delta^+_{\bfx}\right)\rangle\cal{N}[h^{-1}\mu_s^h](\tridiff)\\
   \le \langle E^\beta,\mu_0^h\rangle+
  h\int_0^t\int_{\cal{J}}\left(\langle E,\delta^+_{\bfx}\rangle^\beta-\langle E^\beta,\delta^+_{\bfx}\right)\rangle\cal{N}[h^{-1}\mu_s^h](\tridiff).
\end{multline*}  
Since $\langle E,\delta^+_{\bfx}\rangle^\beta-\langle E^\beta,\delta^+_{\bfx}\rangle$ is non-negative when $\beta>1$ and  $\cal{K}^{P}[\mu^h]-\kappa^{h,P}[\mu^h]$ is a non-negative measure (See Proposition~\ref{youngsch}), we can have the following Gr\"onwall type inequality,
\begin{multline*}
\E^{h,P}_{\mu_0^h}\sup_{0\le r\le t}\langle E^\beta,\mu_r^h\rangle\le\langle E^\beta,\mu_0^h\rangle\\\hfill
+
h\E^{h,P}_{\mu_0^h}\int_0^t\int_{\cal{J}}\left(\langle E,\delta^+_{\bfx}\rangle^\beta-\langle E^\beta,\delta^+_{\bfx}\right)\rangle\cal{N}[h^{-1}\mu_s^h](\tridiff)\\
\hfill
\le \langle E^\beta,\mu_0^h\rangle+
\E^{h,P}_{\mu_0^h}\int_0^t\int_{\cal{J}}\sup_{0\le r\le s}\left(\langle E,\delta^+_{\bfx}\rangle^\beta-\langle E^\beta,\delta^+_{\bfx}\right)\rangle\cal{K}^{P}[\mu^h_r](\diff s\diff\bfx\diff\bfy).\\
\end{multline*}
We complete the proof of relation~\eqref{betamoment} by following the proof of Proposition~4.2 in Kolokoltsov~\cite{kololln}.

To prove~\eqref{expbetamoment}, we shall use a function $f$ given by,
\[
f(t,x):=\alpha e^{-\lambda t}(1+E^\beta).
\]
Then it is easy to see that
\begin{multline*}
  \cal{M}^{h,(P^f|P)}_t[\mu^h]\\=\exp\Bigg\{\alpha h^{-1}\left(e^{-\lambda t}\langle 1+E^\beta,\mu^h_t \rangle-\langle 1+E^\beta,\mu^h_0 \rangle+\int_0^t\lambda e^{-\lambda s}\langle 1+E^\beta,\mu^h_s\rangle\diff s \right)\\
  -h^{-1}\int_0^t\int_{\cal{J}}\left(e^{\alpha e^{-\lambda s}\langle  (1+E^\beta),\delta^+_\bfy-\delta^+_\bfx\rangle}-1\right)\kappa^{h,P}[\mu^h](\tridiff)\Bigg\}
\end{multline*}
is a supermartingale under $\P^{h,P}_{\mu_0^h}$. By Proposition~\ref{youngsch}, we see that $\P^{h,P}_{\mu_0^h}$ almost surely, $(\cal{K}^P[\mu^h]-\kappa^{h,P}[\mu^h])$ is a positive measure.  If $\alpha\le 1$, by using $a(e^{\alpha  y/a}-1)\le \alpha (e^y-1)$ for all $y\in \R$, all $a>0$, we get
\begin{multline*}
  \cal{M}^{h,(P^f|P)}_t[\mu^h]\\\ge \exp\Bigg\{\alpha h^{-1}\left(e^{-\lambda t}\langle 1+E^\beta,\mu^h_t \rangle-\langle 1+E^\beta,\mu^h_0 \rangle+\int_0^t\lambda e^{-\lambda s}\langle 1+E^\beta,\mu^h_s\rangle\diff s \right)\\
  -h^{-1}\alpha \int_0^t\int_{\cal{J}} e^{-\lambda s} \left(e^{ |\langle(1+E^\beta),\delta^+_\bfy-\delta^+_\bfx\rangle|}-1\right)\cal{K}^P[\mu^h](\tridiff)\Bigg\}.
\end{multline*}

By the assumption~\eqref{jumpcontrol},
 let $\lambda\ge   \lambda_0 k^2\left(C_{\#}\right)^{k-1}$, then
\begin{multline*}
\int_0^t e^{-\lambda s} \int_{\cal{J}} \left(e^{|\langle(1+E^\beta),\delta^+_\bfy-\delta^+_\bfx\rangle|}-1\right)\cal{K}^P[\mu^h](\tridiff)\\
  \le \sum_{\ell=1}^k\sum_{m=0}^{d_\ell} \int_0^t e^{-\lambda s}  \int_{\cal{J}^{\ell}_m}  \left(e^{  |\langle(1+E^\beta),\delta^+_\bfy-\delta^+_\bfx\rangle|}-1\right)P(s,\bfx,\diff\bfy)(\mu^h_s)^{\tilde{\otimes}\ell}(\diff\bfx)\diff s\\
  \le\int_0^t e^{-\lambda s}  \lambda_0\sum_{\ell=1}^k\frac{d_\ell}{(\ell-1)!}\langle 1+E^\beta,\mu_s^h\rangle\left(\langle 1,\mu_s^h\rangle\right)^{\ell-1}\diff s\\
  \le \lambda\int_0^t e^{-\lambda s}\langle 1+E^\beta,\mu_s^h\rangle\diff s 
\end{multline*}
holds $\P^{h,P}_{\mu_0^h}$ almost surely. Therefore,
\begin{equation*}
  \cal{M}^{h,(P^f|P)}_t[\mu^h]\ge \exp\left\{ h^{-1}\left(\alpha e^{-\lambda t}\langle 1+E^\beta,\mu^h_t \rangle-\langle 1+E^\beta,\mu^h_0 \rangle\right)\right\}\qquad \P^{h,P}_{\mu_0^h}-a.s.
\end{equation*}
Hence, for any $\gamma\in(0,e^{-T\lambda_0k^2(C_\#)^{k-1}}]$,
\[
\E^{h,P}_{\mu_0^h}\exp\left\{\gamma h^{-1}\langle 1+E^\beta,\mu^h_t \rangle\right\}
\le \exp\left\{\gamma e^{T\lambda_0 k^2\left(C_\#\right)^{k-1}} h^{-1}\langle 1+E^\beta,\mu^h_0 \rangle\right\}.
\]

\end{proof}

\begin{lemma}[Initial Moments]\label{hlogexp}
 Suppose  one of the initial condition holds, then one has
  \begin{equation}\label{initial:kmoment}
    \sup_h\E\left(\langle 1+E,\mu_0^h\rangle\right)^{k}<\infty,
  \end{equation}
  and 
  \begin{equation}\label{initial:log}
\limsup_{h\to 0}h\log \E\exp\left\{h^{-1}\alpha\left(\langle 1+E,\mu_0^h\rangle\right)^{k}\right\}<\infty.
  \end{equation}
Specially, in the chaotic initial  case with $\beta>1$, one has additionally,
  \begin{equation}\label{initial:betalog}
\limsup_{h\to 0}h\log \E\exp\left\{h^{-1}\alpha\left(\langle 1+E^\beta,\mu_0^h\rangle\right)\right\}<\infty.
  \end{equation}
\end{lemma}

\begin{proof}[Proof of Lemma~\ref{hlogexp}]
  In the deterministic case, the proof is trivial. In the chaotic case, by noticing that $\mu_0^h$ is a probability measure and then using Jensen's inequality, we have
  \[
  \E\left(\langle 1+E,\mu_0^h\rangle\right)^{k}\le
    \E\langle \left(1+E\right)^{k},\mu_0^h\rangle= \langle (1+E)^{k},\nu\rangle<\infty,
    \]
    and
    \begin{multline*}
      h\log \E\exp\left\{h^{-1}\alpha\left(\langle 1+E,\mu_0^h\rangle\right)^{k}\right\}\\
      \le h\log \E\exp\left\{h^{-1}\alpha\langle (1+E)^{k},\mu_0^h\rangle\right\}=
      \log \langle e^{\alpha (1+E)^{k}},\nu\rangle,
    \end{multline*}
      \begin{equation*}
      h\log \E\exp\left\{h^{-1}\alpha\left(\langle 1+E^\beta,\mu_0^h\rangle\right)\right\}=
      \log \langle e^{\alpha (1+E^\beta)},\nu\rangle.
    \end{equation*}
\end{proof}

  \begin{prop}\label{diffexp}

    Suppose  one of the initial condition holds.
    Let $f$ be any fixed bounded Borel measurable function on $[0,T]{\times}\cal{J}$.
   We consider two kernels $P$ and $P'$ that both satisfy either condition $(\otimes)$ or  condition  $(\oplus)$ and  let $(\mu^{h})$ be the  coordinate processes under $\P^{h, P}$,
   then
   \begin{equation}\label{bd}
\sup_{h}\E^{h,P}\left(\left|\left\langle f,\cal{K}^{P'}[\mu^h]\right\rangle\right|+\left|\left\langle f,\kappa^{h,P'}[ \mu^h]\right\rangle\right|\right)<\infty
   \end{equation}
\begin{equation}\label{kcvg}
\lim_{h\to 0} \E^{h,P}\left|\left\langle f,\left(\cal{K}^{P'}[\mu^h]-\kappa^{h,P'}[ \mu^h]\right)\right\rangle\right|=0,
\end{equation}
and
\begin{equation}\label{expkcvg}
\lim_{h\to 0} h\log \left(\E^{h,P}\exp\left\{h^{-1}\left|\left\langle f,\left(\cal{K}^{P'}[\mu^h]-\kappa^{h,P'}[ \mu^h]\right)\right\rangle\right|\right\}\right)=0.
\end{equation}
Moreover, if there is a sub-sequence of $(\mu^{h})$  converges  to $\sigma\in D([0,T];\cal{M}^+(X))$ in the sense that for all bounded Borel measurable function $g$ on $X$,
  \[
 \sup_{t\in[0,T]}\left|\langle g,\mu^{h}_t \rangle-\langle g,\sigma_t \rangle\right|\to 0,
  \]
    in probability under $\P^{h,P}$ as $h\to 0$.
  Then one has  that 
\begin{equation}\label{wkcvg}
\int_0^T\int_{\cal{J}}f(s,\bfx,\bfy)\kappa^{h,P'}[ \mu^{h}](\diff s\diff\bfx\diff\bfy)\to\int_0^T\int_{\cal{J}}f(s,\bfx,\bfy)\cal{K}^{P'}[\sigma](\tridiff)
\end{equation}
   in probability under $\E^{h,P}$  as $h\to 0$.
  \end{prop}

  \begin{proof}[Proof of Proposition~\ref{diffexp}]
    Thanks to relations~\eqref{1+emoment} and~\eqref{initial:kmoment}, the boundedness~\eqref{bd} is obvious,
    \begin{multline*}
      \sup_{h}\E^{h,P}\left(\left|\left\langle f,\cal{K}^{P'}[\mu^h]\right\rangle\right|+\left|\left\langle f,\kappa^{h,P'}[ \mu^h]\right\rangle\right|\right)
\\
  \le 2 \|f\|_\infty  \|P'\|_{(1+E)^{\tilde{\otimes}},\infty}\E^{h,P}\int_0^T \sum_{\ell=1}^k\frac{1}{\ell!}\left(\langle (1+E),  \mu^h_s\rangle\right)^\ell\diff s\\
  \le  
 2^{k} \|f\|_\infty k T\left(1+\sup_h\E\left(c_0\langle (1+E),\mu^h_0\rangle\right)^k\right)<\infty.
\end{multline*}

In the case that both kernels $P$ and $P'$ satisfy condition $(\otimes)$, by lemma~\ref{diffk} and relation~\eqref{1+emoment}, there is a constant $C>0$,  depend on $f, P'$, such that
\begin{multline*}
\frac{1}{C}\E^{h,P}\left|\left\langle f,\left(\cal{K}^{P'}[\mu^h]-\kappa^{h,P'}[ \mu^h]\right)\right\rangle\right|\\
  \le
  \eps \left( 1+\E\left(\langle 1+E,\mu_0^h\rangle \right)^k\right)+ (1+b(\eps))^kh\left( 1+\E\left(\langle 1+E,\mu_0^h\rangle \right)^k\right),
\end{multline*}
and
\begin{multline*}
  h\log \left(\E^{h,P}\exp\left\{h^{-1}\left|\left\langle f,\left(\cal{K}^{P'}[\mu^h]-\kappa^{h,P'}[ \mu^h]\right)\right\rangle\right|\right\}\right)
  \le
  \eps C+ (1+b(\eps))^kh C\\
  +h\log \left(\E\exp\left\{h^{-1}\left(
 \eps C\left(\langle 1+E,\mu_0^h\rangle \right)^k+ (1+b(\eps))^khC\left(\langle 1+E,\mu_0^h\rangle \right)^k
 \right)\right\}\right)
\end{multline*}
By using H\"older's inequality, the last line does not exceed
\begin{multline*}
  \eps C\frac{1}{2h^{-1} \eps C}\log \left(\E\exp\left\{2h^{-1}\eps C\left(\langle 1+E,\mu_0^h\rangle \right)^k\right\}\right)\\
  +h(1+b(\eps))^kC\frac{1}{2(1+b(\eps))^kC}\log \left(\E\exp\left\{
 2(1+b(\eps))^kC\left(\langle 1+E,\mu_0^h\rangle \right)^k
\right\}\right).
\end{multline*}
By Lemma~\ref{hlogexp}, when we take
\[
h\le \eps^2 \land\frac{1}{(1+b(\eps))^{k+1}}
\]
and let $\eps\to 0$, we can obtain the convergences~\eqref{kcvg} and~\eqref{expkcvg}.

In the case that both kernels $P$ and $P'$ satisfy condition $(\oplus)$, by Lemma~\ref{diffk} and relation~\eqref{1bound}, there is a constant $C>0$,  depend on $f, P'$, and by Lemma~\ref{boundness} relation~\eqref{betamoment} a constant $C'$ depends on $\beta$ such that
  \begin{multline*}
  \frac{1}{C}  \E^{h,P}  \left|\left\langle f,\left(\cal{K}^{P'}[\mu^h]-\kappa^{h,P'}[ \mu^h]\right)\right\rangle\right|\\
      \le (1+b)^kh
      +\frac{C'T}{(1+b)^{\beta-1}} \E \left(\left(\langle 1+E^\beta,\mu^h_0\rangle+1\right)\right).
     \end{multline*}
By considering $h\le 1/(1+b)^{k+1}$ and let $b\to \infty$, we obtain the limit~\eqref{kcvg}.
To prove the limit~\eqref{expkcvg}, we first using Lemma~\ref{diffk} and get
\begin{multline*}
  h\log \left(\E^{h,P}\exp\left\{h^{-1}\left|\left\langle f,\left(\cal{K}^{P'}[\mu^h]-\kappa^{h,P'}[ \mu^h]\right)\right\rangle\right|\right\}\right)
  \le
(1+b)^kh C\\
    +h\log \left(\E^{h,P}\exp\left\{h^{-1}\left(
   \frac{CC'}{(1+b)^{\beta-1}} \int_0^T\langle 1+E^\beta,\mu_s^h\rangle\diff s
 \right)\right\}\right).
\end{multline*}
By lemma~\ref{boundness}, there exists a constant $c'$, such that for $b$ large enough, for any $h$, the last term is not larger than
\begin{multline*}
h\log \left(\frac{1}{T}\int_0^T\E^{h,P}\exp\left\{h^{-1}\left(
   \frac{C'}{(1+b)^{\beta-1}} T\langle 1+E^\beta,\mu_s^h\rangle
   \right)\right\}\diff s\right)\\
   \le
   h\log \left(\E\exp\left\{h^{-1}\left(
   \frac{C'c'}{(1+b)^{\beta-1}} T\langle 1+E^\beta,\mu_0^h\rangle
   \right)\right\}\right).
\end{multline*}
Let $h\le \frac{1}{(1+b)^{2(\beta-1)}}$, then when $b\to \infty$, this term is vanishing.

The proof of convergence~\eqref{wkcvg} is quite similar to the proof of Proposition~\ref{ctnfk}. We only need to check that in the $(\otimes)$ case, by inequalities~\eqref{tail1} and~\eqref{1+emoment}, the tail integration under $\P^{h,P}$ has a uniform upper bound over $h>0$,
\begin{multline*}
\E^{h,P}\left|\int_0^T\int_{\cal{J}}f(s,\bfx,\bfy)\ind{\langle E,\delta^+_\bfx\rangle> b(\eps)}\cal{K}^{P}[\mu^{h}](\tridiff)\right|\\
\le 2^{k-1}\eps\|f\|_\infty kT\left(1+c_0\sup_{h}\E\left(\langle 1+E,\mu_0^h\rangle\right)^k\right),
\end{multline*}
which is vanishing as $\eps\to 0$.

In the $(\oplus)$ case, by inequalities~\eqref{tail2}~\eqref{1bound} and~\eqref{betamoment}, the tail integration under $\P^{h,P}$ has a uniform upper bound over $h>0$,
\begin{multline*}
\E^{h,P}\left|\int_0^T\int_{\cal{J}}f(s,\bfx,\bfy)\ind{\max\{E(\bfx)\}> b}\cal{K}^{P}[\mu^{h}](\tridiff)\right|\\
\le \frac{C(1+C_\#^{k-1})}{(1+b)^{(\beta-1)}} T\left(c_\beta\sup_{h}\E\langle E^{\beta}, \mu_0^h\rangle+C_\#+b_\beta\right),
\end{multline*}
which is vanishing as $b\to \infty$.
  \end{proof}
 
\subsection{Passing to the limit: LLN \& Fluctuation}

\begin{proof} [Proof of Theorem~\ref{lln}] Later we shall prove in Lemma~\ref{exptight}, the probabilities $\P^{h,P}$ are exponential tight on $D_T(\cal{M}^+(X))$. As a byproduct we have
the tightness of $(\langle g,\mu^{h}_t \rangle)$ for all bounded Borel measurable function $g$ on $X$ and the uniform norm for processes.
  Combine with Proposition~\ref{diffexp}, we have the law of large number. See also Kolokoltsov~\cite{kololln}.

\end{proof}

\begin{proof}[Proof of Corollary~\ref{appsub}]
Clearly, the law of large number holds for the kernel $\eta P$ when $\|\eta\|_\infty<\infty$. Moreover, since $\eta\ge 0$
\[
\|\tau^{*}\left(\eta (s,\bfx,\bfy)-1\right)\|_\infty\le \tau^{*}\left(\|\eta\|_\infty-1\right)+1<\infty,
\]
thus by using Proposition~\ref{diffexp}, relation~\eqref{re1} has limit
\begin{multline*}
\lim_{h\to 0}  h\cdot H(\P^{h,\eta P}|\P^{h,P})=
  \lim_{h\to 0}\E^{h,\eta P} \int_0^T\int_{\cal{J}} \tau^{*}\left(\eta (s,\bfx,\bfy)-1\right)\kappa^{h,P}[ \mu^h](\diff s\diff\bfx\diff\bfy)\\
  = \int_0^T\int_{\cal{J}} \tau^{*}\left(\eta(s,\bfx,\bfy)-1\right)\cal{K}^{P}[\sigma^{\eta }](\diff s\diff\bfx\diff\bfy).
\end{multline*}

\end{proof}

\begin{proof}[Proof of Corollary~\ref{mb}]
  For all bounded Borel measurable function $f$ on $[0,T]\times \cal{J}$, under $\P^{h,P}$, the process
  $$\left(\sqrt{h}\int_0^t\int_{\cal{J}} f(s,\bfx,\bfy)\tilde{\cal{N}}^{h,P}[\mu^h](\diff s \diff \bfx\diff\bfy) \right)$$
  is a martingale with previsible increasing process
  \[
 \left(\int_0^t\int_{\cal{J}} f(s,\bfx,\bfy)^2 \kappa^{h,P}[ \mu^h](\diff s\diff\bfx\diff\bfy)\right).
 \]
By Theorem~\ref{lln} and Proposition~\ref{diffexp}, we have
\begin{equation*}
  \int_0^t\int_{\cal{J}} f(s,\bfx,\bfy)^2 \kappa^{h,P}[ \mu](\diff s\diff\bfx\diff\bfy)
  \to 
\int_0^t\int_{\cal{J}} f(s,\bfx,\bfy)^2 \cal{K}^P[\sigma](\diff s\diff\bfx\diff\bfy)
\end{equation*}
in probability as $h\to 0$. The proof completes according to Lemma~\ref{lemmawalsh}.

\end{proof}

\section{Proofs of Large deviation}\label{sec:ldp}
In this section we first prove Theorem~\ref{ratethm} for the alternative representation of the rate function $\cal{R}^{P}_{\rm upper}$ defined in~\eqref{defrate} by following the idea of L\'eonard~\cite{leonard1995large}. In order to simplify the proof, we consider the duality in the spaces $ \cal{L}^{\infty}([0,T]{\times} \cal{J};\cal{K}^P[\pi])$ and  $\cal{L}^\infty([0,T]{\times} \cal{J};\cal{K}^P[\pi])^*$ instead of the duality in the Orlicz spaces $\cal{L}_{\tau}([0,T]{\times} \cal{J};\cal{K}^P[\pi])$ and $\cal{L}_{\tau^*}([0,T]{\times} \cal{J};\cal{K}^P[\pi])$, where for $\phi=\tau,\tau^*$, (See Definition~\ref{deftau})
\begin{multline*}
  \cal{L}_{\phi}([0,T]\times \cal{J};\cal{K}^P[\pi])\\
  =\left\{z:[0,T]\times \cal{J}\mapsto \R\textrm{ measurable}\bigg|\inf\left\{a>0:\left\langle\phi\left(\frac{z}{a}\right),\cal{K}^P[\pi]\right\rangle\right\}<\infty\right\}.
\end{multline*}
Our result is based on an extension of the paper Rockafellar~\cite{rock1}. The proof is included in the Appendix~\ref{sec:app}. For more about the duality in the Orlicz spaces, we refer to  L\'eonard~\cite{leonard2001convex}.

The rest of this section is devoted to the proof of the large deviation results, Theorem~\ref{main:ldp}. We first prove the exponential tightness, then prove the upper bound by using the $f$-perturbed particle system (See Lemma~\ref{fcoup}) and finally the lower bound by using the general change of measure (See Proposition~\ref{coupling}) combined with Theorem~\ref{ratethm}.

\begin{proof}[Proof of Theorem~\ref{ratethm}]
  When  kernel $P$ is $(1+E)^{\tilde{\otimes}}$-bounded and $\sup_{0\le t\le T}\langle 1+E,\pi_t\rangle<\infty,$ the measure $\cal{K}^P[\pi]$ is finite on the space $[0,T]{\times}\cal{J}$. Then by Theorem~\ref{thm:ro} (Rockafellar's theorem) in the Appendix~\ref{sec:app}, the functional $I_{\tau}^P(\pi,\cdot)$ is a well-defined finite convex function on $\cal{L}^\infty([0,T]{\times}\cal{J};\cal{K}^P[\pi])$ and it is everywhere continuous with respect
  to the uniform norm. Consider the subspace
  \[
L\cal{C}:=  \left\{L[g]\in  \cal{L}^{\infty}([0,T]{\times} \cal{J};\cal{K}^{P}[\pi])\bigg|\forall  g\in \cal{C}_b^{1,0}([0,T]{\times} X)\right\},
\]
for any $h\in L\cal{C}$, we pick any $g\in \cal{C}_b^{1,0}([0,T]{\times} X)$ such that $L[g]=h$ and define
\[
\jmath^P_T(\pi,h):=\gamma_T^P(\pi,g).
\]
We first claim that $\jmath_T^P(\pi,\cdot)$ is a well-defined continuous linear functional on $L\cal{C}$ when $\cal{R}^{P}_{\rm upper}(\pi)<\infty$. Suppose that there exists two $g_1,g_2$ such that $L[g_1]=L[g_2]$ and $\gamma_T^P(\pi,g_1)>\gamma_T^P(\pi,g_2)$, Then for any $a>0$,
\begin{multline*}
  \cal{R}^{P}_{\rm upper}(\pi)\ge \gamma_T^P(\pi,g_2+a(g_1-g_2))-I_{\tau}^P(\pi,L[g_2+a(g_1-g_2)])\\
  =\gamma_T^P(\pi,g_2)+a\left(\gamma_T^P(\pi,g_1)-\gamma_T^P(\pi,g_2)\right)-I_{\tau}^P(\pi,L[g_2]),
\end{multline*}
that leads to a contradiction by letting $a\to \infty$. The linearity is trivial. To prove the continuity, for any $h\in L\cal{C}$, we pick any $g\in \cal{C}_b^{1,0}([0,T]{\times} X)$ such that $h=L[g]$, then we have $\|h\|_\infty\le 2k\|g\|_\infty<\infty$ and
\[
  \frac{1}{\|h\|_\infty}\jmath_T^P(\pi,h)=\gamma_T^P(\pi,\frac{g}{\|L[g]\|_\infty})\le \cal{R}^{P}_{\rm upper}(\pi)+
  I_\tau^{P}(\pi,\frac{L[g]}{\|L[g]\|_\infty}).
\]
Since
\[
\left\|\tau\left(\frac{L[g]}{\|L[g]\|_\infty}\right)\right\|_\infty\le e+2<\infty,
\]
and
\[
\langle e+2,\cal{K}^P[\pi]\rangle<\infty,
\]
then $\jmath_T^P(\pi,\cdot)\in ( L\cal{C})^*$.
Finally, by  using the Corollary~\ref{subconj} in the Appendix~\ref{sec:app},
\[
  \cal{R}^{P}_{\rm upper}(\pi)=(I^P_\tau(\pi,\cdot))^*(\jmath_T^P(\pi,\cdot))=\inf_{z\in\cal{O}^P[\pi]}I_{\tau^*}^P(\pi,z)
  \]
  where
\[
    \cal{O}^P[\pi]
    =\left\{z\in \cal{L}^1([0,T]{\times} \cal{J};\cal{K}^{P}[\pi])\bigg|
\jmath^{P}_T(\pi,h)=\langle z,h\rangle_{\cal{K}^{P}[\pi]},~\forall h\in L\cal{C}\right\}.
\]
Hence, $\cal{R}^P_{\rm upper}(\pi)<\infty$ if and only if there exists a $\eta\in \cal{L}^1([0,T]{\times} \cal{J};\cal{K}^{P}[\pi])$ that $\langle \tau^*(\eta-1),\cal{K}^P[\pi] \rangle<\infty$, such that
\[
\gamma^{\eta P}_t(\pi,g)=0\qquad \forall a.e.-t\in[0,T],\forall g\in \cal{C}_b^{1,0}([0,T]{\times} X).
\]
Moreover, if the condition~\eqref{vanish} holds, then by Corollary~\ref{co:ro} in the Appendix~\ref{sec:app},  we have
$$\cal{R}^{P}_{\rm upper}(\pi)=I_{\tau^*}^P(\pi,z)$$ and $ \cal{O}^P[\pi]=\{z\}$.

The proofs of $(a)$ and $(b)$ are trivial. To prove~\eqref{rate2}, we first note that for any $f,g\in \cal{C}_b^{1,0}([0,T]{\times} X)$,
  \[
  \gamma_t^P(\pi,g)-\gamma_t^{e^{L[f]}P}(\pi,g)=
  \int_0^t\int_{\cal{J}}L[g]\left(e^{L[f]}-1\right)(s,\bfx,\bfy)\cal{K}^P[\pi](\tridiff).
  \]
  If $\gamma_T^{e^{L[f]}P}(\pi,g)=0$, then
  \[
  \gamma_T^P(\pi,g)=\left\langle \left(e^{L[f]}-1\right),L[g] \right\rangle_{\cal{K}^P[\pi]},\qquad \forall g\in \cal{C}_b^{1,0}([0,T]{\times} X).
  \]  
  By using relation~\eqref{rate1} and the fact $e^{L[f]}-1\in \cal{L}^1([0,T]{\times} \cal{J};\cal{K}^P[\pi])$,
  \[
\left\langle \left|e^{L[f]}-1\right|,\cal{K}^P[\pi]\right\rangle\le \left(e^{2k\|f\|_\infty}+1\right)\left\langle 1,\cal{K}^P[\pi]\right\rangle<\infty,
  \]
  we have
  \[
\cal{R}^P_{\rm upper}(\pi)\le I_{\tau^*}^{P}(\pi,e^{L[f]}-1).
  \]
  On the other hand, by definition~\eqref{defrate}, we take $g=f$ then obtain
  \[
  \cal{R}^{P}_{\rm upper}(\pi)\ge \gamma_{T}^{P}(\pi,f)
  -I_\tau^{P}(\pi,L[f])= I_{\tau^*}^{P}(\pi,e^{L[f]}-1).
  \]
  Hence by the definition~\eqref{defrate2} of $\cal{R}^{P}_{\rm lower}(\pi)$,
  \[
  I_{\tau^*}^{P}(\pi,e^{L[f]}-1)\le   \cal{R}^{P}_{\rm upper}(\pi)\le
    \cal{R}^{P}_{\rm lower}(\pi)\le   I_{\tau^*}^{P}(\pi,e^{L[f]}-1),
  \]
  the proof completes.
\end{proof}

\begin{lemma}
  [Exponential tightness]\label{exptight}
 Under the assumptions in Theorem~\ref{main:ldp}, there exists a sequence of subsets $(\cal{E}_q)_{q\in\N}$ of $D_T(\cal{M}^+(X))$ such that for all $q\in \N$, the closure of $\cal{E}_q$ is compact and
  \[
\limsup_{h\to 0}h\log \P^{h,P}(\mu\notin \cal{E}_q)\le -q.
\]
Moreover, for any $q>0$, there exists $n_0$, such that for all $n>n_0$, if $P$ satisfies $(\oplus)$ condition, then
\[
\limsup_{h\to 0}h\log \P^{h,P}(\mu\notin \tilde{\cal{D}}^n_T)\le -q,
\]
and if $P$ satisfies $(\otimes)$ condition, then
\[
\limsup_{h\to 0}h\log \P^{h,P}(\mu\notin \cal{D}^n_T)\le -q,
\]
where $\cal{D}^n_T,\tilde{\cal{D}}^n_T$ are the sets defined in~\eqref{dn} and \eqref{dn2}.
\end{lemma}
\begin{proof}
We define two decreasing sequences $(\delta_r)_{r\in\N}$ and $(\eps_r)_{r\in\N}$, for $r\in\N$, let
  \begin{multline*}
\delta_r:=\min\bigg\{\alpha\left(\left(e^{4r^2 k}+1\right)\|P\|_{(1+E)^{\tilde{\otimes}},\infty}e^{c_0}\right)^{-1},\\
\left(2\left(e^{4r^2 k}+1\right)\|P\|_{(1+E)^{\tilde{\otimes}},\infty}e^{c_0}\left(1+\frac{1}{\alpha}
\sup_h h\log \E \exp\left\{h^{-1}\alpha\left(\langle 1+E,\mu_0^h\rangle\right)^k\right\}
   \right)\right)^{-1}\bigg\}
\end{multline*}
where $\alpha$ is the constant in relation~\eqref{initial:log} and
\[
\eps_r=\frac{1}{r}.
\]
Let $(f_j)_{j\in\N}$ be a sequence of functions dense in $C_b(X)$ and such that $\|f_j\|_{\infty}\le j$. We define a sequence of sets
\[
\cal{F}_q:=\left\{\mu\in D_T(\cal{M}^+(X))\bigg|\inf_{t_i}\max_i\sup_{t_i\le s<t\le t_{i+1}}\sum_{j=1}^q\frac{1}{4^j}\left|\langle \mu_t-\mu_s,f_j\rangle\right|\le \eps_q\right\}
\]
where the infimum extends over all finite partitions $(t_j)$ of $[0,T]$ satisfying $0=t_0<t_1<\dots<t_n=T$, $t_{i+1}-t_i\ge \delta_q$. Then by  Chapter~3 of Ethier and Kurtz~\cite{ethier} , the closure sets of 
\[
\cal{E}_q:=\bigcap_{r\ge q}\cal{F}_r
\]
are compact in $D_T(\cal{M}^+(X))$ for all $q\in\N$.

We fix any function $f\in C_b(X)$, for any $\lambda>0$,
any $0\le s\le t\le T$, $|s-t|\le \delta$, then by Lemma~\ref{cdtc0},  under $\P_{\mu_0^h}^P$,
\begin{multline*}
  \int_s^t\int_{\cal{J}}\left(e^{\langle \lambda f,\delta^+_\bfy-\delta^+_\bfx\rangle}-1\right)\kappa^{h,P}[ \mu](\diff u\diff  \bfx\diff\bfy)\\
  \le \delta \left(e^{2\lambda k \|f\|_{\infty}}+1\right)\|P\|_{(1+E)^{\tilde{\otimes}},\infty}\sum_{\ell=1}^k\frac{1}{\ell!}\left(c_0\langle 1+E,\mu_0^h\rangle\right)^\ell\\
  \le  \delta \left(e^{2\lambda k \|f\|_{\infty}}+1\right)\|P\|_{(1+E)^{\tilde{\otimes}},\infty}e^{c_0}\left(1+\left(\langle 1+E,\mu_0^h\rangle\right)^k\right)\\
  \end{multline*}
By using the $f$-perturbed exponential martingale $\cal{M}_{s,t}^{h,(P^f|P)}$ (Definition~\ref{defmg}), for any $\eps>0$, we get
\begin{multline*}
 \P^{h,P}\left(\sup_{|s-t|<\delta}\langle f,\mu_t\rangle-\langle f,\mu_s\rangle\ge \eps\right)\\
 \le \P^{h,P}\bigg(\sup_{|s-t|<\delta}\cal{M}_{s,t}^{h,(P^f|P)}[\mu]
 \ge
  \exp\left\{h^{-1}\lambda\eps
  -h^{-1}\delta C_1(\lambda,f)\left(1+\left(\langle 1+E,\mu_0^h\rangle\right)^k\right)\right\}\bigg)\\
  \le \E \exp\left\{-h^{-1}\lambda \eps +h^{-1}\delta C_1(\lambda,f)\left(1+\left(\langle 1+E,\mu_0^h\rangle\right)^k\right)\right\},
\end{multline*}
where
\[
C_1(\lambda,f)=\left(e^{2\lambda k \|f\|_{\infty}}+1\right)\|P\|_{(1+E)^{\tilde{\otimes}},\infty}e^{c_0}.
\]
We get the same upper bound for $-f$. Thus
\begin{multline*}
  h\log  \P^{h,P}\left(\sup_{|s-t|<\delta}\left|\langle f,\mu_t\rangle-\langle f,\mu_s\rangle\right|\ge \eps\right)\\
  \le
  -\lambda \eps+\delta C_1(\lambda,f)
  +h\log \E \exp\left\{h^{-1}\delta C_1(\lambda,f)\left(\langle 1+E,\mu_0^h\rangle\right)^k\right\}\\
   \le
  -\lambda \eps+\delta C_1(\lambda,f)
+\frac{\delta C_1(\lambda,f)}{\alpha}h\log \E \exp\left\{h^{-1}\alpha\left(\langle 1+E,\mu_0^h\rangle\right)^k\right\}
\end{multline*}
if $\delta C_1(\lambda,f)<\alpha$.

For each $r$ and $j$, we chose $\lambda=r^2j$, by using $\|j f_j2^{-j}\|\le j^{2}2^{-j}\le 2$, then we have
\[
C_1(r^2j,f_j2^{-j})\le \left(e^{4r^2 k}+1\right)\|P\|_{(1+E)^{\tilde{\otimes}},\infty}e^{c_0}.
\]
Hence, we get
\begin{multline*}
   h\log  \P^{h,P}\left(\sup_{|s-t|\le 2\delta_r,s,t\in[0,T]}\left|\langle \mu_t-\mu_s,f_j\rangle\right|\ge 2^j \eps_r \right)\\
   \le -r^2j\eps_r +2\delta_rC_1(r^2j,f_j2^{-j})\left(1+\frac{1}{\alpha}
h\log \E \exp\left\{h^{-1}\alpha\left(\langle 1+E,\mu_0^h\rangle\right)^k\right\}
   \right)\\
  \le (1- jr).
\end{multline*}

Therefore,
\begin{multline*}
  \P^{h,P}\left(\overline{\cal{E}_q}^c\right)\le \sum_{r\ge q}\P^{h,P}\left({\cal{F}_r}^c\right)\\\le
  \sum_{r\ge q}  \P^{h,P}\left(\sup_{|s-t|\le 2\delta_r,s,t\in[0,T]}\sum_{j=1}^\infty\frac{1}{4^j}\left|\langle \mu_t-\mu_s,f_j\rangle\right|\ge \left(\sum_{j\ge 1}\frac{1}{2^j}\right) \eps_r \right)\\
  \le  \sum_{r\ge q} \sum_{j\ge 1} \P^{h,P}\left(\sup_{|s-t|\le 2\delta_r,s,t\in[0,T]}\left|\langle \mu_t-\mu_s,f_j\rangle\right|\ge 2^j \eps_r \right)\\
\le \sum_{r\ge q} e^{1/h}\sum_{j\ge 1}
e^{- jr/h}\le  \sum_{r\ge q} e^{1/h}\frac{e^{-r/h}}{1-e^{-r/h}}\le \frac{e^{-q/h+1/h}}{(1-e^{-q/h})^2}.
\end{multline*}
In conclusion,
\[
\limsup_{h\to 0} h\log \P^{h,P}\left(\overline{\cal{E}_q}^c\right)\le -q+1.
\]

Now suppose that $P$ satisfies condition~$(\oplus)$ then
\[
\sup_{0\le t\le T}\langle 1,\mu_t\rangle\le C_\#\qquad \P^{h,P}-a.s.
\]
where $C_\#$ is the constant in~\eqref{1bound}. Hence for any $n$,
\[
\P^{h,P}(\mu\notin \tilde{\cal{D}}^n_T)=
\P^{h,P}\left(\int_0^T\langle 1+E^\beta,\mu_t \rangle \diff t>n\right).
\]
By Jensen's inequality and  Lemma~\ref{boundness},   we have
\begin{multline*}
  \P^{h,P}\left(\int_0^T\langle 1+E^\beta,\mu_t \rangle\diff t>n\right)\le
  e^{-h^{-1}\gamma n/T}\E^{h,P}\exp\left\{h^{-1}\gamma \frac{1}{T}\int_0^T\langle 1+E^\beta,\mu_t \rangle\diff t\right\}\\
  \le   e^{-h^{-1}\gamma n/T}\frac{1}{T}\int_0^T\E^{h,P}\exp\left\{h^{-1}\gamma \langle 1+E^\beta,\mu_t \rangle\right\}\diff t\\
  \le e^{-h^{-1}\gamma n/T}\sup_{0\le t\le T}\E^{h,P}\exp\left\{h^{-1}\gamma \langle 1+E^\beta,\mu_t \rangle\right\}.
\end{multline*}
Therefore, by Lemma~\ref{hlogexp}, we have
\begin{multline*}
  \limsup_{h\to 0}h\log \P^{h,P}\left(\int_0^T\langle 1+E^\beta,\mu_t \rangle\diff t>n\right)\\
  \le
-\frac{\gamma n}{T}+\frac{c'\gamma}{\alpha}\limsup_{h\to 0}h\log \E \exp\left\{\alpha h^{-1}\langle 1+E^\beta,\mu^h_0 \rangle\right\}
\end{multline*}
that tends to $-\infty$ as $n$ goes to $\infty$.

When $P$ satisfies condition $(\otimes)$, by relation~\eqref{1+emoment}, we have
\[
\P^{h,P}(\mu\notin {\cal{D}}^n_T)=
\P^{h,P}\left(\sup_{0\le t\le T}\langle 1+E,\mu_t \rangle >n\right)
\le \P^{h,P}\left(\langle 1+E,\mu_t \rangle >n/c_0\right).
\]
The proof completes by using Lemma~\ref{hlogexp} again.
\end{proof}

\begin{proof}[Proof of Theorem~\ref{main:ldp}]

  \proofstep{(Upper bound)}
  For any $f\in C^{1,0}_b([0,T]{\times} X)$, any measurable set $A$ in $D([0,T];\cal{M}^+(X))$, any $\mu_0^h\in \cal{M}^+_{h\delta }(X)$ the Radon-Nikodym~\eqref{rndf} of the $f$-perturbed probability $\P^{h,P^f}_{\mu_0^h}$ with respect to the probability $\P^{h,P}_{\mu_0^h}$ gives
  \begin{multline*}
    \P^{h,P}_{\mu_0^h}(A)
    =\E^{h,P^f}_{\mu_0^h}\left(\frac{\diff \P^{h,P}_{\mu_0^h}}{\diff \P^{h,P^f}_{\mu_0^h}}[\mu]\ind{\mu\in A}\right)\\
    =\E^{h,P^f}_{\mu_0^h}\Bigg(
\exp\bigg\{
  - h^{-1}\left(\langle f_T,\mu_T\rangle-\langle f_0,\mu_0\rangle-
    \int_0^T\langle \partial_s f_s,\mu_s\rangle\diff s-\int_0^T\Lambda^{h,P}_s\langle f_s,\mu_s\rangle\diff s\right)
    \\+h^{-1}\int_0^T\int_{\cal{J}}\tau\left(\langle f_s,\delta^+_\bfy-\delta^+_\bfx\rangle\right)\kappa^{h,P}[\mu](\diff s\diff\bfx\diff\bfy)\bigg\}
    \ind{\mu\in A}
    \Bigg)
  \end{multline*}
  For the chaotic initial condition, by using notations~\eqref{rategamma} and~\eqref{ratei}, we get the following upper bound
   \begin{multline*}
    h   \log\left( \P^{h,P}(A)\right)
    \le -\inf_{\mu\in A}\left(
    \gamma^P_T(\mu,f)-I_\tau^P(\mu,L[f])
+\langle \phi,\mu_0\rangle-\log \langle e^{\phi},\nu \rangle
    \right)\\
+h\log \left(\E^{h,P^f}
\exp\bigg\{
h^{-1}\left|\left\langle e^{L[f]}-1,\left(\kappa^{h,P}[\mu]-\cal{K}^P[\mu]\right)\right\rangle\right|\bigg\}
\right)
   \end{multline*}
holds for all $f\in C^{1,0}_b([0,T]{\times} X)$ and $\phi\in C_b(X)$.

  By Proposition~\ref{diffexp}, we get
   \begin{multline*}
     \limsup_{h\to 0}    h   \log\left( \P^{h,P}(A)\right)\\
     \le -\sup_{\substack{f\in C^{1,0}_b([0,T]{\times} X)\\ \phi\in C_b(X)}}\inf_{\mu\in A}\Bigg(
\gamma^P_T(\mu,f)-I_\tau^P(\mu,L[f])+\langle \phi,\mu_0\rangle-\log \langle e^{\phi},\nu \rangle
    \Bigg).
   \end{multline*}
Recall Proposition~\ref{ctnfk}, $I_\tau^P(\cdot,L[f])$ is continuous on $\cal{D}^n_T$ (\emph{resp.} on $\tilde{\cal{D}}^n_T$) for $(\otimes)$ (\emph{resp.} $(\oplus)$)  type kernel. By  following the standard method (See Appendix~2 Lemma~3.2 in Kipnis and Landim~\cite{kipnis1999sli} for example), in the $(\otimes)$ case, we can exchange the infimum with the maximum in the upper bound for the set $A\cap \cal{D}^n_T$ when $A$ is a compact set and get
     \[
\limsup_{h\to 0}    h   \log\left( \P^{h,P}(A\cap \cal{D}^n_T)\right)\le -\inf_{\mu\in A\cap \cal{D}^n_T}\cal{I}^P[\mu]\le  -\inf_{\mu}\cal{I}^P[\mu].
\]
By Lemma~\ref{exptight}, we can choose $n$ sufficient large, such that
\[
\limsup_{h\to 0} h\log \P^{h,P}(\mu\notin  D_T^{n})\le  -\inf_{\mu}\cal{I}^P[\mu].
\]
By using 
 \begin{multline*}
    h   \log\left( \P^{h,P}(A)\right)\le 
    h   \log\left(\P^{h,P}(A\cap (\cal{D}^{n}_T))+\P^{h,P}(\mu\notin  D_T^{n})\right)\\
    \le h\log 2+\max\left\{h   \log\left(\P^{h,P}(A\cap (\cal{D}^{n}_T))\right),
h\log \P^{h,P}(\mu\notin  D_T^{n})
    \right\},
\end{multline*}
we have
\[
    \limsup_{h\to 0}    h   \log\left( \P^{h,P}(A)\right)
\le  -\inf_{\mu\in A}\cal{I}^P[\mu].
\]
 Combine with the exponential tightness result, the upper bound holds for all closed set.   Similar proof holds for the deterministic initial condition and also the $(\oplus)$ kernels.

\proofstep{(Lower bound)}
Here we prove the lower bound for $(\otimes)$ type kernels ($\beta=1$). For $(\oplus)$ case, the proof holds when we replace $\cal{D}^{n}_T$ by $\tilde{\cal{D}}^{n}_T$ defined in~\eqref{dn2}.

For any open set $\cal{O}$ of $D_T(\cal{M}^+(X))$,
if $\sigma\in \cal{O}$ or $\cal{O}\cap H^P_0=\emptyset$, then there is nothing to prove. We suppose that $\cal{O}\cap H^P_0\neq\emptyset$. Then by definition of $H^P_0$,  there exists a measure $\nu'\in\cal{M}^+(E)$ such that $\langle 1+E^\beta,\nu'\rangle<\infty$, a non-negative bounded function $\eta\in\cal{L}^\infty([0,T]{\times}\cal{J})$ and a constant $\eps>0$, such that the open neighborhood $B_{\eps}(\sigma^\eta)\subset \cal{O}\cap H^P_0$ where $\sigma^\eta\in H_0^{P}[\nu']$ with   Radon-Nikodym derivative $\eta$ and
\[
B_{\eps}(\sigma^\eta):=\left\{\mu\in D_{T}(\cal{M}^+(E))\bigg|\sup_{0\le t\le T}\|\mu_t-\sigma^\eta_t\|<\eps\right\}.
\]
For all $n\in\N$, let
$B_{\eps,n}(\sigma^\eta):=B_{\eps}(\sigma^\eta)\cap \cal{D}^{n}_T$ where $\cal{D}^{n}_T$ defined in~\eqref{dn}. Then by Proposition~\ref{ctnfk}, for all $n$ large enough, for any $\delta>0$, we can choose $\eps$ small enough, such that
\[
\sup_{\mu\in B_{\eps,n}(\sigma^\eta)}\left|\langle \tau^*\left(\eta-1\right),\cal{K}^{P}[\mu]\rangle-
\langle \tau^*\left(\eta-1\right),\cal{K}^{P}[\sigma^{\eta}]\rangle\right|\le \delta.
\]

For such a $\delta$, we define a subset of $D_T(\cal{M}^+(X))$ by
\begin{multline*}
  V^{h,\eta}_\delta:=
  \Bigg\{\mu\in D_T(\cal{M}^+(X))\Bigg|\left|\langle \tau^*(\eta-1),\left(\cal{K}^P[\mu]-\kappa^{h,P}[ \mu]\right)\rangle\right|\le \delta,\\
  \left|\langle \log \left(\eta\right),\tilde{\cal{N}}^{h,\eta  P}[\mu]\rangle\right|\le h^{-1}\delta
  \Bigg\},
\end{multline*}
where $\tilde{\cal{N}}^{h,\eta P}[\mu]$ is defined in~\eqref{defitildeN}. 

For any $\pi\in\cal{M}^{+}_{h\delta}(X)$, from Proposition~\ref{coupling}, we see that on the event $B_{\eps,n}(\sigma^\eta)\cap V^{h,\eta}_\delta$
\begin{multline*}
\frac{\diff \P^{h, P}_{\pi}}{\diff\P^{h,\eta P}_{\pi}}[\mu]=\exp\big\{-
\langle \log \left(\eta\right),\tilde{\cal{N}}^{h,\eta  P}[\mu]\rangle+h^{-1}\langle \tau^*(\eta-1),\left(\cal{K}^P[\mu]-\kappa^{h,P}[ \mu]\right)\rangle
\\-h^{-1}\langle \tau^*\left(\eta-1\right),\cal{K}^{P}[\mu]\rangle\big\}\\
\ge \exp\left\{-3h^{-1}\delta-h^{-1}\langle \tau^*\left(\eta-1\right),\cal{K}^{P}[\sigma^\eta]\rangle\right\}.
\end{multline*}

In the case of chaotic initial condition, let  $\tilde{\P}^{h,\eta P}$ be the martingale solution with kernel $\eta P$ and the chaotic initial particles with common distribution $\nu'$, then
\begin{multline}\label{chaopf}
  \P^{h,P}(B_{\eps,n}(\sigma^\eta)\cap V^{h,\eta}_\delta)\\
  \ge
\exp\left\{-4h^{-1}\delta-I_0(\nu')-h^{-1}\langle \tau^*\left(\eta-1\right),\cal{K}^{P}[\sigma^{\eta}]\rangle\right\}
\tilde{\P}^{h,\eta P}\left(B_{\eps,n}(\sigma^\eta)\cap V^{h,\eta}_\delta\right).
\end{multline}

By using   Proposition~\ref{diffexp}, relation~\eqref{kcvg}, as $h\to 0$, we have
\begin{multline*}
  \tilde{\P}^{h,\eta P}\left(\left|\langle \tau^*(\eta-1),\left(\cal{K}^P[\mu]-\kappa^{h,P}[ \mu]\right)\rangle\right|> \delta\right)\\
  \le\frac{1}{\delta}\tilde{\E}^{h,\eta P}
\left|\langle \tau^*(\eta-1),\left(\cal{K}^P[\mu]-\kappa^{h,P}[ \mu]\right)\rangle\right|\to 0.
\end{multline*}
On the other hand, according to Proposition~\ref{cox}, under probability $\tilde{\P}^{h,\eta P}$, the process
\[
\left(\sqrt{h}\int_0^t\int_{\cal{J}}\log \left(\eta(s,\bfx,\bfy)\right)\tilde{\cal{N}}^{h,\eta  P}[\mu](\diff s\diff\bfx\diff\bfy )\right)
\]
is a local martingale whose predictable process is given by
\begin{multline*}
  \left(\int_0^t\int_{\cal{J}}\left(\log \eta(s,\bfx,\bfy)\right)^2\kappa^{h,\eta P}[ \mu](\diff s\diff\bfx\diff\bfy )\right)\\
  =
\left(\int_0^t\int_{\cal{J}}\eta\left(\log \eta\right)^2(s,\bfx,\bfy)\kappa^{h, P}[ \mu](\diff s\diff\bfx\diff\bfy )\right).
\end{multline*}
Note that 
the function $f(a)=a(\log a)^2$ is bounded on the interval $[0,\|\eta\|_\infty]$, thus $\eta(\log (\eta))^2\in\cal{L}^\infty([0,T]\times \cal{J})$. According to relation~\eqref{bd}, we have
\[
\sup_h\tilde{\E}^{h,\eta P}\langle \eta(\log (\eta))^2,\kappa^{h, P}[ \mu]\rangle<\infty.
\]
Therefore, we have
\[
\tilde{\P}^{h,\eta P}\left(  \left|\langle \log \left(\eta\right),\tilde{\cal{N}}^{h,\eta  P}[\mu]\rangle\right|>h^{-1}\delta\right)\le \frac{\sqrt{h}}{\delta}\sup_{h}\tilde{\E}^{h,\eta P}\langle \eta(\log (\eta))^2,\kappa^{h, P}[ \mu]\rangle,
\]
that converges to $0$ as $h\to 0$.
Combining with the law of large numbers and Lemma~\ref{exptight},  for  $n$ sufficient large, we have 
\[
  \lim_{h\to 0} \tilde{\P}^{h,\eta P}(B_{\eps,n}(\sigma^\eta)\cap V^{h,\eta}_\delta)\ge 1
  - \lim_{h\to 0} \tilde{\P}^{h,\eta P}((\cal{D}_T^n)^c)=1.
\]
In conclusion,
\begin{multline*}
\liminf_{h\to 0}h\log\left(
 \P^{h,P}(\cal{O})
 \right)\ge
 \liminf_{h\to 0}h\log\left(
 \P^{h,P}(B_{\eps,n}(\sigma^\eta)\cap V^{h,\eta}_\delta)
\right)
\\\ge \left(-4\delta-I_0(\nu')-\langle \tau^*\left(\eta-1\right),\cal{K}^{P}[\sigma^{\eta}]\rangle\right),
\end{multline*}
holds for all $\delta>0$. By letting $\delta\to 0$, we obtain the lower bound,
\[
\liminf_{h\to 0}h\log\left(
 \P^{h,P}(\cal{O})
 \right)\ge -I_0(\nu')-I_{\tau^*}^P(\sigma^\eta,\eta-1).
\]
Note that the  proof above holds for all bounded measurable non-negative $\eta'$ such that
\[
\gamma^{\eta'P}_t(\sigma^\eta,g)=0\qquad \forall a.e.-t\in[0,T],\forall g\in \cal{C}_b^{1,0}([0,T]{\times} X),
\]
thus
\[
\liminf_{h\to 0}h\log\left(
 \P^{h,P}(\cal{O})
 \right)\ge -I_0(\nu')-\cal{R}^{P}_{\rm lower}(\sigma^\eta).
 \]
 With a trivial modification in relation~\eqref{chaopf}, the proof also holds for the deterministic initial condition.
\end{proof}

\section{Proofs of the good ``gelling'' solution}\label{sec:cm}

\begin{proof}[Proof of Theorem~\ref{gel}]
Suppose that both $P$ and $\eta P$ are $(1+E)^{\tilde{\otimes}}$ bounded.  
  Let $\sigma^{n,\eta}_t$ be the solution of equation~\eqref{llneq} with a bounded cut-off kernel $\eta_n P$ (See Definition~\ref{kernel2}) starting from $\nu$ such that $\langle 1+ E^2,\nu\rangle<\infty$.  We then have
\begin{equation*}
    \langle 1+ E^2,\sigma^{n,\eta}_t\rangle\le     \langle 1+ E^2,\nu\rangle
    +\int_0^t\int_{\cal{J}}\langle 1+E^2,\delta^+_{\bfy}-\delta^+_{\bfx}\rangle \eta_n \cal{K}^P[\sigma^{n,\eta}](\tridiff)
\end{equation*}
Note that the increasment of the number of particles is no more than $k$ and  only happens when the fragmentation happens, \emph{i.e.}, $\langle 1,\delta^+_{\bfy}-\delta^+_{\bfx}\rangle\ge 0 $ only if $|\bfx|=1$. Hence we have
\begin{multline*}
  \int_0^t\int_{\cal{J}}\langle 1,\delta^+_{\bfy}-\delta^+_{\bfx}\rangle \eta_n \cal{K}^P[\sigma^{n,\eta}](\tridiff)\le k\|\eta P\|_{(1+E)^{\otimes},\infty}\int_{0}^t\langle 1+ E,\sigma^{n,\eta}_s\rangle\diff s\\
  \le  k\|\eta P\|_{(1+E)^{\otimes},\infty}tc_0\langle 1+ E,\nu\rangle,
\end{multline*}
where $c_0=1$ or $(1+1/\eps_0)$ coming from   the ``no dust'' condition.
By using the fact kernel $P$ is $E$ non-increasing and $\langle E^2,\delta^+_\bfy\rangle\le \langle E,\delta^+_\bfy\rangle^2$, we have
\begin{multline*}
  \int_0^t\int_{\cal{J}}\langle E^2,\delta^+_{\bfy}-\delta^+_{\bfx}\rangle \eta_n \cal{K}^P[\sigma^{n,\eta}](\tridiff)\\
  \le \|\eta P\|_{(1+E)^{\otimes},\infty}
  \sum_{\ell=1}^k\int_0^t\int_{SX^\ell}\left(\left(\langle E,\delta^+_\bfx\rangle\right)^2-\langle E^2,\delta^+_{\bfx}\rangle\right) (1+E)^{\tilde{\otimes}}(\bfx) (\sigma^{n,\eta}_s)^{\tilde{\otimes}\ell}(\diff\bfx)\diff s\\
  \le \|\eta P\|_{(1+E)^{\otimes},\infty}
  \sum_{\ell=2}^k\frac{1}{(\ell-2)!}\int_0^t\langle (1+E)^2,\sigma^{n,\eta}_s\rangle^2\langle (1+E),\sigma^{n,\eta}_s\rangle^{\ell-2}\diff s\\
  \le  2\|\eta P\|_{(1+E)^{\otimes},\infty}
\sum_{\ell=2}^k\frac{1}{(\ell-2)!}(c_0 \langle 1+ E,\nu\rangle)^{\ell-2}\int_0^t\langle 1+E^2,\sigma^{n,\eta}_s\rangle^2\diff s.
\end{multline*}
Let
\begin{multline*}
  C(\|\eta P\|_{(1+E)^{\otimes},\infty},\nu)
  :=\\
  \max\left\{  k\|\eta P\|_{(1+E)^{\otimes},\infty}c_0\langle 1+ E,\nu\rangle,2\|\eta P\|_{(1+E)^{\otimes},\infty}
\sum_{\ell=2}^k\frac{1}{(\ell-2)!}(c_0 \langle 1+ E,\nu\rangle)^{\ell-2}\right\}
\end{multline*}
then
\begin{multline*}
  \frac{\diff}{\diff t}  \langle 1+ E^2,\sigma^{n,\eta}_t\rangle\le C(\|\eta P\|_{(1+E)^{\otimes},\infty},\nu)\left(1+\langle 1+ E^2,\sigma^{n,\eta}_t\rangle^2\right)\\
  \le C(\|\eta P\|_{(1+E)^{\otimes},\infty},\nu)\left(1+\langle 1+ E^2,\sigma^{n,\eta}_t\rangle\right)^2
\end{multline*}
which implies
\[
 \langle 1+ E^2,\sigma^{n,\eta}_t\rangle\le \frac{1}{\frac{1}{ \langle 1+ E^2,\nu\rangle+1}-C(\|\eta P\|_{(1+E)^{\otimes},\infty},\nu)t}-1.
 \]
Therefore, if 
 \[
T<T_*:=\frac{1}{( \langle 1+ E^2,\nu\rangle+1)C(\|\eta P\|_{(1+E)^{\otimes},\infty},\nu)},
\]
we have
\[
\sup_n\sup_{t\le T} \langle 1+ E^2,\sigma^{n,\eta}_t\rangle<\infty.
\]
    
Therefore,
\begin{multline*}
 \sup_{n}\sup_{0\le t\le T} \left|\int_0^t\int_{\cal{J}}\langle g,\delta^+_{\bfy}-\delta^+_{\bfx}\rangle \eta_n \cal{K}^P[\sigma^{n,\eta}](\tridiff)\ind{\max\{E(\bfx)\}\ge K}\right|\\
  \le \frac{2k\|g\|_\infty\|\eta P\|_{(1+E)^{\tilde{\otimes}},\infty}T}{K}\left(\sum_{\ell=1}^k\frac{1}{\ell!}\left(c_0\langle 1+E,\nu\rangle\right)^{\ell-1}\right)\left(\langle E,\nu\rangle +\sup_n\sup_{0\le s\le T}\langle E^2,\sigma^{n,\eta}_s\rangle\right)\\
  \to 0\qquad a.s. ~~~K\to\infty.
\end{multline*}
On the other hand,
\begin{multline*}
 \max_{1\le \ell\le k,0\le m\le d_\ell} \sup_{s\in[0,T]}\sup_{\bfx\in SX^\ell}\left|\int_{SX^m}\langle g,\delta^+_{\bfy}-\delta^+_{\bfx}\rangle\eta_nP(s,\bfx,\diff \bfy)\ind{\max\{E(\bfx)\}<K}\right|\\
  \le 2k\|g\|_\infty\|\eta P\|_{(1+E)^{\tilde{\otimes}},\infty}k(1+K)^k<\infty.
\end{multline*}
The tightness of $\langle g,\sigma^{n,\eta}_t\rangle$ for any bounded function $g$ follows a  standard proof. Let $\sigma^\eta$ be a limit, then  as $n\to \infty$, 
\begin{multline*}
  \int_{0}^t\int_{\cal{J}}\langle g,\delta^+_{\bfy}-\delta^+_{\bfx}\rangle\eta_n\cal{K}^{P}[\sigma^{n,\eta}](\tridiff)\ind{\max\{E(\bfx)\}<K}\\
  \to
\int_{0}^t\int_{\cal{J}}\langle g,\delta^+_{\bfy}-\delta^+_{\bfx}\rangle\eta\cal{K}^{P}[\sigma^{\eta}](\tridiff)\ind{\max\{E(\bfx)\}<K}.
\end{multline*}
The proof completes by letting $K\to \infty$.
When we take $\eta=1$, the above results reduce to the existence of the solution of equation~\eqref{llneq} for all $(1+E)^{\tilde{\otimes}}$ bounded kernel $P$.

\end{proof}

\appendix
\section{Conjugation of integral functionals on a subspace of $L^\infty$}\label{sec:app}

Let $T$ be a measure space with a $\sigma$-finite measure $\Lambda(\diff t)$.
Consider the space $L^\infty(T)$ of bounded measurable functions on $T$.

\begin{theorem}\label{thm:ro}
  (Theorem~4 in Rockafellar~\cite{rock1})
  Let $T$ be of finite measure. Let $f(t,x)$ be a finite convex function of $x$ for each $t$ and a bounded measurable function of $t$ for each $x$.

  Then, $I_f: L^\infty(T)\to \R$
  \[
I_f(u)=\int_T f(t,u(t))\Lambda(\diff t),\qquad u\in L^\infty(T),
\]
is a well-defined finite convex function and is everywhere continuous with respect to the uniform norm.

Moreover, the conjugate $(I_f)^*$ of $I_f$ on $L^\infty(T)^*$ is given by, for any $v\in (L^\infty(T))^*$
\begin{displaymath}
(I_f)^*(v) = \left\{ \begin{array}{ll}
I_{f^*}(u^*), & \textrm{if }\exists~u^*\in L^1(T),~s.t.~v(u)=\int_T\langle u(t),u^*(t) \rangle\Lambda(\diff t);\\
+\infty, & \textrm{otherwise.}
\end{array} \right.
\end{displaymath}

The convex set
\[
\left\{u^*\in L^1(T)\bigg|(I_{f^*})(u^*)+\langle a, u^*\rangle_{\Lambda}+\alpha\le 0\right\}
  \]
  is weak* compact in the $\sigma(L^1(T),L^\infty(T))$ topology for any $a\in L^\infty(T)$ and $\alpha\in \R$.
\end{theorem}
Let $D$ be a subspace of $L^\infty(T)$ supplied with a locally convex topology at least as strong as the uniform norm topology. Let $D^*$ be the space of continuous linear functionals on $D$.

\begin{corollary}\label{co:ro}
  (Corollary~2 in Rockafellar~\cite{rock1})
  Suppose that no nonzero linear functional on $L^\infty(T)$ of the form
  \[
  u\to \int_T\langle u(t), u^*(t)\rangle\diff t,\qquad u^*\in L^1(T)
  \]
  vanishes throughout $D$. Then, under the hypothesis of Theorem~\ref{thm:ro}, $I_f$ is a continuous finite convex function on $D$, and the conjugate $(I_f)^*$ of $(I_f)$ on $D^*$ is given by $I_{f^*}$, in the sense that if $v\in D^*$ corresponds to some $u^*\in L^1(T)$ as above one have $(I_f)^*(v)=I_{f^*}(u^*)$ whereas other-wise $(I_f)^*(v)=+\infty$.
\end{corollary}

Define the equivalence relation $\sim$ on $L^1(T)$: $u^*_1\sim u^*_2$ iff
$u_1^*-u_2^*\in D_0$ where
\[
D_0:=\left\{u^*\in L^1(T)\bigg|\langle u,u^*\rangle_{\Lambda}=0,~\forall u\in D\right\}.
\]
\begin{lemma}\label{compact}
  
The convex sets
  \[
C_\alpha:=\left\{[u^*]\in L^1(T)/D_0\bigg|\inf_{x\in L^1(T),x\sim u^{*}}(I_{f^*})(x)\le \alpha\right\},\qquad \forall~\alpha\in \R
\]
are compact in the quotient topology $\sigma(L^1(T), L^{\infty}(T))/D_0$ on $L^1(T)/D_0$.
\end{lemma}
\begin{proof}
  Let $p$ be the surjective map from $L^1(T)$ to $L^1(T)/D_0$: $p(u^*)=[u^*]$. For any $\alpha \in \R$, let
  \[c_\alpha:=\left\{u^*\in L^1(T)\bigg|I_{f^*}(u^*)\le \alpha\right\}.\]
It is easy to see $p(c_\alpha)\subset C_\alpha.$
Since $I_{f^*}$ is a proper l.s.c. convex functional on $L^1(T)$ and $\{x\in L^1(T)|x\sim u^*\}$ is a norm-closed convex set in $L^1(T)$, then 
\[
p(c_\alpha)= C_\alpha.
\]
Let $\{U_i\}_{i\in I}$ be any open cover of $C_\alpha$ in the quotient topology $\sigma(L^1(T), L^{\infty}(T))/D_0$, $\{p^{-1}(U_i)\}_{i\in I}$ is an open cover of $c_\alpha$ in the weak topology $\sigma(L^1(T), L^{\infty}(T))$. By Theorem~\ref{thm:ro}, $c_\alpha$ is compact set, then there is a finite cover $\{p^{-1}(U_i)\}_{i\in I_0}$ of $c_\alpha$. Then $\{U_i\}_{i\in I_0}=\{p(p^{-1}(U_i))\}$ is a finite open cover of $C_\alpha=p(c_\alpha)$. Therefore, $C_\alpha$ is compact.
\end{proof}

\begin{corollary}\label{subconj}
Under the hypothesis of Theorem~\ref{thm:ro}, $I_f$ is a continuous finite convex function on $D$, and the conjugate $(I_f)^*$ of $I_f$ on $D^*$ is given by, for any $v\in D^*$
  \[
(I_f)^*(v) =\inf_{u^*\in\cal{O}_v}I_{f^*}(u^*)
  \]
  where
  \[
\cal{O}_v:=\left\{u^*\in L^1(T)\bigg| v(u)=\int_T\langle u(t),u^*(t)\rangle\Lambda(\diff t),~\forall~u\in D\right\},
\]
and with convention $\inf_{\emptyset}=+\infty$.
\end{corollary}
\begin{proof}

  Define functional $J: D^*\mapsto (-\infty,+\infty]$,
  \[
J(v)=\inf_{u^*\in\cal{O}_v}I_{f^*}(u^*).
\]
For any $v_1,v_2\in D^*$, $\lambda\in[0,1]$, using the fact
\[
\left\{\lambda u^*_1+(1-\lambda)u^*_2\bigg|u^*_1\in\cal{O}_{v_1},u^*_2\in\cal{O}_{v_2}\right\}\subset\cal{O}_{\lambda v_1+(1-\lambda)v_2},
\]
one has
\begin{multline*}
  J(\lambda v_1+(1-\lambda)v_2)=\inf_{\cal{O}_{\lambda v_1+(1-\lambda)v_2}}I_{f^*}(u^*)\\
  \le
  \inf_{u_1^*\in\cal{O}_{v_1},u_2^*\in\cal{O}_{v_2}}I_{f^*}(\lambda u^*_1+(1-\lambda)u^*_2)
  \le
  \lambda J(v_1)+(1-\lambda)J(v_2).
\end{multline*}

For any $u\in D$, the conjugate
\begin{multline*}
  (J)^*(u)=\sup_{v\in D^*}\{v(u)-J(v)\}=\sup_{v\in D^*}\{v(u)-\inf_{u^*\in\cal{O}_v}I_{f^*}(u^*)\}\\
  =\sup_{v\in D^*}\sup_{u^*\in\cal{O}_v}\{\langle u,u^*\rangle_{\Lambda}-I_{f^*}(u^*)\}\\
  \le \sup_{v\in D^*}\sup_{u^*\in L^1(T)}\{\langle u,u^*\rangle_{\Lambda}-I_{f^*}(u^*)\}=I_f(u).
\end{multline*}
For any $u^*\in L^1(T)$, since $\langle\cdot, u^* \rangle_\Lambda\in D^*$, one has
\[
 (J)^*(u)\ge \langle u,u^*\rangle-J(\langle\cdot,u^* \rangle_{\Lambda})\ge \langle u,u^*\rangle-I_{f^*}(u^*),
\]
then
\[
 (J)^*(u)\ge \sup_{u^*\in L^1(T)}\{\langle u,u^*\rangle_{\Lambda}-I_{f^*}(u^*)\}=I_f(u).
\]
Therefore, $(J)^*|_{D}=I_f$.
In the following, we will show that $J$ is  lower semi-continuous on $D^*$: for any $\alpha\in \R$
\[
j_\alpha:=\{v\in D^*\bigg| J(v)\le \alpha\}
\]
is a weak* compact set in the topology $\sigma(D^*, D)$,
then one has
\[
J=J^{**}=(I_f)^*.
\]

 Let $D^0$ denote the orthogonal of $D$:
  \[
D^0:=\{v\in L^\infty(T)^*|v(u)=0,~\forall u\in D\},
\]
then the spaces $D^*$ and $L^\infty(T)^*/D^0$ are isomorphic. Let $\phi$ denote the isomorphism from $D^*$ to $L^\infty(T)^*/D^0$.
Moreover, the topology $\sigma(L^\infty(T)^*,L^\infty(T))/D^0$ is finer than the topology $\sigma(L^\infty(T)^*/D^0, D)$.  (See Proposition 35.5 and 35.6 in Treves~\cite{TVS}). Therefore, any  $\sigma(L^\infty(T)^*/D^0, D)-$open cover of $\phi(j_\alpha)$ is a $\sigma(L^\infty(T)^*,L^\infty(T))/D^0$-open cover where
\[
\phi(j_\alpha)=\left\{[v]\in L^{\infty}(T)^*/D^0\bigg|J(\phi^{-1}([v]))\le \alpha\right\}.
\]
If one considers $L^1(T)$ as a subspace of $L^\infty(T)^*$, then for any $u_1^*,u_2^*\in L^1(T)$, $u_1^*-u_2^*\in D_0\Leftrightarrow u_1^*-u_2^*\in D^0$ and therefore  $L^1(T)/D_0\subset L^\infty(T)^*/D^0$. Therefore,
\[
\phi(j_\alpha)\subset {\rm dom}(J\circ\phi^{-1})\subset L^1(T)/D_0.
\]
Note that $\phi(j_\alpha)=C_\alpha$, since for any $[v]\in L^1(T)/D_0$,
\[
J(\phi^{-1}([v]))=\inf_{u^*\in L^{1}(T), [u^*]=[v]}I_{f^*}(u^*).
\]
Any $\sigma(L^\infty(T)^*,L^\infty(T))/D^0-$open cover $\{U_i\}_{i\in I}$ of $\phi(j_\alpha)$, $\{U_i\cap L^1(T)\}_{i\in I}$ is a $\sigma(L^1(T), L^{\infty})/D_0-$open cover. Thanks to Lemma~\ref{compact}, there exists a finite open cover $\{U_i\cap L^1(T)\}_{i\in I_0}$. Then $\{U_i\}_{i\in I_0}$ is also a finite open cover in the topology  $\sigma(L^\infty(T)^*,L^\infty(T))/D^0$.  In conclusion, $j_\alpha$ is a weak* compact set in the topology $\sigma(D^*,D)$.

\end{proof}

\bibliographystyle{amsplain}
\bibliography{ref}
\bigskip

\end{document}